# Integral of Map into Abelian Ω-group

Aleks Kleyn

Aleks_Kleyn@MailAPS.org
http://AleksKleyn.dyndns-home.com:4080/
http://sites.google.com/site/AleksKleyn/
http://arxiv.org/a/kleyn_a_1
http://AleksKleyn.blogspot.com/

Abstract. The common in ring, module and algebra is that they are Abelian group with respect to addition. This property is enough to study integration. I treat integral of measurable map into normed Abelian $\Omega$-group. Theory of integration of maps into $\Omega$-group has a lot of common with theory of integration of functions of real variable. However I had to change some statements, since they implicitly assume either compactness of range or total order in $\Omega$-group.

# Contents

CHAPTER 1

# Preface

## 1.1. Preface to Version 1

The theory of measurable functions and integral developed in the middle of XX century solves a lot of mathematical problems. However this theory considers real functions.

When I initiated study of algebra with continuous basis, I realized that I need similar theory of integration in normed vector space. I considered integral of map into ring, module and algebra. The common in these algebraic structures is that they are Abelian group with respect to addition. This property is enough to study integration.

It was natural to consider the problem from the general point of view. I recalled that some time ago I faced similar task. So I decided to explore normed $\Omega$-group and the theory of integration of maps into $\Omega$-group.

There are a lot of references to similar definitions and theorems; so a reader sees that theory of integration of maps into $\Omega$-group has a lot of common with theory of integration of functions of real variable. However I had to change some statements, since they implicitly assume either compactness of range or total order in $\Omega$-group.

I have devoted much attention to the topology of normed $\Omega$-group. To determine the integral of measurable map $f$, I have to consider a sequence of simple maps uniformly convergent to the map $f$.

October, 2013

## 1.2. Preface to Version 2

When I was writing the text of version 1, I had to stop on Fubini's theorem. The proof of the theorem depends significantly on properties of real numbers. This is not acceptable in the case of an arbitrary $\Omega$-group. However, we cannot extract matrix of linear transformation of algebra with continuous basis if we do not use Fubini's theorem. So I returned to the paper to prove Fubini's theorem.

The section 4.3 is very important now for me. From theorems of the section 4.3, it follows that properties of integral are similar to properties of measure. It is attractive to consider more general definition of measure. The text of paper is ready to generalization. However the following problem should be resolved before. Let $\mu$ be measure with value in $\Omega$-group. Let $\mu(A) \ne 0$, $\mu(B) \ne 0$. We require that $\mu(A \cup B) \ne \emptyset$. We can make this request as part of definition.

Since I consider integral of map into $\Omega$-group $A$, I had to define interaction of real number and $A$-number. I consider the representation of real field in $\Omega$-group $A$. This consideration is essential when I study algebra with continuous basis.

March, 2014

## 1.3. Conventions

CONVENTION 1.3.1. *Element of $\Omega$-group $A$ is called $A$-**number***. *For instance, complex number is also called $C$-number, and quaternion is called $H$-number.* □

CONVENTION 1.3.2. *Let $A$ be $\Omega_1$-algebra. Let $B$ be $\Omega_2$-algebra. Notation*

$$A \longrightarrow_{*} B$$

*means that there is representation of $\Omega_1$-algebra $A$ in $\Omega_2$-algebra $B$.* □

Without a doubt, the reader may have questions, comments, objections. I will appreciate any response.

CHAPTER 2

# Measure

## 2.1. Algebra of Sets

DEFINITION 2.1.1. A nonempty system of sets $\mathcal{S}$ is called **semiring of sets**,[2.1] if

2.1.1.1: $\emptyset \in \mathcal{S}$
2.1.1.2: If $A, B \in \mathcal{S}$, then $A \cap B \in \mathcal{S}$
2.1.1.3: If $A, A_1 \in \mathcal{S}, A_1 \subset A$, then the set $A$ can be represented as

$$(2.1.1) \qquad A = \bigcup_{i=1}^{n} A_i \quad A_i \in \mathcal{S}$$

where $i \neq j => A_i \cap A_j = \emptyset$

The representation (2.1.1) of the set $A$ is called **finite expansion of set** $A$. □

DEFINITION 2.1.2. A nonempty system of sets $\mathcal{R}$ is called **ring of sets**,[2.2] if condition $A, B \in \mathcal{R}$ imply $A \triangle B, A \cap B \in \mathcal{R}$. A set $E \in \mathcal{R}$ is called **unit of ring of sets** if

$$A \cap E = A$$

A ring of sets with unit is called **algebra of sets**. □

REMARK 2.1.3. For any $A$, $B$

$$A \cup B = (A \Delta B) \Delta (A \cap B)$$
$$A \setminus B = A \Delta (A \cap B)$$

Therefore, if $A, B \in \mathcal{R}$, then $A \cup B \in \mathcal{R}, A \setminus B \in \mathcal{R}$. □

THEOREM 2.1.4. *Ring of sets $\mathcal{R}$ is semiring.*

PROOF. Let $A, A_1 \in \mathcal{R}, A_1 \subset A$. Then

$$A = A_1 \cup A_2$$

where

$$A_2 = A \setminus A_1 \in \mathcal{R}$$

□

THEOREM 2.1.5. *The intersection $\mathcal{R} = \bigcap \mathcal{R}_i$ of any set of rings is ring.*[2.3]

PROOF. Let $A, B \in \mathcal{R}$. Then for any $i$, $A, B \in \mathcal{R}_i$. According to the definition 2.1.2, for any $i$, $A \triangle B, A \cap B \in \mathcal{R}_i$. Therefore, $A \triangle B, A \cap B \in \mathcal{R}$. According to the definition 2.1.2, the set $\mathcal{R}$ is a ring of sets. □

---

[2.1]See also the definition [1]-2, page 32.
[2.2]See also the definition [1]-1, page 31.
[2.3]See also the theorem [1]-1 on page 32.

THEOREM 2.1.6. *For any nonempty system of sets $\mathcal{C}$, there exists a unique ring of sets[2.4] $\mathcal{R}(\mathcal{C})$ containing $\mathcal{C}$ and contained in any ring $\mathcal{R}$ such that $\mathcal{C} \subset \mathcal{R}$.*

PROOF. Let rings of sets $\mathcal{R}_1$, $\mathcal{R}_2$, $\mathcal{R}_1 \neq \mathcal{R}_2$, satisfy to condition of the theorem. By the theorem 2.1.5, the set $\mathcal{R}_1 \cap \mathcal{R}_2$, is the ring of sets satisfying to condition of the theorem. Therefore, if ring of sets $\mathcal{R}(\mathcal{C})$ exists, then it is unique.

Let
$$R = \bigcup_{X \in \mathcal{C}} X$$
The set $\mathcal{B}(R)$ of all subsets of the set $R$ is ring of sets and $\mathcal{C} \subseteq \mathcal{B}(R)$. Let $\Sigma$ be set of such rings of sets $\mathcal{R}$ that $\mathcal{C} \subseteq \mathcal{R} \subseteq \mathcal{B}(R)$. Then, according to the theorem 2.1.5, the set $\mathcal{P} = \bigcap_{\mathcal{R} \in \Sigma} \mathcal{R}$ is ring of sets which satisfies the theorem. □

THEOREM 2.1.7. *Let $\mathcal{C}$ be nonempty system of sets. Let*

(2.1.2) $$R = \bigcup_{X \in \mathcal{C}} X$$

*If*

(2.1.3) $$R \in \mathcal{C}$$

*then the ring of sets $\mathcal{R}(\mathcal{C})$ is algebra of sets.*

PROOF. The set $\mathcal{B}(R)$ of all subsets of the set $R$ is ring of sets and $\mathcal{C} \subseteq \mathcal{B}(R)$. Let $\Sigma$ be set of such rings of sets $\mathcal{R}$ that

(2.1.4) $$\mathcal{C} \subseteq \mathcal{R} \subseteq \mathcal{B}(R)$$

From (2.1.3), (2.1.4), it follows that $R \in \mathcal{R}$ for any $\mathcal{R} \in \Sigma$. Then, according to the theorem 2.1.5, the set $\mathcal{P} = \bigcap_{\mathcal{R} \in \Sigma} \mathcal{R}$ is the smallest ring of sets such that $R \in \mathcal{P}$, $\mathcal{C} \subseteq \mathcal{P}$. From (2.1.2), it follows that $R \cap A = A$ for any set $A \in \mathcal{P}$. According to the definition 2.1.2, the ring of sets $\mathcal{P}$ is algebra of sets. □

THEOREM 2.1.8. *For any nonempty system of sets $\mathcal{C}$, there exists a unique algebra of sets $\mathcal{A}(\mathcal{C})$ containing $\mathcal{C}$ and contained in any algebra $\mathcal{R}$ such that $\mathcal{C} \subset \mathcal{R}$.*

PROOF. Let
$$R = \bigcup_{X \in \mathcal{C}} X$$
The theorem follows from theorems 2.1.6, 2.1.7, if we assume
$$\mathcal{A}(\mathcal{C}) = \mathcal{R}(\{R\} \cup \mathcal{C})$$
□

LEMMA 2.1.9. *Let[2.5] $\mathcal{S}$ be semiring. Let $A$, $A_1$, ..., $A_n \in \mathcal{S}$, $i \neq j \Rightarrow A_i \cap A_j = \emptyset$. Then there exists finite expansion of set $A$*
$$A = \bigcup_{i=1}^{s} A_i \quad s \geq n$$

---

[2.4]See also the theorem [1]-2 on page 32.

[2.5]Lemmas 2.1.9, 2.1.10 and the theorem 2.1.11 are similar to lemmas 1, 2 and the theorem 3, [1], pages 33, 34.

PROOF. For $n = 1$, the lemma follows from the statement 2.1.1.3.

Let the lemma hold for $n = m$. Let sets $A_1, ..., A_{m+1}$ satisfy condition of the lemma. According to the assumption

$$(2.1.5) \qquad A = A_1 \cup ... \cup A_m \cup B_1 \cup ... \cup B_p$$

where $B_i \in \mathcal{S}$, $i = 1, ..., p$, $i \neq j \Rightarrow A_i \cap A_j = \emptyset$, $A_i \cap B_j = \emptyset$, $i \neq j \Rightarrow B_i \cap B_j = \emptyset$. According to the statement 2.1.1.2

$$(2.1.6) \qquad B_{i \cdot 1} = A_{m+1} \cap B_i \in \mathcal{S}$$

According to the statement 2.1.1.3

$$(2.1.7) \qquad B_i = B_{i \cdot 1} \cup ... \cup B_{i \cdot r_s} \quad B_{i \cdot j} \in \mathcal{S}$$

Since $A_{m+1} \subset B_1 \cup ... \cup B_p$, then

$$A = A_1 \cup ... \cup A_m \cup A_{m+1} \cup \bigcup_{i=1}^{p} \bigcup_{j=2}^{r_i} B_{i \cdot j}$$

follows from (2.1.5), (2.1.6), (2.1.7). Therefore, the lemma holds for $n = m + 1$.

According to mathematical induction, the lemma holds for any $n$. □

LEMMA 2.1.10. *For any finite system of sets* $A_1, ..., A_n \in \mathcal{S}$ *there exists finite system of sets* $B_1, ..., B_t \in \mathcal{S}$, $i \neq j \Rightarrow B_i \cap B_j = \emptyset$, *such that*

$$(2.1.8) \qquad A_i = \bigcup_{j \in M_i} B_j$$

*where* $M_i \subset \{1, ..., t\}$.

PROOF. For $n = 1$, the lemma is evident, since we let $t = 1$, $B_1 = A_1$.

Let the lemma hold for $n = m$. Consider the system of sets $A_1, ..., A_{m+1}$. Let $B_i \in \mathcal{S}$, $i = 1, ..., p$, be sets satisfying the condition of the lemma with respect to $A_1, ..., A_m$. According to the statement 2.1.1.2

$$(2.1.9) \qquad B_{i \cdot 1} = A_{m+1} \cap B_i \in \mathcal{S}$$

According to the statement 2.1.1.3

$$(2.1.10) \qquad B_i = B_{i \cdot 1} \cup ... \cup B_{i \cdot r_s} \quad B_{i \cdot j} \in \mathcal{S}$$

The expansion

$$(2.1.11) \qquad A_{m+1} = \bigcup_{s=1}^{t} B_{s \cdot 1} \cup \bigcup_{p=1}^{q} B'_p \quad B'_p \in \mathcal{S}$$

follows from the lemma 2.1.9. The expansion

$$A_i = \bigcup_{j \in M_i} \bigcup_{p=1}^{r_j} B_{j \cdot p}$$

follows from equations (2.1.8), (2.1.10). From equations (2.1.9), (2.1.11), it follows that $B_i \cap B'_p = \emptyset$. Therefore, from the equation (2.1.10), it follows that $B_{i \cdot j} \cap B'_p = \emptyset$. Therefore, sets $B_{i \cdot j}$, $B'_p$ satisfy the conditions of the lemma with respect to $A_1, ..., A_{m+1}$. Therefore, the lemma holds for $n = m + 1$.

According to mathematical induction, the lemma holds for any $n$. □

THEOREM 2.1.11. *Let $\mathcal{S}$ be semiring of sets. The system $\mathcal{R}$ of sets $A$ which have finite expansion*

$$A = \bigcup_{i=1}^{n} A_i \quad A_i \in \mathcal{S}$$

*is* **ring of sets generated by the semiring of sets $\mathcal{S}$.**

PROOF. Let $A, B \in \mathcal{R}$ Then

(2.1.12) $$A = \bigcup_{i=1}^{n} A_i \quad A_i \in \mathcal{S}$$

(2.1.13) $$B = \bigcup_{i=1}^{m} B_i \quad B_i \in \mathcal{S}$$

From (2.1.12), (2.1.13) and from the statement 2.1.1.2, it follows that

$$C_{ij} = A_i \cap B_j \in \mathcal{S}$$

According to the lemma 2.1.9,

(2.1.14)
$$A_i = \bigcup_{j=1}^{m} C_{ij} \cup \bigcup_{k=1}^{r_i} D_{ik} \quad D_{ik} \in \mathcal{S}$$
$$B_j = \bigcup_{i=1}^{n} C_{ij} \cup \bigcup_{k=1}^{s_j} E_{kj} \quad E_{kj} \in \mathcal{S}$$

From (2.1.14), it follows that

$$A \cup B = \bigcup_{ij} C_{ij} \in \mathcal{R}$$
$$A \Delta B = \bigcup_{ik} D_{ik} \cup \bigcup_{jl} E_{jl} \in \mathcal{R}$$

Therefore, $\mathcal{R}$ is ring of sets. $\square$

THEOREM 2.1.12. *Let semiring of sets $\mathcal{C}$ contain unit. Then the ring of sets $\mathcal{R}(\mathcal{C})$ is an algebra of sets.*

PROOF. The theorem follows from theorems 2.1.8, 2.1.11. $\square$

DEFINITION 2.1.13. *The ring of sets $\mathcal{R}$ is called $\sigma$-**ring of sets**,[2.6] if condition $A_i \in \mathcal{R}$, $i = 1, ..., n, ...,$ imply*

$$\bigcup_{n} A_n \in \mathcal{R}$$

*$\sigma$-Ring of sets with unit is called $\sigma$-**algebra of sets**.* $\square$

THEOREM 2.1.14. *The intersection $\mathcal{R} = \bigcap \mathcal{R}_i$ of any set of $\sigma$-rings is $\sigma$-ring.[2.7]*

---

[2.6]See similar definition in [1], p. 35, definition 3.
[2.7]See also the theorem [1]-1 on page 32.

PROOF. Since $A_i \in \mathcal{R}$, $i = 1, ..., n, ...,$ then for any $i, j$, $A_i \in \mathcal{R}_j$. Therefore, for any $j$,
$$\bigcup_i A_i \in \mathcal{R}_j$$
Therefore,
$$\bigcup_i A_i \in \mathcal{R}$$
According to the definition 2.1.13, the set $\mathcal{R}$ is a $\sigma$-ring of sets. $\square$

THEOREM 2.1.15. *For any nonempty system of sets $\mathcal{C}$, there exists a unique $\sigma$-ring of sets[2.8] $\mathcal{R}_\sigma(\mathcal{C})$ containing $\mathcal{C}$ and contained in any $\sigma$-ring $\mathcal{R}$ such that $\mathcal{C} \subset \mathcal{R}$.*

PROOF. Let $\sigma$-rings of sets $\mathcal{R}_1, \mathcal{R}_2, \mathcal{R}_1 \neq \mathcal{R}_2$, satisfy to condition of the theorem. By the theorem 2.1.14, the set $\mathcal{R}_1 \cap \mathcal{R}_2$, is the $\sigma$-ring of sets satisfying to condition of the theorem. Therefore, if $\sigma$-ring of sets $\mathcal{R}_\sigma(\mathcal{C})$ exists, then it is unique.

Let
$$R = \bigcup_{X \in \mathcal{C}} X$$
The set $\mathcal{B}(R)$ of all subsets of the set $R$ is $\sigma$-ring of sets and $\mathcal{C} \subseteq \mathcal{B}(R)$. Let $\Sigma$ be set of such $\sigma$-rings of sets $\mathcal{R}$ that $\mathcal{C} \subseteq \mathcal{R} \subseteq \mathcal{B}(R)$. Then, according to the theorem 2.1.5, the set $\mathcal{P} = \bigcap_{\mathcal{R} \in \Sigma} \mathcal{R}$ is $\sigma$-ring of sets which satisfies the theorem. $\square$

THEOREM 2.1.16. *For any nonempty system of sets $\mathcal{C}$, there exists a unique $\sigma$-algebra of sets $\mathcal{A}_\sigma(\mathcal{C})$ containing $\mathcal{C}$ and contained in any $\sigma$-algebra $\mathcal{R}$ such that $\mathcal{C} \subset \mathcal{R}$.*

PROOF. Let
$$R = \bigcup_{X \in \mathcal{C}} X$$
The theorem follows from theorems 2.1.15, 2.1.7, if we assume
$$\mathcal{A}_\sigma(\mathcal{C}) = \mathcal{R}_\sigma(\{R\} \cup \mathcal{C})$$
$\square$

## 2.2. Measure

DEFINITION 2.2.1. Let $\mathcal{C}_m$ be semiring of sets.[2.9] The map
$$m : \mathcal{C}_m \to R$$
is called **measure**, if

2.2.1.1: $m(A) \geq 0$

---
[2.8]See also the theorem [1]-2 on page 32.
[2.9]See also the definition [1]-1 on page 270.

2.2.1.2: The map $m$ is additive map. If a set $A \in \mathcal{C}_m$ has finite expansion

$$A = \bigcup_{i=1}^{n} A_i \quad A_i \in \mathcal{C}_m$$

where $i \neq j => A_i \cap A_j = \emptyset$, then

$$m(A) = \sum_{i=1}^{n} m(A_i)$$

$\square$

THEOREM 2.2.2. $m(\emptyset) = 0$.

PROOF. Since $\emptyset = \emptyset \cup \emptyset$, $\emptyset \cap \emptyset = \emptyset$, then the theorem follows from the statements 2.1.1.1, 2.2.1.2. $\square$

DEFINITION 2.2.3. A measure $\mu$ is called[2.10] **extension of measure** $m$, if $\mathcal{C}_m \subseteq \mathcal{C}_\mu$ and $\mu(A) = m(A)$, $A \in \mathcal{C}_m$. $\square$

THEOREM 2.2.4. Let[2.11] $\mathcal{R}(\mathcal{C}_m)$ be ring of sets generated by the semiring of sets $\mathcal{C}_m$. Measure $m$ defined on the semiring of sets $\mathcal{C}_m$ has a unique extension $\mu$ defined on the ring of sets $\mathcal{R}(\mathcal{C}_m)$.

PROOF. By the theorem 2.1.11, every set $A \in \mathcal{R}(\mathcal{C}_m)$ has a finite expansion

$$(2.2.1) \qquad A = \bigcup_{i=1}^{n} A_i \quad A_i \in \mathcal{C}_m$$

where $i \neq j => A_i \cap A_j = \emptyset$. Let

$$(2.2.2) \qquad \mu(A) = \sum_{i=1}^{n} m(A_i)$$

Let $A$ have two finite expansions

$$A = \bigcup_{i=1}^{n} A_i = \bigcup_{j=1}^{m} B_j \quad A_i, B_j \in \mathcal{C}_m$$

Since $A_i \cap B_j \in \mathcal{C}_m$ according to the statement 2.1.1.2, then

$$\sum_{i=1}^{n} m(A_i) = \sum_{i=1}^{n} \sum_{j=1}^{m} m(A_i \cap B_j) = \sum_{j=1}^{m} m(B_j)$$

follows from the statement 2.2.1.2. Therefore, the value $\mu(A)$ defined by the equation (2.2.2) is independent of the finite expansion (2.2.1).

Therefore, we designed the map

$$\mu : \mathcal{R}(\mathcal{C}_m) \to R$$

which satisfies the definition 2.2.1.

---

[2.10] See also the definition [1]-2 on page 271.
[2.11] See also the theorem [1]-3 on page 271.

For any extension $\mu'$ of measure $m$ and for the finite expansion (2.2.1), from the statement 2.2.1.2 and the definition 2.2.3, it follows that

$$\mu'(A) = \sum_{i=1}^{n} \mu'(A_i) = \sum_{i=1}^{n} \mu(A_i) = \mu(A)$$

Therefore, the measure $\mu'$ coinside with the measure $\mu$ defined by the the equation (2.2.2). □

DEFINITION 2.2.5. A measure $\mu$ is called **complete measure**, if conditions $B \subset A$, $\mu(A) = 0$ imply that $B$ is measurable set. □

THEOREM 2.2.6. *Let $\mu$ be a measure defined on the set $X$. Let*[2.12]

(2.2.3) $$A \subset B \quad A, B \in \mathcal{C}_\mu$$

*Then*

(2.2.4) $$\mu(A) \leq \mu(B)$$

PROOF. The equation

(2.2.5) $$B = A \cup (B \setminus A)$$

follows from the statement (2.2.3). According to the remark 2.1.3, from the statement (2.2.3), it follows that

$$B \setminus A \in \mathcal{C}_\mu$$

The equation

(2.2.6) $$\mu(B) = \mu(A) + \mu(B \setminus A)$$

follows from the equation (2.2.5). According to the statement 2.2.10.1,

(2.2.7) $$\mu(B \setminus A) \geq 0$$

The statement (2.2.4) follows from the statements 2.2.10.2, (2.2.6). □

THEOREM 2.2.7. *Let $\mu$ be a measure defined on the set $X$. Let*

$$A \subset \bigcup_i A_i \quad A, A_i \in \mathcal{C}_\mu$$

*where $\{A_i\}$ is finite or countable system of sets. Then*

$$\mu(A) \leq \mu(B)$$

PROOF. We prove the theorem in case of $\{2\}$ sets. It is easy to prove general case using proof by induction. Since

$$A \subseteq A_1 \cup A_2$$

then the theorem follows from the theorem 2.2.6, from the statement

$$\mu(A_1 \cup A_2) = \mu(A_1) + \mu(A_2) - \mu(A_1 \cap A_2) \leq \mu(A_1) + \mu(A_2)$$

and from the statement 2.2.10.1. □

THEOREM 2.2.8. *Let $\mu$ be a complete measure defined on the set $X$. Let*

$$A \subset B \quad A, B \in \mathcal{C}_\mu$$

*Since $\mu(B) = 0$, then $\mu(A) = 0$.*

---
[2.12]See also the theorem [1]-2, p. 257.

PROOF. The theorem follows from the definition 2.2.5 and from the theorem 2.2.6. □

REMARK 2.2.9. We assume that considered measure is complete measure. □

DEFINITION 2.2.10. Let $\mathcal{C}_\mu$ be $\sigma$-algebra of sets of set $F$.[2.13] The map

$$\mu : \mathcal{C}_\mu \to R$$

into real field $R$ is called $\sigma$-**additive measure**, if, for any set $X \in \mathcal{C}_\mu$, following conditions are true.

2.2.10.1: $\mu(X) \geq 0$
2.2.10.2: Let

$$X = \bigcup_i X_i \quad i \neq j => X_i \cap X_j = \emptyset$$

be finite or countable union of sets $X_n \in \mathcal{C}_\mu$. Then

$$\mu(X) = \sum_i \mu(X_i)$$

where series on the right converges absolutly.

□

THEOREM 2.2.11. *Let $m$ be $\sigma$-additive measure defined on the semiring of sets $\mathcal{C}_m$. Extension $\mu$ of measure $m$ defined on the ring of sets $\mathcal{R}(\mathcal{C}_m)$ is $\sigma$-additive measure.*

PROOF. Let $A, B_n \in \mathcal{R}(\mathcal{C}_m)$, $n = 1, 2, ...$, $i \neq j => B_i \cap B_j = \emptyset$. Let

$$A = \bigcup_{n=1}^\infty B_n$$

According to the theorem 2.1.11, there exist finite expansions

(2.2.8) $$A = \bigcup_j A_j \quad B_n = \bigcup_j B_{nj}$$

where

$$A_k \cap A_l = \emptyset \quad B_{nk} \cap B_{nl} = \emptyset \quad k \neq l$$

Let $C_{nil} = B_{ni} \cap A_l$. It is evident that sets $C_{nil}$ are pairwise disjoint and

(2.2.9) $$A_j = \bigcup_{n=1}^\infty \bigcup_i C_{nij} \quad B_{ni} = \bigcup_j C_{nij}$$

From the equation (2.2.9) and the statement 2.2.10.2, it follows that

(2.2.10) $$m(A_j) = \sum_{n=1}^\infty \sum_i m(C_{nij})$$
$$m(B_{ni}) = \sum_j m(C_{nij})$$

---

[2.13]See similar definitions in [1], definition 1 on page 270 and definition 2 on page 272.

From the equation (2.2.8) and the theorem 2.2.4, it follows that

$$\mu(A) = \sum_j m(A_j)$$
(2.2.11)
$$\mu(B_n) = \sum_i m(B_{ni})$$

The equation

$$\mu(A) = \sum_{n=1}^{\infty} \mu(B_n)$$

follows from (2.2.10), (2.2.11). □

THEOREM 2.2.12 (Continuity of $\sigma$-additive measure). *Let*[2.14]

(2.2.12) $$A_1 \supset A_2 \supset ...$$

*be a sequence of $\mu$-measurable sets. Then*

$$\mu(A) = \lim_{n \to \infty} \mu(A_n)$$

*where* $A = \bigcap_n A_n$ .

PROOF. We consider the case when $A = \emptyset$. General case reduces to this case if we replace $A_n$ by $A_n \setminus A$. According to the statement (2.2.12)

(2.2.13) $$A_n = (A_n \setminus A_{n+1}) \cup (A_{n+1} \setminus A_{n+2}) \cup ... \quad n = 1, ...$$

According to the statement 2.2.10.2, the equation

(2.2.14) $$\mu(A_n) = \sum_{k=n}^{\infty} \mu(A_k \setminus A_{k+1}) \quad n = 1, ...$$

follows from (2.2.13). Since the series (2.2.14) for $n = 1$ converges, its remainder (2.2.14) converges to 0 when $n \to \infty$. Therefore,

$$\lim_{n \to \infty} \mu(A_n) = 0$$

□

THEOREM 2.2.13 (Continuity of $\sigma$-additive measure). *Let*[2.15]

$$A_1 \subset A_2 \subset ...$$

*be a sequence of $\mu$-measurable sets. Then*

$$\mu(A) = \lim_{n \to \infty} \mu(A_n)$$

*where* $A = \bigcup_n A_n$ .

PROOF. The theorem follows from the theorem 2.2.12, if we consider sets $X \setminus A_n$. □

---

[2.14] See also the theorem [1]-11 on page 266.
[2.15] See also corollary of the theorem [1]-11 on pages 266, 267.

## 2.3. Lebesgue Extension of Measure

DEFINITION 2.3.1. Let $m$ be $\sigma$-additive measure on a semiring $\mathcal{S}$ with unit $E$.[2.16] The **outer measure** of a set $A \subseteq E$ is defined by the equation
$$\mu^*(A) = \inf_{A \subseteq \bigcup_k B_k} \sum_k m(B_k)$$
where the lower bound is taken over all coverings of $A$ by a finite or countable system of sets $B_n \in \mathcal{S}$. □

THEOREM 2.3.2 (countable subadditivity). If[2.17]
$$(2.3.1) \qquad A \subseteq \bigcup_n A_n$$
where $\{A_n\}$ is a finite or countable system of sets, then
$$(2.3.2) \qquad \mu^*(A) \leq \sum_n \mu^*(A_n)$$

PROOF. By the definition 2.3.1, for any $A_n$, there exists finite or countable system of sets $\{P_{nk}\}$ such that
$$(2.3.3) \qquad A_n \subseteq \bigcup_k P_{nk}$$
$$(2.3.4) \qquad \sum_k m(P_{nk}) \leq \mu^*(A_n) + \frac{\epsilon}{2^n}$$
From (2.3.1), (2.3.3), it follows that
$$(2.3.5) \qquad A \subseteq \bigcup_n \bigcup_k P_{nk}$$
According to the statement (2.3.5) and to the definition 2.3.1, inequation
$$(2.3.6) \qquad \mu^*(A) \leq \sum_n \sum_k m(P_{nk}) \leq \sum_n \mu^*(A_n) + \epsilon$$
follows from the inequation (2.3.4). (2.3.2) follows from (2.3.6), since $\epsilon$ is arbitrary. □

THEOREM 2.3.3. *For any sets $A$, $B$,*
$$(2.3.7) \qquad |\mu^*(A) - \mu^*(B)| \leq \mu^*(A \Delta B)$$

PROOF. Since $\mu^*(A) \geq \mu^*(B)$, then the statement
$$(2.3.8) \qquad \mu^*(A) \leq \mu^*(B) + \mu(A \Delta B)$$
folows from the statement
$$A \subseteq B \cup (A \Delta B)$$
and the theorem 2.3.2. Since $\mu^*(B) \geq \mu^*(A)$, then the statement
$$(2.3.9) \qquad \mu^*(B) \leq \mu^*(A) + \mu(A \Delta B)$$

---
[2.16] See definition [1]-1, page 276.
[2.17] See also theorems [1]-4 on page 259, [1]-1 on page 276.

folows from the statement
$$B \subseteq A \cup (A\Delta B)$$
and the theorem 2.3.2. The statement (2.3.7) folows from statements (2.3.8), (2.3.9). □

DEFINITION 2.3.4. The set $A$ is called **Lebesgue measurable**, if for any $\epsilon \in R$, $\epsilon > 0$, there exists $B \in \mathcal{R}(\mathcal{S})$ such that[2.18]
$$\mu^*(A\Delta B) < \epsilon$$
□

Let $\mathcal{C}_\mu$ be the system of Lebesgue measurable sets.

THEOREM 2.3.5. *Let $m$ be $\sigma$-additive measure on a semiring $\mathcal{S}$ with unit $E$.[2.19] If a set $A$ is Lebesgue measurable, then the set $E \setminus A$ is also Lebesgue measurable.*

PROOF. Let $A$ be Lebesgue measurable set. By the definition 2.3.4, for any $\epsilon \in R$, $\epsilon > 0$, there exists $B \in \mathcal{R}(\mathcal{S})$ such that
(2.3.10) $$\mu^*(A\Delta B) < \epsilon$$
By the remark 2.1.3, $E \setminus B \in \mathcal{R}(\mathcal{S})$. The statement
$$\mu^*((E \setminus A)\Delta(E \setminus B)) < \epsilon$$
follows from the statement (2.3.10) and from the equation
$$A\Delta B = (E \setminus A)\Delta(E \setminus B)$$
□

THEOREM 2.3.6. *Let $A_1$, $A_2 \in \mathcal{C}_\mu$. Then $A = A_1 \setminus A_2 \in \mathcal{C}_\mu$.*

PROOF. By the definition 2.3.4, for any $\epsilon \in R$, $\epsilon > 0$, there exist $B_1$, $B_2 \in \mathcal{R}(\mathcal{S})$ such that
(2.3.11) $$\mu^*(A_1\Delta B_1) < \frac{\epsilon}{2}$$
$$\mu^*(A_2\Delta B_2) < \frac{\epsilon}{2}$$
By the remark 2.1.3, $B = B_1 \setminus B_2 \in \mathcal{R}(\mathcal{S})$. The statement
$$\mu^*(A\Delta B) < \epsilon$$
follows from the statement (2.3.11), from the statement
$$(A_1 \setminus A_2)\Delta(B_1 \setminus B_2) \subseteq (A_1\Delta B_1) \cup (A_2\Delta B_2)$$
and from the theorem 2.3.2. □

THEOREM 2.3.7. *Let $m$ be $\sigma$-additive measure on a semiring $\mathcal{S}$ with unit $E$.[2.20] Let $\mu$ be extension of measure $m$ on the ring of sets $\mathcal{R}(\mathcal{S})$. Evry set $A \in \mathcal{R}(\mathcal{S})$ is Lebesgue measurable and*
(2.3.12) $$\mu^*(A) = \mu(A)$$

---

[2.18] See theorem [1]-3, page 277.
[2.19] See also remark in [1] after the definition 4 on page 276.
[2.20] See also theorems [1]-5 on page 259, [1]-2 on page 277.

PROOF. Let $A \in \mathcal{R}(\mathcal{S})$. According to the definition 2.1.1 and the theorem 2.1.11, the set $A$ can be represented as

$$A = \bigcup_{i=1}^{n} A_i \quad A_i \in \mathcal{S}$$

where $i \neq j => A_i \cap A_j = \emptyset$. According to the theorem 2.2.4,

$$\mu(A) = \sum_{i=1}^{n} m(A_i) \tag{2.3.13}$$

According to the definition 2.3.1, since sets $A_i$ cover set $A$, then

$$\mu^*(A) \leq \sum_{i=1}^{n} m(A_i) = \mu(A) \tag{2.3.14}$$

follows from (2.3.13).

Let $\{Q_i\}$, $Q_i \in \mathcal{S}$, be finite or countable system of sets covering the set $A$. According to the theorems 2.2.4, 2.2.7,

$$\mu(A) \leq \sum_{j} m(Q_j) \tag{2.3.15}$$

According to the definition 2.3.1,

$$\mu(A) \leq \mu^*(A) \tag{2.3.16}$$

follows from (2.3.15).

(2.3.12) follows from (2.3.14), (2.3.16). □

THEOREM 2.3.8. *The system $\mathcal{C}_\mu$ of Lebesgue measurable sets is algebra of sets.*

PROOF. From the theorem 2.3.6, from the definition 2.1.2 and from equations

$$A_1 \cap A_2 = A_1 \setminus (A_1 \setminus A_2)$$
$$A_1 \cup A_2 = E \setminus ((E \setminus A_1) \cap (E \setminus A_2))$$
$$A_1 \Delta A_2 = (A_1 \cup A_2) \setminus (A_1 \cap A_2)$$

it follows that $\mathcal{C}_\mu$ is ring of sets. The ring of sets $\mathcal{C}_\mu$ is algebra of sets because $E \in \mathcal{C}_\mu$ is unit of the ring of sets $\mathcal{C}_\mu$. □

THEOREM 2.3.9. *The map $\mu^*(A)$ is additive on the algebra of sets $\mathcal{C}_\mu$.*[2.21]

PROOF. It is enough to prove the theorem for two sets. Let $A_1, A_2 \in \mathcal{C}_\mu$. By the definition 2.3.4, for any $\epsilon \in R$, $\epsilon > 0$, there exist $B_1, B_2 \in \mathcal{R}(\mathcal{S})$ such that

$$\mu^*(A_1 \Delta B_1) < \frac{\epsilon}{2}$$
$$\mu^*(A_2 \Delta B_2) < \frac{\epsilon}{2} \tag{2.3.17}$$

Let $A = A_1 \cup A_2$, $B = B_1 \cup B_2$. According to the theorem 2.3.8, $A \in \mathcal{C}_\mu$. Since $A_1 \cap A_2 = \emptyset$, then

$$B_1 \cap B_2 \subseteq (A_1 \Delta B_1) \cup (A_2 \Delta B_2) \tag{2.3.18}$$

The statement

$$\mu(B_1 \cap B_2) \leq \epsilon \tag{2.3.19}$$

---

[2.21]See also theorems [1]-5 on page 259, [1]-2 on page 277.

follows from the statement (2.3.18) and from the theorems 2.3.2, 2.3.7. The statement

$$|\mu(B_1) - \mu^*(A_1)| \leq \frac{\epsilon}{2}$$
$$|\mu(B_2) - \mu^*(A_2)| \leq \frac{\epsilon}{2}$$
(2.3.20)

follows from the statement (2.3.17) and from the theorems 2.3.3, 2.3.7. Since measure is additive on algebra of sets $\mathcal{R}(\mathcal{S})$, then the statement

(2.3.21) $\quad \mu(B) = \mu(B_1) + \mu(B_2) - \mu(B_1 \cap B_2) \geq \mu^*(A_1) + \mu^*(A_2) - 2\epsilon$

follows from statements (2.3.19), (2.3.20). The statement

(2.3.22) $\quad \mu^*(A) \geq \mu(B) - \mu^*(A \Delta B) \geq \mu(B) - 2\epsilon \geq \mu^*(A_1) + \mu^*(A_2) - 3\epsilon$

follows from the statement

$$A \Delta B \subseteq (A_1 \Delta B_1) \cup (A_2 \Delta B_2)$$

The statement

(2.3.23) $\quad \mu^*(A_1) + \mu^*(A_2) \geq \mu^*(A) \geq \mu^*(A_1) + \mu^*(A_2) - 3\epsilon$

follows from the statement (2.3.22) and from the theorem 2.3.2. Since $\epsilon$ can be made arbitrary small, then

$$\mu^*(A) = \mu^*(A_1) + \mu^*(A_2)$$

follows from the statement (2.3.23). $\square$

DEFINITION 2.3.10. If a set $A$ is Lebesgue measurable,[2.22] then the value $\mu(A) = \mu^*(A)$ is called **Lebesgue measure**. The map $\mu$ defined on the algebra of sets $\mathcal{C}_\mu$ is called **Lebesgue extension of measure** $m$. $\square$

THEOREM 2.3.11. *The algebra of sets $\mathcal{C}_\mu$ is $\sigma$-algebra[2.23] with unit $E$.*

PROOF. Let $A_1, ...,$ be countable system of Lebesgue measurable sets. Let

(2.3.24) $$A = \bigcup_{n=1}^{\infty} A_n$$

Let

(2.3.25) $$A'_n = A_n \setminus \bigcup_{k=1}^{n-1} A_k$$

From (2.3.24), (2.3.25), it folows that

(2.3.26) $$A = \bigcup_{n=1}^{\infty} A'_n$$

where $i \neq j \implies A'_i \cap A'_j = \emptyset$. According to the theorem 2.3.8, the definition 2.1.2 and remark 2.1.3, all sets $A'_n$ are Lebesgue measurable sets. According to the theorem 2.3.8 and to the definition 2.3.1, for any $n$

$$\sum_{k=1}^{n} \mu(A'_k) = \mu\left(\bigcup_{k=1}^{n} A'_k\right) \leq \mu^*(A)$$

---

[2.22]See definition [1]-4, page 276.
[2.23]See also theorems [1]-9 on page 264, [1]-7 on page 277.

Therefore, the series $\sum_{n=1}^{\infty} \mu(A'_n)$ convergies. Therefore, for any $\epsilon \in R$, $\epsilon > 0$, there exists $N$ such that

(2.3.27) $$\sum_{n>N} \mu(A'_n) < \frac{\epsilon}{2}$$

According to the theorem 2.3.8,

(2.3.28) $$C = \sum_{n=1}^{N} \mu(A'_n) \in \mathcal{C}_\mu$$

According to the definition 2.3.4, from the statement (2.3.28), it follows that there exists $B \in \mathcal{R}(\mathcal{S})$ such that

(2.3.29) $$\mu^*(C \Delta B) < \frac{\epsilon}{2}$$

Since
$$A \Delta B \subseteq (C \Delta B) \cup \left( \bigcup_{n>N} A'_n \right)$$

then, from (2.3.27), (2.3.29) and from the theorem 2.3.2, it follows that
$$\mu^*(A \Delta B) < \epsilon$$

According to the definition 2.3.4, $A \in \mathcal{C}_\mu$.

Since
$$\bigcap_n A_n = E \setminus \bigcup_n (E \setminus A_n)$$

then the theorem follows from the theorem 2.3.5. □

THEOREM 2.3.12. *The map $\mu(A)$ is $\sigma$-additive[2.24] on the algebra of sets $\mathcal{C}_\mu$.*

PROOF. Let
$$A = \bigcup_{i=1}^{\infty} A_i \quad A_i \in \mathcal{C}_\mu$$

where $i \neq j => A_i \cap A_j = \emptyset$. By the theorem 2.3.11, $A \in \mathcal{C}_\mu$. By the theorem 2.3.2,

(2.3.30) $$\mu(A) \leq \sum_i \mu(A_i)$$

By the theorem 2.3.9, for any $N$
$$\mu(A) \geq \mu \left( \bigcup_{i=1}^{N} A_i \right) = \sum_{i=1}^{N} \mu(A_i)$$

and therefore

(2.3.31) $$\mu(A) \geq \sum_i \mu(A_i)$$

The equation
$$\mu(A) = \sum_i \mu(A_i)$$

follows from (2.3.30), (2.3.31). □

---

[2.24]See also the theorem [1]-6 on page 277.

THEOREM 2.3.13. Let[2.25] $A \in \mathcal{C}_\mu$. Then there exist set $B$ such that

(2.3.32) $$A \subseteq B$$

(2.3.33) $$\mu(A) = \mu(B)$$

(2.3.34) $$B = \bigcap_n B_n$$

(2.3.35) $$B_1 \supseteq B_2 \supseteq ... \supseteq B_n \supseteq ...$$

(2.3.36) $$B_n = \bigcup_k B_{nk}$$

(2.3.37) $$B_{nk} \in \mathcal{R}(\mathcal{S})$$

(2.3.38) $$\mu(B_{nk}) < \mu(A) + \frac{1}{n}$$

(2.3.39) $$B_{n1} \subseteq B_{n2} \subseteq ... \subseteq B_{nk} \subseteq ...$$

PROOF. According to the definition 2.3.4 and to the theorem 2.1.11, for any $n$, there exists the set $C_n$ such that

(2.3.40) $$A \subseteq C_n$$

(2.3.41) $$\mu(C_n) < \mu(A) + \frac{1}{n}$$

(2.3.42) $$C_n = \bigcup_r \Delta_{nr} \quad \Delta_{nr} \in \mathcal{S}$$

Let

(2.3.43) $$B_n = \bigcap_{k=1}^n C_k$$

The statement (2.3.35) follows from (2.3.43). The statement

(2.3.44) $$A \subseteq B_n$$

follows from (2.3.40), (2.3.43). The statement (2.3.32) follows from (2.3.34), (2.3.44). From (2.3.42), (2.3.43), it follows that

(2.3.45) $$B_n = \bigcup_r \delta_{nr} \quad \delta_{nr} \in \mathcal{S}$$

Let

(2.3.46) $$B_{nk} = \bigcup_{r=1}^k \delta_{nr}$$

The statement (2.3.39) follows from (2.3.46). The statement (2.3.36) follows from (2.3.45), (2.3.46). The statement (2.3.37) follows from (2.3.45), (2.3.46) and from the theorem 2.1.11.

From the statement (2.3.32) and from the theorem 2.2.6, it follows that

(2.3.47) $$\mu(A) \leq \mu(B)$$

---

[2.25] See also the theorem [1]-8 on pages 277, 278.

From statements (2.3.34), (2.3.41), (2.3.43) and from the theorem 2.2.6, it follows that

(2.3.48) $$\mu(B) \leq \mu(B_n) \leq \mu(C_n) < \mu(A) + \frac{1}{n}$$

Since $n$ is arbitrary, then the statement

(2.3.49) $$\mu(B) \leq \mu(A)$$

follows from the statement (2.3.48). The statement (2.3.33) follows from statements (2.3.47), (2.3.49). The statement (2.3.38) follows from statements (2.3.36), (2.3.48) and from the theorem 2.2.6. □

CHAPTER 3

# Measurable Map into Abelian $\Omega$-Group

### 3.1. Measurable Map

DEFINITION 3.1.1. Minimal $\sigma$-algebra $\mathcal{B}(A)$ generated by the set of all open balls of normed $\Omega$-group $A$ is called **Borel algebra**.[3.1] A set belonging to Borel algebra is called **Borel set** or $B$-**set**. $\square$

DEFINITION 3.1.2. Let $\mathcal{C}_X$ be $\sigma$-algebra of sets of set $X$. Let $\mathcal{C}_Y$ be $\sigma$-algebra of sets of set $Y$. The map
$$f : X \to Y$$
is called $(\mathcal{C}_X, \mathcal{C}_Y)$-measurable[3.2] if for any set $C \in \mathcal{C}_Y$
$$f^{-1}(C) \in \mathcal{C}_X$$

$\square$

EXAMPLE 3.1.3. Let $\mu$ be a $\sigma$-additive measure defined on the set $X$. Let $\mathcal{C}_\mu$ be $\sigma$-algebra of sets measurable with respect to measure $\mu$. Let $\mathcal{B}(A)$ be Borel algebra of normed $\Omega$-group $A$. The map
$$f : X \to A$$
is called $\mu$-**measurable**[3.3] if for any set $C \in \mathcal{B}(A)$
$$f^{-1}(C) \in \mathcal{C}_\mu$$

$\square$

EXAMPLE 3.1.4. Let $\mathcal{B}(A)$ be Borel algebra of normed $\Omega_1$-group $A$. Let $\mathcal{B}(B)$ be Borel algebra of normed $\Omega_2$-group $B$. The map
$$f : A \to B$$
is called **Borel-measurable**[3.4] if for any set $C \in \mathcal{B}(B)$
$$f^{-1}(C) \in \mathcal{B}(A)$$

$\square$

THEOREM 3.1.5. *Let $\mathcal{C}_X$ be $\sigma$-algebra of sets of set $X$. Let $\mathcal{C}_Y$ be $\sigma$-algebra of sets of set $Y$. Let $\mathcal{C}_Z$ be $\sigma$-algebra of sets of set $Z$.*

---

[3.1]See remark in [1], p. 36. According to the remark 2.1.3, the set of closed balls also generates Borel algebra.

[3.2]See similar definition in [1], page 284.

[3.3]See similar definition in [1], pp. 284, 285, definition 1. If the measure $\mu$ is defined on the set $X$ by context, then we also call the map
$$f : X \to A$$
**measurable**.

[3.4]See similar definition in [1], pp. 284, example 1.

3.1.5.1: *Let the map*
$$f : X \to Y$$
*be $(\mathcal{C}_X, \mathcal{C}_Y)$-measurable map.*

3.1.5.2: *Let the map*
$$g : Y \to Z$$
*be $(\mathcal{C}_Y, \mathcal{C}_Z)$-measurable map.*

*Then the map*
$$g \circ f : X \to Z$$
*is $(\mathcal{C}_X, \mathcal{C}_Z)$-measurable map.*[3.5]

PROOF. Let $A \in \mathcal{C}_Z$. According to the definition 3.1.2 and statement 3.1.5.2,
$$g^{-1}(A) \in \mathcal{C}_Y$$
According to the definition 3.1.2 and statement 3.1.5.1,
$$f^{-1}(g^{-1}(A)) = (gf)^{-1}(A) \in \mathcal{C}_X$$
Therefore, the map
$$g \circ f : X \to Z$$
is $(\mathcal{C}_X, \mathcal{C}_Z)$-measurable map. □

## 3.2. Simple Map

DEFINITION 3.2.1. *Let $\mu$ be a $\sigma$-additive measure defined on the set $X$. The map*
$$f : X \to A$$
*into normed $\Omega$-group $A$ is called **simple map**, if this map is $\mu$-measurable and its range is finite or countable set.* □

THEOREM 3.2.2. *Let range $y_1$, $y_2$, ... of the map*
$$f : X \to A$$
*be finite or countable set. The map $f$ is $\mu$-measurable iff all sets*
$$A_n = \{x : f(x) = y_n\}$$
*are $\mu$-measurable.*[3.6]

PROOF. Each set $\{y_n\}$ is a Borel set. Since $A_n$ is preimage of the set $\{y_n\}$, then $A_n$ is $\mu$-measurable if map $f$ is $\mu$-measurable. Therefore, the condition of the theorem is necessary.

Let all sets $A_n$ be $\mu$-measurable. The preimage $f^{-1}(B)$ of Borel set $B \subset A$ is $\mu$-measurable, since it is union
$$\bigcup_{y_n \in B} A_n$$
of no more than countably many $\mu$-measurable sets $A_n$. Therefore, map $f$ is $\mu$-measurable and the condition of the theorem is sufficient. □

---

[3.5]See similar theorem in [1], page 286, theorem 3.
[3.6]See similar theorem in [1], page 286, theorem 4.

THEOREM 3.2.3. *Let $\mu$ be a $\sigma$-additive measure defined on the set $X$. Let*
$$f : X \to A$$
*be simple map into normed $\Omega$-group $A$. Let the map $f$ be $\mu$-measurable on the set $X_i \subset X$, $i = 1, 2, \dots$ . Then the map $f$ is $\mu$-measurable on the set $\bigcup_i X_i$ .*

PROOF. Let $y_1, y_2, \dots$ be the range of the map $f$. According to the theorem 3.2.2, the set
$$Y_{i \cdot k} = \{x \in X_i : f(x) = y_k\}$$
is $\mu$-measurable. According to definitions 2.1.13, 2.2.10, the set
$$Y_k = \bigcup_i Y_{i \cdot k} = \{x \in \bigcup_i X_i : f(x) = y_k\}$$
is $\mu$-measurable. According to the theorem 3.2.2, the map $f$ is $\mu$-measurable on the set $\bigcup_i X_i$ . $\square$

THEOREM 3.2.4. *Let $\mu$ be a $\sigma$-additive measure defined on the set $X$. Let*
$$f : X \to A$$
$$g : X \to A$$
*be simple maps into normed $\Omega$-group $A$. Then map*
$$h = f + g$$
*is simple map.*[3.7]

PROOF. According to the definition 3.2.1, simple maps $f$ and $g$ have finite or countable ranges. Let $y_1, y_2, \dots$ be the range of the map $f$. Let $z_1, z_2, \dots$ be the range of the map $g$. Then the range of the map $h$ consists of values
$$c_{ij} = y_i + z_j$$
and it is finite or countable set. For each $c_{ij}$, the set
$$\{x : h(x) = c_{ij}\} = \bigcup_{y_i + z_j = c_{ij}} \{x : f(x) = y_i\} \cap \{x : g(x) = z_j\}$$
is $\mu$-measurable. Therefore, map $h$ is simple map. $\square$

THEOREM 3.2.5. *Let $\mu$ be a $\sigma$-additive measure defined on the set $X$. Let $\omega \in \Omega$ be n-ari operation. Let*
$$f_i : X \to A \quad i = 1, \dots, n$$
*be simple maps into normed $\Omega$-group $A$. Then map*
$$h = f_1 \dots f_n \omega$$
*is simple map.*

---

[3.7]See similar theorem in [1], page 287, theorem 6.

PROOF. According to the definition 3.2.1, simple maps $f_i$, $i = 1, ..., n$, have finite or countable ranges. Let $y_{i \cdot 1}, y_{i \cdot 2}, ...$ be the range of the map $f_i$. Then the range of the map $h$ consists of values

$$y_{i_1...i_n} = y_{1 \cdot i_1}...y_{n \cdot i_n} \omega$$

and it is finite or countable set. For each $y_{i_1...i_n}$, the set

$$\{x : h(x) = y_{i_1...i_n}\} = \bigcup_{y_{1 \cdot i_1}...y_{n \cdot i_n}\omega = y_{i_1...i_n}} \bigcap_{j=1}^{n} \{x : f_j(x) = y_{j \cdot j_i}\}$$

is $\mu$-measurable. Therefore, map $h$ is simple map. $\square$

THEOREM 3.2.6. Let

$$f : A_1 \dashrightarrow A_2$$

be representation of $\Omega_1$-group $A_1$ with norm $\|x\|_1$ in $\Omega_2$-group $A_2$ with norm $\|x\|_2$. Let $\mu$ be a $\sigma$-additive measure defined on the set $X$. Let

$$g_i : X \to A_i \quad i = 1, 2$$

be simple map. Then map

$$h = f_X(g_1)(g_2)$$

is simple map.

PROOF. According to the definition 3.2.1, simple maps $g_1$ and $g_2$ have finite or countable ranges. Let $y_{1 \cdot 1}, y_{1 \cdot 2}, ...$ be the range of the map $g_1$. Let $y_{2 \cdot 1}, y_{2 \cdot 2}, ...$ be the range of the map $g_2$. Then the range of the map $h$ consists of values

$$c_{ij} = f(y_{1 \cdot i})(y_{2 \cdot j})$$

and it is finite or countable set. For each $c_{ij}$, the set

$$\{x : h(x) = c_{ij}\} = \bigcup_{f(y_{1 \cdot i})(y_{2 \cdot j}) = c_{ij}} \{x : g_1(x) = y_{1 \cdot i}\} \cap \{x : g_2(x) = y_{2 \cdot j}\}$$

is $\mu$-measurable. Therefore, map $h$ is simple map. $\square$

### 3.3. Operations on Measurable Maps

THEOREM 3.3.1. Let $A$ be normed $\Omega$-group. Let $\{f_n\}$ be a sequence of $\mu$-measurable maps

$$f_n : X \to A$$

Let

$$f : X \to A$$

be such map that

(3.3.1) $$f(x) = \lim_{n \to \infty} f_n(x)$$

for every $x$. Then the map $f$ is $\mu$-measurable map.

PROOF. Let us prove the following equation

(3.3.2) $$\{x : f(x) \in B_c(a, R)\} = \bigcup_k \bigcup_n \bigcap_{m > n} \{x : f_m(x) \in B_o(a, R + 1/k)\}$$

- Let $f(x) \in B_c(a, R)$. According to the definition [4]-2.1.15

(3.3.3) $$\|f(x) - a\| \leq R$$

From the equation (3.3.1), it follows that for any $k$ there exists such $n$ that for $m > n$

(3.3.4) $$\|f(x) - f_m(x)\| < 1/k$$

From equations (3.3.3), (3.3.4), it follows that

(3.3.5) $$\begin{aligned}\|f_m(x) - a\| &= \|f_m(x) - f(x) + f(x) - a\| \\ &\leq \|f_m(x) - f(x)\| + \|f(x) - a\| \\ &< R + 1/k\end{aligned}$$

From the definition [4]-2.1.14 and the equation (3.3.5), it follows that for any $k$ there exists such $n$ that for $m > n$

$$f_m(x) \in B_o(a, R + 1/k)$$

Therefore, we proved that

(3.3.6) $$\{x : f(x) \in B_c(a, R)\} \subseteq \bigcup_k \bigcup_n \bigcap_{m>n} \{x : f_m(x) \in B_o(a, R + 1/k)\}$$

- Let
$$x \in \bigcup_k \bigcup_n \bigcap_{m>n} \{x : f_m(x) \in B_o(a, R + 1/k)\}$$

Then there exist $k$, $n$ such that for $m > n$

(3.3.7) $$f_m(x) \in B_o(a, R + 1/k)$$

From the definition [4]-2.1.14 and the equation (3.3.7), it follows that

(3.3.8) $$\|f_m(x) - a\| < R + 1/k$$

From the theorem [5]-2, page 56, the equation (3.3.1) and the inequality (3.3.8), it follows that

(3.3.9) $$\|f(x) - a\| \leq R$$

From the definition [4]-2.1.15 and inequality (3.3.9), it follows that $f(x) \in B_c(a, R)$. Therefore, we proved that

(3.3.10) $$\bigcup_k \bigcup_n \bigcap_{m>n} \{x : f_m(x) \in B_o(a, R + 1/k)\} \subseteq \{x : f(x) \in B_c(a, R)\}$$

- The equation (3.3.2) follows from equations (3.3.6), (3.3.10).

Since maps $f_n$ are measurable maps, then, according to the example 3.1.3, sets

$$\{x : f_m(x) \in B_o(a, R + 1/k)\}$$

are measurable sets. Since the set of measurable sets is $\sigma$-algebra, then, from the equation (3.3.2), it follows that the set

$$\{x : f(x) \in B_c(a, R)\}$$

is measurable set. According to the example 3.1.3, the map $f$ is measurable map. □

THEOREM 3.3.2. *Let A be normed Ω-group. Let range of the map*

$$f : X \to A$$

*be compact set. The map f is µ-measurable iff it can be represented as the limit of the uniformly convergent sequence of simple maps.*[3.8]

PROOF. Let map $f$ be limit of a convergent sequence of simple maps. According to the definition 3.2.1, simple map is $\mu$-measurable. According to the theorem 3.3.1, limit of sequence of simple maps is $\mu$-measurable. Therefore, the map $f$ is $\mu$-measurable.

Let the map $f$ be $\mu$-measurable. For given $n$, consider a set of open balls

$$B_f = \{B_o(y, 1/n) : \exists x, y = f(x)\}$$

The set $B_f$ is open cover of the range of the map $f$. Therefore, there exists finite set $B'_f$, $B'_f \subset B_f$, which is open cover of the range of the map $f$. For $x \in X$, we choose an open ball

(3.3.11) $$B_o(y, 1/n) \in B_f \quad f(x) \in B_o(y, 1/n)$$

and we assume

(3.3.12) $$f_n(x) = y$$

The map $f_n$ is simple map. From equations (3.3.11), (3.3.12) and the definition [4]-2.1.14, it follows that

(3.3.13) $$\|f_n(x) - f(x)\| < \frac{1}{n}$$

From the equation (3.3.13) it follows that the sequence of simple maps $f_n$ uniformly converges to the map $f$. □

THEOREM 3.3.3. *Let $\mu$ be a $\sigma$-additive measure defined on the set $X$. Let $A$ be normed Ω-group. Let range of the map*

$$f : X \to A$$

*be compact set. Let the map $f$ be $\mu$-measurable on the set $X_i \subset X$, $i = 1, ..., n$,*

$$i \neq j \implies X_i \cap X_j = \emptyset$$

*Then the map $f$ is $\mu$-measurable on the set $\bigcup_i X_i$.*

PROOF. According to the theorem 3.3.2, for each $i$, there exists sequence of simple maps

$$f_{i \cdot k} : X_i \to A$$

which uniformly converges to map $f$ on set $X_i$. For each $k$, consider the map

$$f_k : \bigcup_i X_i \to A$$

defined by rule

$$x \in X_i \implies f_k(x) = f_{i \cdot k}(x)$$

---

[3.8] See similar theorem in [1], page 286, theorem 5.

Since the map $f_k$ has finite or countable range over set $X_i$, then the map $f_k$ has finite or countable range over set $\bigcup_i X_i$. According to the theorem 3.2.3, the map $f_k$ is measurable on the set $\bigcup_i X_i$.

Since sequence of maps $f_k$ uniformly converges to map $f$ on set $X_i$, then, according to the definition [4]-2.6.3, for given $\epsilon \in R,\ \epsilon > 0$, there exists $K_i$ such that the condition $k > K_i$ implies that
$$\|f_k(x) - f(x)\| < \epsilon$$
for any $x \in X_i$. Let
$$K = \max(K_1, ..., K_n)$$
Then for given $\epsilon \in R,\ \epsilon > 0$, the condition $k > K$ implies that
$$\|f_k(x) - f(x)\| < \epsilon$$
for any $x \in \bigcup_i X_i$. According to the definition [4]-2.6.3, the sequence of maps $f_k$ uniformly converges to map $f$ on set $\bigcup_i X_i$. According to the theorem 3.3.2, the map $f$ is $\mu$-measurable on the set $\bigcup_i X_i$. □

THEOREM 3.3.4. *Let $\mu$ be a $\sigma$-additive measure defined on the set $X$. Let $A$ be complete $\Omega$-group. Let*
$$f : X \to A$$
$$g : X \to A$$
*be $\mu$-measurable maps with compact range. Then map*
(3.3.14) $$h = f + g$$
*is $\mu$-measurable map with compact range.*[3.9]

PROOF. According to the theorem 3.3.2, there exists the sequence of simple maps $f_n$ uniformly converging to the map $f$ and there exists the sequence of simple maps $g_n$ uniformly converging to the map $g$. For any $n$, the map
$$h_n = f_n + g_n$$
is simple map according to the theorem 3.2.4. According to the theorem [4]-2.6.7, the sequence of maps $h_n$ converges uniformly to the map $h$. According to the theorem 3.3.1, the map $h$ is $\mu$-measurable map.

Consider open cover $O_h$ of the range of the map $h$. According to the definition [4]-2.2.2, open cover $O_h$ contains a set of open balls
$$B_h = \{B_o(y, 2\epsilon_x) : \exists x, y = h(x)\}$$
Consider a set of open balls
$$B_f = \{B_o(y, \epsilon_x) : \exists x, y = f(x)\}$$
The set $B_f$ is open cover of the range of the map $f$. Therefore, there exists finite set $B'_f,\ B'_f \subset B_f$, which is open cover of the range of the map $f$. Let
$$I_f = \{x \in X : y = f(x), B_o(y, \epsilon_x) \in B'_f\}$$

---

[3.9]See similar theorem in [1], page 287, theorem 6.

Consider a set of open balls
$$B_g = \{B_o(y, \epsilon_x) : \exists x, y = g(x)\}$$
The set $B_g$ is open cover of the range of the map $g$. Therefore, there exists finite set $B'_g$, $B'_g \subset B_g$, which is open cover of the range of the map $g$. Let
$$I_g = \{x \in X : y = g(x), B_o(y, \epsilon_x) \in B'_g\}$$
Let
$$I_h = I_f \cup I_g$$
For any $x \in X$, there exists $x' \in I_h$ such that

(3.3.15) $$f(x) \in B_o(f(x'), \epsilon_{x'}) \quad g(x) \in B_o(g(x'), \epsilon_{x'})$$

According to the theorem [4]-2.4.2, the statement
$$h(x) \in B_o(h(x'), 2\epsilon_{x'})$$
follows from (3.3.14), (3.3.15). Therefore, the set of open balls
$$B'_h = \{B_o(y, 2\epsilon_x) : \exists x \in I_h, y = h(x)\} \subseteq B_h \subseteq O_h$$
is finite open cover of the range of the map $h$. According to the definition [4]-2.2.15, the range of the map $h$ is compact set. □

THEOREM 3.3.5. *Let $\mu$ be a $\sigma$-additive measure defined on the set $X$. Let $\omega \in \Omega$ be n-ari operation. Let*
$$f_i : X \to A \quad i = 1, ..., n$$
*be $\mu$-measurable maps into complete $\Omega$-group $A$. Then map*
$$h = f_1...f_n\omega$$
*is $\mu$-measurable map.*

PROOF. According to the theorem 3.3.2, there exists the sequence of simple maps $f_{i \cdot m}$ uniformly converging to the map $f_i$. For any $m$, the map
$$h_m = f_{1 \cdot m}...f_{n \cdot m}\omega$$
is simple map according to the theorem 3.2.5. According to the theorem [4]-2.6.8, the sequence of maps $h_n$ converges uniformly to the map $h$. □

THEOREM 3.3.6. *Let $\mu$ be a $\sigma$-additive measure defined on the set $X$. Let*
$$f : A_1 \dashrightarrow A_2$$
*be representation of $\Omega_1$-group $A_1$ with norm $\|x\|_1$ in $\Omega_2$-group $A_2$ with norm $\|x\|_2$. Let*
$$g_i : X \to A_i \quad i = 1, 2$$
*be $\mu$-measurable map. Then map*
$$h = f_X(g_1)(g_2)$$
*is $\mu$-measurable map.*

PROOF. According to the theorem 3.3.2, there exists the sequence of simple maps $f_{i \cdot m}$ uniformly converging to the map $f_i$. For any $m$, the map $f_X(g_{1 \cdot n})(g_{2 \cdot n})$ is simple map according to the theorem 3.2.6. According to the theorem [4]-3.2.2, the sequence of maps $h_n$ converges uniformly to the map $h$. □

## 3.4. Convergence Almost Everywhere

DEFINITION 3.4.1. Let $\mu$ be a $\sigma$-additive measure defined on the set $X$. The sequence[3.10]
$$f_n : X \to A$$
of $\mu$-measurable map into $\Omega$-group $A$ **converges almost everywhere** if
(3.4.1) $$f(x) = \lim_{n \to \infty} f_n(x)$$
for almost all $x \in X$, i. e. if the set of $x$ for which (3.4.1) fails to hold has measure zero. □

THEOREM 3.4.2. *Let $\mu$ be a $\sigma$-additive measure defined on the set $X$. If the sequence*
$$f_n : X \to A$$
*of $\mu$-measurable map into $\Omega$-group $A$ converges to map*
$$f : X \to A$$
*almost everywhere, then the map $f$ is also $\mu$-measurable map.*

PROOF. Let
$$A = \{x \in X : \lim_{n \to \infty} f_n(x) = f(x)\} \quad B = X \setminus A$$
According to the definition 3.4.1, $\mu(B) = 0$.

LEMMA 3.4.3. *The map $f$ is $\mu$-measurable on the set $B$.*

PROOF. The statement follows from the definition 3.1.2 and the theorem 2.2.8. ⊙

LEMMA 3.4.4. *The map $f$ is $\mu$-measurable on the set $A$.*

PROOF. According to the remark 2.1.3, the set $A$ is $\mu$-measurable set. According to the theorem 3.3.1, the map $f$ is $\mu$-measurable on the set $A$. ⊙

From the theorem 3.3.3 and from lemmmas 3.4.3, 3.4.4 it follows that $f$ is $\mu$-measurable. □

THEOREM 3.4.5 (Dmitri Fyodorovich Egorov). *Let[3.11] a sequence of $\mu$-measurable maps*
$$f_n : X \to A$$
*converge almost everywhere on measurable set $E$ to a map $f$. Then, for any $\delta > 0$, there exists a measurable set $E_\delta \subset E$ such that*

3.4.5.1: $\mu(E_\delta) > \mu(E) - \delta$
3.4.5.2: *The sequence $f_n$ converges uniformly on the set $E_\delta$.*

PROOF. According to the theorem 3.4.2, the map $f$ is $\mu$-measurable on the set $E$. Let
(3.4.2) $$E_n^m = \bigcap_{i \geq n} \left\{ x \in E : \|f_i(x) - f(x)\| < \frac{1}{m} \right\}$$

---

[3.10]See also definition [1]-3 on page 289.
[3.11]See also theorem [1]-12 on page 290.

$$(3.4.3) \qquad E^m = \bigcup_{n=1}^{\infty} E_n^m$$

The statement

$$(3.4.4) \qquad E_1^m \subset E_2^m \subset ...$$

follows from (3.4.2). From the statements (3.4.3), (3.4.4) and from the theorem 2.2.13, it follows that for any $m$ and for any $\delta > 0$, there exists $n_0(m)$ such that

$$(3.4.5) \qquad \mu(E^m - E_{n_0(m)}^m) < \frac{\delta}{2^m}$$

Let

$$(3.4.6) \qquad E_\delta = \bigcap_{m=1}^{\infty} E_{n_0(m)}^m$$

If $x \in E_\delta$, then, for any $m = 1, 2, ...$ and for any $i > n_0(m)$, the inequality

$$(3.4.7) \qquad \|f_i(x) - f(x)\| < \frac{1}{m}$$

follows from (3.4.2), (3.4.6). The statement 3.4.5.2 follows from the inequality (3.4.7) and the definition [4]-2.6.3.

If $x \in E \setminus E^m$, then, from (3.4.2), (3.4.3), it follows that there are arbitrary large values of $i$ such that

$$(3.4.8) \qquad \|f_i(x) - f(x)\| > \frac{1}{m}$$

Therefore, the sequence $f_n(x)$ at the point $x \in E \setminus E^m$ does not converge to $f(x)$. Since $f_n$ converges to $f$ almost everywhere, then, from the theorem 2.2.8, it follows that

$$(3.4.9) \qquad \mu(E \setminus E^m) = 0$$

The inequality

$$(3.4.10) \qquad \mu(E - E_{n_0(m)}^m) = \mu(E^m - E_{n_0(m)}^m) < \frac{\delta}{2^m}$$

follows from (3.4.5), (3.4.9). The statement 3.4.5.1 follows from the inequality

$$\mu(E \setminus E_\delta) = \mu\left(E \setminus \bigcap_{m=1}^{\infty} E_{n_0(m)}^m\right) = \mu\left(\bigcup_{m=1}^{\infty} (E \setminus E_{n_0(m)}^m)\right)$$
$$\leq \sum_{m=1}^{\infty} \mu(E \setminus E_{n_0(m)}^m) < \sum_{m=1}^{\infty} \frac{\delta}{2^m} = \delta$$

based on the inequality (3.4.10). $\qquad\square$

CHAPTER 4

# Integral of Map into Abelian $\Omega$-Group

Let $\mu$ be a $\sigma$-additive measure defined on the set $X$. Let effective representation of real field $R$ in complete Abelian $\Omega$-group $A$ be defined.[4.1] In this chapter, we consider the definition of Lebesgue integral of $\mu$-measurable map

$$f : X \to A$$

## 4.1. Integral of Simple Map

DEFINITION 4.1.1. Let $a_i$ be a sequence of $A$-numbers. If

$$\sum_{i=1}^{\infty} \|a^i\| < \infty$$

then we say that the **series**

$$\sum_{i=1}^{\infty} a^i$$

**converges normally**.[4.2] □

DEFINITION 4.1.2. For simple map

$$f : X \to A$$

consider series

(4.1.1) $$\sum_n \mu(F_n) f_n$$

where

- The set $\{f_1, f_2, ...\}$ is domain of the map $f$
- Since $n \ne m$, then $f_n \ne f_m$
- $F_n = \{x \in X : f(x) = f_n\}$

Simple map

$$f : X \to A$$

is called **integrable map** over the set $X$ if series (4.1.1) converges normally.[4.3] Since the map $f$ is integrable map, then sum of series (4.1.1) is called **integral of map** $f$ over the set $X$

(4.1.2) $$\int_X d\mu(x) f(x) = \sum_n \mu(F_n) f_n$$

□

---

[4.1]In other words, $\Omega$-group $A$ is $R$-vector space.
[4.2]See also the definition of normal convergence of the series on page [6]-12.
[4.3]See similar definition in [1], p. 294.

THEOREM 4.1.3. *Let*
$$f : X \to A$$
*be simple map. Let $f$ get value*[4.4] *$f_n$ on the set $F_n \subset X$. Let $X = \bigcup_n F_n$, $F_n \cap F_m = \emptyset$. The map $f$ is integrable map iff series*

(4.1.3) $$\sum_n \mu(F_n) f_n$$

*converges normally. Then*

(4.1.4) $$\int_X d\mu(x) f(x) = \sum_n \mu(F_n) f_n$$

PROOF. Let
$$X_i = \{x \in X : f(x) = f_i\}$$
Then

(4.1.5) $$X_i = \bigcup_{f(X_i) = f(F_n)} F_n$$

From (4.1.5), it follows that

(4.1.6) $$\mu(X_i) = \sum_{f(X_i) = f(F_n)} \mu(F_n)$$

Since series (4.1.3) converges normally, then

(4.1.7) $$\sum_n \mu(F_n) f_n = \sum_i \left( \sum_{f(X_i) = f(F_n)} \mu(F_n) \right) f(X_i) = \sum_i \mu(X_i) f(X_i)$$

follows from (4.1.6). (4.1.4) follows from (4.1.2), (4.1.7). □

THEOREM 4.1.4. *Let*
$$f : X \to A$$
$$g : X \to A$$
*be simple maps.*[4.5] *Since there exist integrals*
$$\int_X d\mu(x) f(x)$$
$$\int_X d\mu(x) g(x)$$
*then there exists integral*
$$\int_X d\mu(x)(f(x) + g(x))$$
*and*

(4.1.8) $$\int_X d\mu(x)(f(x) + g(x)) = \int_X d\mu(x) f(x) + \int_X d\mu(x) g(x)$$

---

[4.4]We do not request $f_n \neq f_m$ since $n \neq m$. However we request $F_n \cap F_m = \emptyset$. See also lemma in [1], pages 294, 295.

[4.5]See similar statement in [1], theorem 1, p. 295.

PROOF. Let $f$ get value $f_n$ on the set $F_n \subset X$. Let
$$(4.1.9) \qquad F_n \cap F_m = \emptyset$$
Let $g$ get value $g_k$ on the set $G_k \subset X$. Let
$$(4.1.10) \qquad G_k \cap G_l = \emptyset$$
The equation
$$(4.1.11) \qquad \mu(F_n) = \sum_k \mu(F_n \cap G_k)$$
follows from the equation
$$F_n = \bigcup_k F_n \cap G_k$$
and condition (4.1.10). The equation
$$(4.1.12) \qquad \mu(G_k) = \sum_n \mu(F_n \cap G_k)$$
follows from the equation
$$G_k = \bigcup_n F_n \cap G_k$$
and condition (4.1.9). The condition
$$(F_n \cap G_k) \cap (F_m \cap G_l) = \emptyset$$
follows from conditions (4.1.9), (4.1.10).

4.1.4.1: Since the map $f$ is integrable map, then the equation
$$(4.1.13) \qquad \int_X d\mu(x) f(x) = \sum_n \mu(F_n) f_n = \sum_n \left( \sum_k \mu(F_n \cap G_k) \right) f_n$$
$$= \sum_n \sum_k \mu(F_n \cap G_k) f_n$$
follows from (4.1.4), (4.1.11) and series
$$\sum_n \sum_k \mu(F_n \cap G_k) f_n$$
converges normally.

4.1.4.2: Since the map $g$ is integrable map, then the equation
$$(4.1.14) \qquad \int_X d\mu(x) g(x) = \sum_k \mu(G_k) g_k = \sum_k \left( \sum_n \mu(G_k \cap F_n) \right) g_k$$
$$= \sum_k \sum_n \mu(G_k \cap F_n) g_k$$
$$= \sum_n \sum_k \mu(F_n \cap G_k) g_k$$
follows from (4.1.4), (4.1.12) and series
$$\sum_k \sum_n \mu(F_n \cap G_k) g_k = \sum_n \sum_k \mu(F_n \cap G_k) g_k$$
converges normally.

From statements 4.1.4.1, 4.1.4.2, it follows that

$$\int_X d\mu(x)f(x) + \int_X d\mu(x)g(x)$$
(4.1.15)
$$= \sum_n \sum_k \mu(F_n \cap G_k)f_n + \sum_n \sum_k \mu(F_n \cap G_k)g_k$$
$$= \sum_n \sum_k \mu(F_n \cap G_k)(f_n + g_k)$$

and series

$$\sum_n \sum_k \mu(F_n \cap G_k)(f_n + g_k)$$

converges normally. Therefore, according to the theorem 4.1.3, the equation (4.1.8) follows from (4.1.15). □

THEOREM 4.1.5. *Let*

$$f : X \to A$$

*be simple map. Integral*

$$\int_X d\mu(x)f(x)$$

*exists iff integral*

$$\int_X d\mu(x)\|f(x)\|$$

*exists. Then*

(4.1.16) $$\left\|\int_X d\mu(x)f(x)\right\| \le \int_X d\mu(x)\|f(x)\|$$

PROOF. The equation

(4.1.17) $$\sum_n \|\mu(F_n)f_n\| = \sum_n \mu(F_n)\|f_n\|$$

follows from the equation

$$\|\mu(F_n)f_n\| = \mu(F_n)\|f_n\|$$

Therefore, series

$$\sum_n \mu(F_n)f_n$$

converges normally iff series

$$\sum_n \mu(F_n)\|f_n\|$$

converges. inequality (4.1.16) follows from the inequality

$$\left\|\sum_n \mu(F_n)f_n\right\| \le \sum_n \|\mu(F_n)f_n\|$$

the equation (4.1.17) and the definition of integral. □

THEOREM 4.1.6. *Let $\omega \in \Omega$ be n-ary operation in Abelian $\Omega$-group $A$. Let simple map*

$$f_i : X \to A \quad i = 1, ..., n$$

*be integrable map. Then map*

$$h = f_1...f_n\omega$$

*is integrable map and*

$$\left\| \int_X d\mu(x) h(x) \right\| \leq \int_X d\mu(x)(\|\omega\| \|f_1(x)\|...\|f_n(x)\|)$$

PROOF. The theorem follows from the theorem 4.1.5 and the inequality [4]-(2.1.6). □

THEOREM 4.1.7. *Let $\mu$ be a $\sigma$-additive measure defined on the set $X$. Let*

$$f : A_1 \relbar\joinrel\twoheadrightarrow A_2$$

*be representation of $\Omega_1$-group $A_1$ with norm $\|x\|_1$ in $\Omega_2$-group $A_2$ with norm $\|x\|_2$. Let*

$$g_i : X \to A_i \quad i = 1, 2$$

*be simple integrable map. Then map*

$$h = f_X(g_1)(g_2)$$

*is integrable map and*

$$\left\| \int_X d\mu(x) h(x) \right\|_2 \leq \int_X d\mu(x)(\|f\| \|g_1(x)\|_1 \|g_2(x)\|_2)$$

PROOF. The theorem follows from the theorem 4.1.5 and the inequality [4]-(3.1.2). □

THEOREM 4.1.8. *Let $\mu$ be a $\sigma$-additive measure defined on the set $X$. Let*

$$f : A_1 \relbar\joinrel\twoheadrightarrow A_2$$

*be effective representation of $\Omega_1$-group $A_1$ with norm $\|x\|_1$ in $\Omega_2$-group $A_2$ with norm $\|x\|_2$. Let transformations of the representation $f$ be morphism of representation*

$$R \relbar\joinrel\twoheadrightarrow A_2$$

*Let*

$$g_2 : X \to A_2$$

*be simple integrable map into normed $\Omega$-group $A_2$. Then map*

$$h = a_1 g_2$$

*is integrable map and*

(4.1.18) $$\int_X d\mu(x)(a_1 g_2(x)) = a_1 \int_X d\mu(x) g_2(x)$$

PROOF. Let $g_2$ get value $g_{2 \cdot k}$ on the set $G_k \subset X$. Let

$$G_k \cap G_l = \emptyset$$

Since the map $g_2$ is integrable map, then the equation

$$(4.1.19) \qquad \int_X d\mu(x) g_2(x) = \sum_k \mu(G_k) g_{2 \cdot k}$$

follows from (4.1.4) and series in (4.1.19) converges normally. According to the theorem [4]-3.1.2,

$$\|a_1 g_{2 \cdot j}\|_2 \le \|f\| \|a_1\|_1 \|g_{2 \cdot j}\|_2$$

Therefore, the series

$$\sum_k \mu(G_k)(a_1 g_{2 \cdot k})$$

converges normally. Since the transformation generated by $a_1 \in A_1$ is morphism of representation[4.6]

$$R \relbar\joinrel\twoheadrightarrow A_2$$

the equation (4.1.18) follows from the equation

$$\sum_k \mu(G_k)(a g_{2 \cdot k}) = \sum_k a(\mu(G_k) g_{2 \cdot k}) = a \sum_k \mu(G_k) g_{2 \cdot k}$$

and the theorem 4.1.3. □

THEOREM 4.1.9. *Let $\mu$ be a $\sigma$-additive measure defined on the set $X$. Let*

$$f : A_1 \relbar\joinrel\twoheadrightarrow A_2$$

*be effective representation of $\Omega_1$-group $A_1$ with norm $\|x\|_1$ in $\Omega_2$-group $A_2$ with norm $\|x\|_2$. Let transformations of the representation*

$$R \relbar\joinrel\twoheadrightarrow A_1$$

*be morphism of representation $f$. Let*

$$g_1 : X \to A_1$$

*be simple integrable map into normed $\Omega$-group $A_1$. Then map*

$$h = g_1 a_2$$

*is integrable map and*

$$(4.1.20) \qquad \int_X d\mu(x)(g_1(x) a_2) = \left( \int_X d\mu(x) g_1(x) \right) a_2$$

PROOF. Let $g_1$ get value $g_{1 \cdot k}$ on the set $G_k \subset X$. Let

$$G_k \cap G_l = \emptyset$$

Since the map $g_1$ is integrable map, then the equation

$$(4.1.21) \qquad \int_X d\mu(x) g_1(x) = \sum_k \mu(G_k) g_{1 \cdot k}$$

---

[4.6] It is important for us that for any $p \in R$, $a_1 \in A_1$, $a_2, b_2 \in A_2$,

$$p(a_1 a_2) = a_1(p a_2)$$
$$a_1(a_2 + b_2) = a_1 a_2 + a_1 b_2$$

See also the definition [2]-2.2.2.

follows from (4.1.4) and series in (4.1.21) converges normally. According to the theorem [4]-3.1.2,
$$\|g_{1\cdot j}a_2\|_1 \le \|f\|\|g_{1\cdot j}\|_1\|a_2\|_2$$
Therefore, the series
$$\sum_k \mu(G_k)(g_{1\cdot k})a_2$$
converges normally. Since the transformation of representation
$$R \xrightarrow{\ *\ } A_1$$
is morphism of representation $f$, the equation[4.7] (4.1.20) follows from the equation
$$\sum_k \mu(G_k)(g_{1\cdot k}a) = \sum_k (\mu(G_k)g_{1\cdot k})a = \left(\sum_k \mu(G_k)g_{1\cdot k}\right)a$$
and the theorem 4.1.3.  □

THEOREM 4.1.10. *Let*
$$f : X \to A$$
*be simple map such that*

(4.1.22) $$\|f(x)\| \le M$$

*Since the measure of the set $X$ is finite, then*

(4.1.23) $$\int_X d\mu(x)\|f(x)\| \le M\mu(X)$$

PROOF. Let $f$ get value $f_n$ on the set $F_n \subset X$. Let $X = \bigcup_n F_n$,
$$F_n \cap F_m = \emptyset$$
According to the theorem 4.1.3, the inequality

(4.1.24) $$\int_X d\mu(x)\|f(x)\| = \sum_n \mu(F_n)\|f_n\| \le M\sum_n \mu(F_n) = M\mu(X)$$

follows from the equation (4.1.4) and from the inequality (4.1.22). The inequality (4.1.23) follows from the inequality (4.1.24).  □

THEOREM 4.1.11. *Let*
$$f : X \to A$$
*be simple map such that*
$$\|f(x)\| \le M$$
*Since the measure of the set $X$ is finite, then map $f$ is integrable and*

(4.1.25) $$\left\|\int_X d\mu(x)f(x)\right\| \le M\mu(X)$$

PROOF. Integrability of the map $f$ follows from theorems 4.1.5, 4.1.10. The inequality (4.1.25) follows from inequalities (4.1.16), (4.1.23).  □

---

[4.7] It is important for us that for any $p \in R$, $a_1, b_1 \in A_1$, $a_2 \in A_2$,
$$p(a_1 a_2) = pf(a_1)(a_2) = f(pa_1)(a_2) = (pa_1)a_2$$
$$(a_1 + b_1)a_2 = a_1 a_2 + b_1 a_2$$
See also the definition [2]-2.2.2.

### 4.2. Integral of Measurable Map over Set of Finite Measure

DEFINITION 4.2.1. $\mu$-measurable map
$$f : X \to A$$
is called **integrable map** over the set $X$,[4.8] if there exists a sequence of simple integrable over the set $X$ maps
$$f_n : X \to A$$
converging uniformly to $f$. Since the map $f$ is integrable map, then the limit

(4.2.1) $$\int_X d\mu(x) f(x) = \lim_{n \to \infty} \int_X d\mu(x) f_n(x)$$

is called **integral of map** $f$ over the set $X$. □

THEOREM 4.2.2. *Let*
$$f : X \to A$$
*be $\mu$-measurable map.*[4.9] *Let the measure of the set $X$ be finite.*

4.2.2.1: *For any uniformly convergent sequence $f_n$ of simple integrable maps*
$$f_n : X \to A$$
*the limit* (4.2.1) *exists.*

4.2.2.2: *The limit* (4.2.1) *does not depend on choice of sequence $f_n$.*

4.2.2.3: *For simple map, the definition 4.2.1 reduces to the definition 4.1.2.*

PROOF. According to the theorem [4]-2.6.5, since the sequence $f_n$ converges uniformly to the map $f$, then for any $\epsilon \in R$, $\epsilon > 0$, there exists $N$ such that

(4.2.2) $$\|f_n(x) - f_m(x)\| < \frac{\epsilon}{\mu(X)}$$

for any $n, m > N$. According to theorems 4.1.4, 4.1.11,

(4.2.3) $$\left\| \int_X d\mu(x) f_m(x) - \int_X d\mu(x) f_n(x) \right\| = \left\| \int_X d\mu(x)(f_m(x) - f_n(x)) \right\|$$
$$\leq \frac{\epsilon}{\mu(X)} \mu(X) = \epsilon$$

follows from the inequality (4.2.2) for any $n, m > N$. According to the definition [4]-2.1.19 and the inequality (4.2.3), the sequence of integrals
$$\int_X d\mu(x) f_n(x)$$
is fundamental sequence. Therefore, there exists limit (4.2.1) and the statement 4.2.2.1 is true.

Let $f_{1 \cdot n}$, $f_{2 \cdot n}$ be fundamental sequences of simple maps uniformly convergent to the map $f$. According to the definition [4]-2.6.3, since the sequence $f_{1 \cdot n}$ converges uniformly to the map $f$, then for any $\epsilon \in R$, $\epsilon > 0$, there exists $N_1$ such that

(4.2.4) $$\|f_{1 \cdot n}(x) - f(x)\| < \frac{\epsilon}{2\mu(X)}$$

---
[4.8]See also the definition in [1], page 296.
[4.9]See also analysis of the definition in [1], pages 296, 297.

for any $n > N_1$.  According to the definition [4]-2.6.3, since the sequence $f_{2 \cdot n}$ converges uniformly to the map $f$, then for any $\epsilon \in R, \epsilon > 0,$ there exists $N_2$ such that

$$\|f_{2 \cdot n}(x) - f(x)\| < \frac{\epsilon}{2\mu(X)} \tag{4.2.5}$$

for any $n > N_2$. Let
$$N = \max(N_1, N_2)$$
From inequalities (4.2.4), (4.2.5), it follows that for given $\epsilon \in R, \epsilon > 0,$ there exists positive integer $N$ depending on $\epsilon$ and such, that

$$\begin{aligned}\|f_{1 \cdot n}(x) - f_{2 \cdot n}(x)\| &= \|f_{1 \cdot n}(x) - f(x) + f(x) - f_{2 \cdot n}(x)\| \\ &\leq \|f_{1 \cdot n}(x) - f(x)\| + \|f_{2 \cdot n}(x) - f(x)\| \\ &< \frac{\epsilon}{\mu(X)}\end{aligned} \tag{4.2.6}$$

for every $n > N$. According to theorems 4.1.4, 4.1.11,

$$\left\|\int_X d\mu(x) f_{1 \cdot n}(x) - \int_X d\mu(x) f_{2 \cdot n}(x)\right\| = \left\|\int_X d\mu(x)(f_{1 \cdot n}(x) - f_{2 \cdot n}(x))\right\| \\ \leq \frac{\epsilon}{\mu(X)}\mu(X) = \epsilon \tag{4.2.7}$$

follows from the inequality (4.2.6) for any $n > N$. According to the theorem [4]-2.1.22 and the definition [4]-2.1.17,

$$\lim_{n \to \infty} \int_X d\mu(x) f_{1 \cdot n}(x) = \lim_{n \to \infty} \int_X d\mu(x) f_{2 \cdot n}(x) \tag{4.2.8}$$

follows from the inequality (4.2.7) for any $n > N$. From the equation (4.2.8), it follows that the statement 4.2.2.2 is true.

Let $f$ be simple map. To prove the statement 4.2.2.3, it is sufficient to consider the sequence in which $f_n = f$ for any $n$. $\square$

THEOREM 4.2.3. *Let*
$$f : X \to A$$
$$g : X \to A$$
*be $\mu$-measurable maps with compact range. Since there exist integrals*
$$\int_X d\mu(x) f(x)$$
$$\int_X d\mu(x) g(x)$$
*then there exists integral*
$$\int_X d\mu(x)(f(x) + g(x))$$
*and*
$$\int_X d\mu(x)(f(x) + g(x)) = \int_X d\mu(x) f(x) + \int_X d\mu(x) g(x) \tag{4.2.9}$$

PROOF. According to the theorem 3.3.2, there exists a sequence of simple integrable over the set $X$ maps

$$f_n : X \to A$$

converging uniformly to $f$. According to the theorem 4.2.2, for any $\epsilon \in R$, $\epsilon > 0$, there exist $N_1$ such that

$$(4.2.10) \qquad \left\| \int_X d\mu(x) f(x) - \int_X d\mu(x) f_n(x) \right\| < \frac{\epsilon}{2}$$

for any $n > N_1$. According to the theorem 3.3.2, there exists a sequence of simple integrable over the set $X$ maps

$$g_n : X \to A$$

converging uniformly to $g$. According to the theorem 4.2.2, for any $\epsilon \in R$, $\epsilon > 0$, there exist $N_2$ such that

$$(4.2.11) \qquad \left\| \int_X d\mu(x) g(x) - \int_X d\mu(x) g_n(x) \right\| < \frac{\epsilon}{2}$$

for any $n > N_2$. Let

$$N = \max(N_1, N_2)$$

According to the theorem 4.1.4, for any $n > N$, there exists integral

$$(4.2.12) \qquad \int_X d\mu(x)(f_n(x) + g_n(x)) = \int_X d\mu(x) f_n(x) + \int_X d\mu(x) g_n(x)$$

According to the theorem 3.3.4, the map

$$(4.2.13) \qquad h = f + g$$

is $\mu$-measurable map and

$$(4.2.14) \qquad h(x) = \lim_{n \to \infty} f_n(x) + g_n(x)$$

If $k, n > N$, then

$$(4.2.15) \quad \begin{aligned} & \left\| \int_X d\mu(x)(f_n(x) + g_n(x)) - \int_X d\mu(x)(f_k(x) + g_k(x)) \right\| \\ &= \left\| \int_X d\mu(x) f_n(x) + \int_X d\mu(x) g_n(x) - \int_X d\mu(x) f_k(x) - \int_X d\mu(x) g_k(x) \right\| \\ &\leq \left\| \int_X d\mu(x) f_n(x) - \int_X d\mu(x) f(x) \right\| + \left\| \int_X d\mu(x) f_k(x) - \int_X d\mu(x) f(x) \right\| \\ &+ \left\| \int_X d\mu(x) g_n(x) - \int_X d\mu(x) g(x) \right\| + \left\| \int_X d\mu(x) g_k(x) - \int_X d\mu(x) g(x) \right\| \\ &< 2\epsilon \end{aligned}$$

follows from (4.2.10), (4.2.11), (4.2.12). According to the definition [4]-2.1.19, from the inequality (4.2.15), it follows that the sequence of integrals (4.2.12) is fundamental sequence. Therefore, according to statements (4.2.13), (4.2.14) and the theorem 4.2.2, there exist integral

$$(4.2.16) \qquad \int_X d\mu(x)(f(x) + g(x)) = \lim_{n \to \infty} \int_X d\mu(x)(f_n(x) + g_n(x))$$

If $n > N$, then

$$(4.2.17) \quad \begin{aligned} &\left\| \int_X d\mu(x)f(x) + \int_X d\mu(x)g(x) - \int_X d\mu(x)(f_n(x) + g_n(x)) \right\| \\ &= \left\| \int_X d\mu(x)f(x) + \int_X d\mu(x)g(x) - \int_X d\mu(x)f_n(x) - \int_X d\mu(x)g_n(x) \right\| \\ &\leq \left\| \int_X d\mu(x)f_n(x) - \int_X d\mu(x)f(x) \right\| + \left\| \int_X d\mu(x)g_n(x) - \int_X d\mu(x)g(x) \right\| \\ &< \epsilon \end{aligned}$$

follows from (4.2.10), (4.2.11), (4.2.12). From the inequality (4.2.17), it follows that the sequence of integrals (4.2.12) is fundamental sequence. According to the definition [4]-2.1.17,

$$(4.2.18) \quad \int_X d\mu(x)f(x) + \int_X d\mu(x)g(x) = \lim_{n \to \infty} \int_X d\mu(x)(f_n(x) + g_n(x))$$

The equation (4.2.9) follows from equations (4.2.16), (4.2.18). □

THEOREM 4.2.4. *Let*

$$f : X \to A$$

*be measurable map. Integral*

$$\int_X d\mu(x)f(x)$$

*exists iff integral*

$$\int_X d\mu(x)\|f(x)\|$$

*exists. Then*

$$(4.2.19) \quad \left\| \int_X d\mu(x)f(x) \right\| \leq \int_X d\mu(x)\|f(x)\|$$

PROOF. The theorem follows from the theorem 4.1.5 and from the definition 4.2.1. □

THEOREM 4.2.5. *Let $\omega \in \Omega$ be n-ary operation in Abelian $\Omega$-group $A$. Let*

$$f_i : X \to A \quad i = 1, ..., n$$

*be $\mu$-measurable map with compact range. Since map $f_i$, $i = 1, ..., n$, is integrable map, then map*

$$h = f_1...f_n\omega$$

*is integrable map and*

$$(4.2.20) \quad \left\| \int_X d\mu(x)h(x) \right\| \leq \int_X d\mu(x)(\|\omega\|\|f_1(x)\|...\|f_n(x)\|)$$

PROOF. Since the range of the map $f_i$ is compact set, then, according to the theorem [4]-2.3.11, the following value is defined

$$(4.2.21) \quad F_i = \sup \|f_i(x)\|$$

From the equation (4.2.21) and the statement [4]-2.1.9.1, it follows that

$$(4.2.22) \quad F_i \geq 0$$

According to the theorem 3.3.2, for $i = 1, ..., n$, there exists a sequence of simple integrable over the set $X$ maps

$$f_{i \cdot m} : X \to A$$

converging uniformly to $f_i$. From the equation

$$f_i(x) = \lim_{m \to \infty} f_{i \cdot m}(x)$$

and the theorem [4]-2.6.5, it follows that for given

(4.2.23) $$\delta_1 \in R, \ \delta_1 > 0$$

there exists $M_i$ such that the condition $m > M_i$ implies that

(4.2.24) $$\|f_{i \cdot m}(x) - f_i(x)\| < \delta_1$$

Let

$$M = \max(M_1, ..., M_n)$$

For $k, m > M$, the inequality

(4.2.25) $$\begin{aligned}\|f_{i \cdot k}(x) - f_{i \cdot m}(x)\| &= \|f_{i \cdot k}(x) - f_i(x) + f_i(x) - f_{i \cdot m}(x)\| \\ &\leq \|f_{i \cdot k}(x) - f_i(x)\| + \|f_i(x) - f_{i \cdot m}(x)\| < 2\delta_1\end{aligned}$$

follows from the inequality (4.2.24) and the statement [4]-2.5.14.3. For $m > M$, the inequality

(4.2.26) $$\|f_{i \cdot m}(x)\| < F_i + \delta_1$$

follows from [4]-(2.1.1), (4.2.21), (4.2.24) and the inequality

$$\|f_{i \cdot m}(x)\| - \|f_i(x)\| \leq \|f_{i \cdot m}(x) - f_i(x)\| < \delta_1$$

According to the theorem 4.1.6, for any $m > M$, there exists integral

(4.2.27) $$\int_X d\mu(x)(f_{1 \cdot m}(x)...f_{n \cdot m}(x)\omega)$$

According to the statement [4]-2.1.9.3,

$$\left\| \int_X d\mu(x)(f_{1\cdot k}(x)...f_{n\cdot k}(x)\omega) - \int_X d\mu(x)(f_{1\cdot m}(x)...f_{n\cdot m}(x)\omega) \right\|$$

$$= \left\| \int_X d\mu(x)(f_{1\cdot k}(x)f_{2\cdot k}(x)...f_{n\cdot k}(x)\omega) \right.$$

$$- \int_X d\mu(x)(f_{1\cdot m}(x)f_{2\cdot k}(x)...f_{n\cdot k}(x)\omega)$$

$$+ \int_X d\mu(x)(f_{1\cdot m}(x)f_{2\cdot k}(x)...f_{n\cdot k}(x)\omega)$$

(4.2.28)
$$- ... - \left. \int_X d\mu(x)(f_{1\cdot m}(x)...f_{n\cdot m}(x)\omega) \right\|$$

$$\le \left\| \int_X d\mu(x)(f_{1\cdot k}(x)f_{2\cdot k}(x)...f_{n\cdot k}(x)\omega) \right.$$

$$- \left. \int_X d\mu(x)(f_{1\cdot m}(x)f_{2\cdot k}(x)...f_{n\cdot k}(x)\omega) \right\|$$

$$+ ... + \left\| \int_X d\mu(x)(f_{1\cdot m}(x)...f_{n-1\cdot m}(x)f_{n\cdot k}(x)\omega) \right.$$

$$- \left. \int_X d\mu(x)(f_{1\cdot m}(x)...f_{n-1\cdot m}(x)f_{n\cdot m}(x)\omega) \right\|$$

According to the definition [4]-2.1.3,

$$\left\| \int_X d\mu(x)(f_{1\cdot k}(x)...f_{n\cdot k}(x)\omega) - \int_X d\mu(x)(f_{1\cdot m}(x)...f_{n\cdot m}(x)\omega) \right\|$$

$$\le \left\| \int_X d\mu(x)(f_{1\cdot k}(x)f_{2\cdot k}(x)...f_{n\cdot k}(x)\omega - f_{1\cdot m}(x)f_{2\cdot k}(x)...f_{n\cdot k}(x)\omega) \right\| + ...$$

(4.2.29)
$$+ \left\| \int_X d\mu(x)(f_{1\cdot m}(x)...f_{n-1\cdot m}(x)f_{n\cdot k}(x)\omega \right.$$

$$- \left. f_{1\cdot m}(x)...f_{n-1\cdot m}(x)f_{n\cdot m}(x)\omega) \right\|$$

$$= \left\| \int_X d\mu(x)((f_{1\cdot k}(x) - f_{1\cdot m}(x))f_{2\cdot k}(x)...f_{n\cdot k}(x)\omega) \right\| + ...$$

$$+ \left\| \int_X d\mu(x)(f_{1\cdot m}(x)...f_{n-1\cdot m}(x)(f_{n\cdot k}(x) - f_{n\cdot m}(x))\omega) \right\|$$

follows from (4.2.9), (4.2.28). According to the theorem 4.1.6, for any $m > M$, the inequality

$$\left\| \int_X d\mu(x)(f_{1\cdot k}(x)...f_{n\cdot k}(x)\omega) - \int_X d\mu(x)(f_{1\cdot m}(x)...f_{n\cdot m}(x)\omega) \right\|$$

(4.2.30)
$$\le \int_X d\mu(x)(\|\omega\| \|f_{1\cdot k}(x) - f_{1\cdot m}(x)\| \|f_{2\cdot k}(x)\|...\|f_{n\cdot k}(x)\|) + ...$$

$$+ \int_X d\mu(x)(\|\omega\| \|f_{1\cdot m}(x)\|...\|f_{n-1\cdot m}(x)\|\|f_{n\cdot k}(x) - f_{n\cdot m}(x)\|)$$

$$\le 2\mu(x)\|\omega\|\delta_1((F_2 + \delta_1)...(F_n + \delta_1) + ... + (F_1 + \delta_1)...(F_{n-1} + \delta_1))$$

follows from inequalities (4.2.25), (4.2.26), (4.2.29). Let

$$(4.2.31) \quad \epsilon_1 = 2\mu(x)\|\omega\|\delta_1((F_2 + \delta_1)...(F_n + \delta_1) + ... + (F_1 + \delta_1)...(F_{n-1} + \delta_1))$$

From statements (4.2.22), (4.2.23) it follows that

$$(4.2.32) \quad \epsilon_1 > 0 \quad \frac{d\epsilon_1}{d\delta_1} > 0$$

From the equation (4.2.31) and the statement (4.2.32), it follows that $\epsilon_1$ is polynomial strictly monotone increasing function of $\delta_1$ such that

$$\delta_1 = 0 \Rightarrow \epsilon_1 = 0$$

According to the theorem [4]-2.3.5, the map (4.2.31) maps the interval $[0, \delta_1)$ into the interval $[0, \epsilon_1)$. According to the theorem [4]-2.3.3, for given $\epsilon > 0$ there exist $\delta > 0$ such that

$$\epsilon_1(\delta) < \epsilon$$

According to construction, a value of $M$ depends on a value of $\delta_1$. We choose the value of $M$ corresponding to $\delta_1 = \delta$. Therefore, for given $\epsilon \in R$, $\epsilon > 0$, there exists $M$ such that the condition $m > M$ implies that

$$(4.2.33) \quad \left\| \int_X d\mu(x)(f_{1 \cdot k}(x)...f_{n \cdot k}(x)\omega) - \int_X d\mu(x)(f_{1 \cdot m}(x)...f_{n \cdot m}(x)\omega) \right\| < \epsilon$$

From the equation (4.2.33), it follows that the sequence of integrals (4.2.27) is fundamental sequence. Therefore, the sequence of integrals (4.2.27) has limit

$$(4.2.34) \quad \int_X d\mu(x)(f_1(x)...f_n(x)\omega) = \lim_{m \to \infty} \int_X d\mu(x)(f_{1 \cdot m}(x)...f_{n \cdot m}(x)\omega)$$

According to the theorem 4.1.6, for any $m > M$,

$$(4.2.35) \quad \left\| \int_X d\mu(x)(f_{1 \cdot m}(x)...f_{n \cdot m}(x)\omega) \right\| \leq \int_X d\mu(x)(\|\omega\|\|f_{1 \cdot m}(x)\|...\|f_{n \cdot m}(x)\|)$$

The inequality (4.2.20) follows from the inequality (4.2.35) using passage to limit $m \to \infty$. □

THEOREM 4.2.6. *Let $\mu$ be a $\sigma$-additive measure defined on the set $X$. Let*

$$f : A_1 \relbar\joinrel\twoheadrightarrow A_2$$

*be representation of $\Omega_1$-group $A_1$ with norm $\|x\|_1$ in $\Omega_2$-group $A_2$ with norm $\|x\|_2$. Let*

$$g_i : X \to A_i \quad i = 1, 2$$

*be integrable map with compact range. Then map*

$$h = f_X(g_1)(g_2)$$

*is integrable map and*

$$(4.2.36) \quad \left\| \int_X d\mu(x)h(x) \right\|_2 \leq \int_X d\mu(x)(\|f\|\|g_1(x)\|_1\|g_2(x)\|_2)$$

PROOF. Since the range of the map $g_i$ is compact set, then, according to the theorem [4]-2.3.11, the following value is defined

(4.2.37) $$G_i = \sup \|g_i(x)\|_i$$

From the equation (4.2.37) and the statement [4]-2.1.9.1, it follows that

(4.2.38) $$G_i \geq 0$$

According to the theorem 3.3.2, for $i = 1$, $2$, there exists a sequence of simple integrable over the set $X$ maps

$$g_{i \cdot n} : X \to A_i$$

converging uniformly to $g_i$. From the equation

$$g_i(x) = \lim_{n \to \infty} g_{i \cdot n}(x)$$

and the theorem [4]-2.6.5, it follows that for given

(4.2.39) $$\delta_1 \in R,\ \delta_1 > 0$$

there exists $N_i$ such that the condition $n > N_i$ implies that

(4.2.40) $$\|g_{i \cdot n}(x) - g_i(x)\|_i < \delta_1$$

Let

$$N = \max(N_1, N_2)$$

For $k$, $n > N$, the inequality

(4.2.41) $$\|g_{i \cdot k}(x) - g_{i \cdot n}(x)\|_i = \|g_{i \cdot k}(x) - g_i(x) + g_i(x) - g_{i \cdot n}(x)\|_i$$
$$\leq \|g_{i \cdot k}(x) - g_i(x)\|_i + \|g_i(x) - g_{i \cdot n}(x)\|_i < 2\delta_1$$

follows from the inequality (4.2.40) and the statement [4]-2.5.14.3. For $n > N$, the inequality

(4.2.42) $$\|g_{i \cdot n}(x)\|_i < G_i + \delta_1$$

follows from [4]-(2.1.1), (4.2.37), (4.2.40) and the inequality

$$\|g_{i \cdot n}(x)\|_i - \|g_i(x)\|_i \leq \|g_{i \cdot n}(x) - g_i(x)\|_i < \delta_1$$

According to the theorem 4.1.7, for any $n > N$, there exists integral

(4.2.43) $$\int_X d\mu(x)(f(g_{1 \cdot n}(x))(g_{2 \cdot n}(x)))$$

According to the statement [4]-2.1.9.3,

(4.2.44) $$\left\| \int_X d\mu(x)(f(g_{1 \cdot k}(x))(g_{2 \cdot k}(x))) - \int_X d\mu(x)(f(g_{1 \cdot n}(x))(g_{2 \cdot n}(x))) \right\|_2$$
$$= \left\| \int_X d\mu(x)(f(g_{1 \cdot k}(x))(g_{2 \cdot k}(x))) - \int_X d\mu(x)(f(g_{1 \cdot n}(x))(g_{2 \cdot k}(x))) \right.$$
$$\left. + \int_X d\mu(x)(f(g_{1 \cdot n}(x))(g_{2 \cdot k}(x))) - \int_X d\mu(x)(f(g_{1 \cdot n}(x))(g_{2 \cdot n}(x))) \right\|_2$$
$$\leq \left\| \int_X d\mu(x)(f(g_{1 \cdot k}(x))(g_{2 \cdot k}(x))) - \int_X d\mu(x)(f(g_{1 \cdot n}(x))(g_{2 \cdot k}(x))) \right\|_2$$
$$+ \left\| \int_X d\mu(x)(f(g_{1 \cdot n}(x))(g_{2 \cdot k}(x))) - \int_X d\mu(x)(f(g_{1 \cdot n}(x))(g_{2 \cdot n}(x))) \right\|_2$$

According to the definition [2]-2.1.4, the map $f$ is homomorphism of Abelian group $A_1$ and, for any $a_1 \in A_1$, a map $f(a_1)$ is homomorphism of Abelian group $A_2$. Therefore,

$$\left\| \int_X d\mu(x)(f(g_{1\cdot k}(x))(g_{2\cdot k}(x))) - \int_X d\mu(x)(f(g_{1\cdot n}(x))(g_{2\cdot n}(x))) \right\|_2$$
$$\leq \left\| \int_X d\mu(x)(f(g_{1\cdot k}(x))(g_{2\cdot k}(x)) - f(g_{1\cdot n}(x))(g_{2\cdot k}(x))) \right\|_2$$
(4.2.45)
$$+ \left\| \int_X d\mu(x)(f(g_{1\cdot n}(x))(g_{2\cdot k}(x))) - f(g_{1\cdot n}(x))(g_{2\cdot n}(x))) \right\|_2$$
$$= \left\| \int_X d\mu(x)(f(g_{1\cdot k}(x) - g_{1\cdot n}(x))g_{2\cdot k}(x))) \right\|_2$$
$$+ \left\| \int_X d\mu(x)(f(g_{1\cdot n}(x))(g_{2\cdot k}(x) - g_{2\cdot n}(x))) \right\|_2$$

follows from (4.2.9), (4.2.44). According to the theorem 4.1.7, for any $n > N$, the inequality

$$\left\| \int_X d\mu(x)(f(g_{1\cdot k}(x))(g_{2\cdot k}(x))) - \int_X d\mu(x)(f(g_{1\cdot n}(x))(g_{2\cdot n}(x))) \right\|_2$$
(4.2.46)
$$\leq \int_X d\mu(x)(\|f\| \, \|g_{1\cdot k}(x) - g_{1\cdot n}(x)\|_1 \, \|g_{2\cdot k}(x)\|_2)$$
$$+ \int_X d\mu(x)(\|f\| \, \|g_{1\cdot n}(x)\|_1 \, \|g_{2\cdot k}(x) - f_{2\cdot n}(x)\|\|_2)$$
$$\leq 2\mu(x)\|f\|\delta_1((G_2 + \delta_1) + (G_1 + \delta_1))$$

follows from inequalities (4.2.41), (4.2.42), (4.2.45). Let

(4.2.47) $\quad\quad\quad \epsilon_1 = 2\mu(x)\|f\|\delta_1((G_2 + \delta_1) + (G_1 + \delta_1))$

From statements (4.2.38), (4.2.39) it follows that

(4.2.48) $\quad\quad\quad\quad\quad \epsilon_1 > 0 \quad \dfrac{d\epsilon_1}{d\delta_1} > 0$

From the equation (4.2.47) and the statement (4.2.48), it follows that $\epsilon_1$ is polynomial strictly monotone increasing function of $\delta_1$ such that

$$\delta_1 = 0 \Rightarrow \epsilon_1 = 0$$

According to the theorem [4]-2.3.5, the map (4.2.47) maps the interval $[0, \delta_1)$ into the interval $[0, \epsilon_1)$. According to the theorem [4]-2.3.3, for given $\epsilon > 0$ there exist $\delta > 0$ such that

$$\epsilon_1(\delta) < \epsilon$$

According to construction, a value of $N$ depends on a value of $\delta_1$. We choose the value of $N$ corresponding to $\delta_1 = \delta$. Therefore, for given $\epsilon \in R$, $\epsilon > 0$, there exists $N$ such that the condition $n > N$ implies that

(4.2.49) $\quad \left\| \int_X d\mu(x)(f(g_{1\cdot k}(x))(g_{2\cdot k}(x))) - \int_X d\mu(x)(f(g_{1\cdot n}(x))(g_{2\cdot n}(x))) \right\|_2 < \epsilon$

From the equation (4.2.49), it follows that the sequence of integrals (4.2.43) is fundamental sequence. Therefore, the sequence of integrals (4.2.43) has limit

$$(4.2.50) \quad \int_X d\mu(x)(f(g_1(x))(g_2(x))) = \lim_{n \to \infty} \int_X d\mu(x)(f(g_{1 \cdot n}(x))(g_{2 \cdot n}(x)))$$

According to the theorem 4.1.7, for any $n > N$,

$$(4.2.51) \quad \left\| \int_X d\mu(x)(f(g_{1 \cdot n}(x))(g_{2 \cdot n}(x))) \right\|_2 \leq \int_X d\mu(x)(\|f\| \|g_{1 \cdot n}(x)\|_1 \|g_{2 \cdot n}(x)\|_2)$$

The inequality (4.2.36) follows from the inequality (4.2.51) using passage to limit $n \to \infty$. □

THEOREM 4.2.7. *Let $\mu$ be a $\sigma$-additive measure defined on the set $X$. Let*

$$f : A_1 \relbar\joinrel\twoheadrightarrow A_2$$

*be effective representation of $\Omega_1$-group $A_1$ with norm $\|x\|_1$ in $\Omega_2$-group $A_2$ with norm $\|x\|_2$. Let transformations of the representation $f$ be morphism of representation*

$$R \relbar\joinrel\twoheadrightarrow A_2$$

*Let*

$$g_2 : X \to A_2$$

*be integrable map into normed $\Omega$-group $A_2$. Then map*

$$h = a_1 g_2$$

*is integrable map and*

$$\int_X d\mu(x)(a_1 g_2(x)) = a_1 \int_X d\mu(x) g_2(x)$$

PROOF. The theorem follows from the theorem 4.1.8 and from the definition 4.2.1. □

THEOREM 4.2.8. *Let $\mu$ be a $\sigma$-additive measure defined on the set $X$. Let*

$$f : A_1 \relbar\joinrel\twoheadrightarrow A_2$$

*be effective representation of $\Omega_1$-group $A_1$ with norm $\|x\|_1$ in $\Omega_2$-group $A_2$ with norm $\|x\|_2$. Let transformations of the representation*

$$R \relbar\joinrel\twoheadrightarrow A_1$$

*be morphism of representation $f$. Let*

$$g_1 : X \to A_1$$

*be integrable map into normed $\Omega$-group $A_1$. Then map*

$$h = g_1 a_2$$

*is integrable map and*

$$\int_X d\mu(x)(g_1(x) a_2) = \left( \int_X d\mu(x) g_1(x) \right) a_2$$

PROOF. The theorem follows from the theorem 4.1.9 and from the definition 4.2.1. □

Theorem 4.2.9. *Let*
$$f : X \to A$$
*be $\mu$-measurable map such that*
$$\|f(x)\| \le M$$
*Since the measure of the set $X$ is finite, then*

(4.2.52) $$\int_X d\mu(x)\|f(x)\| \le M\mu(X)$$

Proof. The theorem follows from the theorem 4.1.10 and from the definition 4.2.1. □

Theorem 4.2.10. *Let*
$$f : X \to A$$
*be $\mu$-measurable map such that*
$$\|f(x)\| \le M$$
*Since the measure of the set $X$ is finite, then map $f$ is integrable and*

(4.2.53) $$\left\|\int_X d\mu(x) f(x)\right\| \le M\mu(X)$$

Proof. Integrability of the map $f$ follows from theorems 4.2.4, 4.2.9. The inequality (4.2.53) follows from inequalities (4.2.19), (4.2.52). □

### 4.3. Lebesgue Integral as Map of the Set

Consider $\mu$-measurable map
$$f : X \to A$$
into $\Omega$-group $A$. Let $\mathcal{C}_X$ be $\sigma$-algebra of measurable sets of set $X$. Since $Y \in \mathcal{C}_X$, then we consider the expression

(4.3.1) $$F(Y) = \int_Y d\mu(x) f(x)$$

as map
$$F : \mathcal{C}_X \to A$$

Lemma 4.3.1. *Let*
$$f : X \to A$$
*be simple map. Let*

(4.3.2) $$X = \bigcup_i X_i \quad i \ne j => X_i \cap X_j = \emptyset$$

*be finite or countable union of sets $X_i$. The map $f$ is integrable map on the set $X$ iff the map $f$ is integrable map on each $X_i$*

(4.3.3) $$\int_X d\mu(x) f(x) = \sum_i \int_{X_i} d\mu(x) f(x)$$

*where series on the right converges normally.*

PROOF. Let the set $\{f_1, f_2, ...\}$ be domain of the map $f$. Let

(4.3.4) $$F_n = \{x \in X : f(x) = f_n\}$$

(4.3.5) $$F_{in} = \{x \in X_i : f(x) = f_n\}$$

Equations

(4.3.6) $$F_{in} = F_n \cap X_i$$

(4.3.7) $$F_n = \bigcup_i F_{in}$$

follow from equations (4.3.2), (4.3.4), (4.3.5). The equation

(4.3.8) $$F_{in} \cap F_{im} = \emptyset$$

follows from equations (4.3.2), (4.3.6). The equation

(4.3.9) $$\mu(F_n) = \sum_n \mu(F_{in})$$

follows from equations (4.3.7), (4.3.8) and from the statement 2.2.10.2. The equation

(4.3.10) $$\int_X d\mu(x) f(x) = \sum_n \mu(F_n) f_n = \sum_n \left( \sum_i \mu(F_{in}) \right) f_n = \sum_n \sum_i \mu(F_{in}) f_n$$

follows from equations (4.1.1), (4.3.9). The inequality

(4.3.11) $$\|\mu(F_{in}) f_n\| \leq \|\mu(F_n) f_n\| \leq \sum_i \|\mu(F_{in}) f_n\|$$

follows from equations (4.3.9), (4.3.6) and the statement [4]-2.1.9.3. The inequality

(4.3.12) $$\sum_n \|\mu(F_{in}) f_i\| \leq \sum_n \|\mu(F_n) f_n\| \leq \sum_n \sum_i \|\mu(F_{in}) f_n\|$$

follows from the inequality (4.3.11). Normal convergence of series $\sum_n \sum_i \|\mu(F_{in}) f_n\|$ means

(4.3.13) $$\sum_n \sum_i \|\mu(F_{in}) f_n\| = \sum_i \sum_n \|\mu(F_{in}) f_n\|$$

(4.3.14) $$\sum_n \sum_i \mu(F_{in}) f_n = \sum_i \sum_n \mu(F_{in}) f_n$$

From (4.3.12), (4.3.13) it follows that series $\sum_n \mu(F_n) f_n$ converges normally iff for each $i$ series $\sum_n \mu(F_{in}) f_n$ converges normally. According to the definition 4.1.2 for each $i$, there exists integral

(4.3.15) $$\int_{X_i} d\mu(x) f(x) = \sum_n \mu(F_{in}) f_n$$

Therefore, the map $f$ is integrable map on the set $X$ iff the map $f$ is integrable map on each $X_i$. The equation (4.3.3) follows from equations (4.3.10), (4.3.14), (4.3.15). □

THEOREM 4.3.2 ($\sigma$-additivity of Lebesgue integral). *Let*
$$f : X \to A$$
*be $\mu$-measurable map. Let*

(4.3.16) $$X = \bigcup_i X_i \quad i \neq j => X_i \cap X_j = \emptyset$$

*be finite or countable union of sets $X_i$. The map $f$ is integrable map on the set $X$ iff the map $f$ is integrable map on each $X_i$*

(4.3.17) $$\int_X d\mu(x) f(x) = \sum_i \int_{X_i} d\mu(x) f(x)$$

*where series on the right converges normally.*

PROOF. According to the definition 4.2.1, for any $\epsilon \in R$, $\epsilon > 0$, there exist a simple map
$$g : X \to A$$
integrable on $X$ and such that for any $x \in X$

(4.3.18) $$\|f(x) - g(x)\| < \frac{\epsilon}{2\mu(X)}$$

According to the lemma 4.3.1,

(4.3.19) $$\int_X d\mu(x) g(x) = \sum_i \int_{X_i} d\mu(x) g(x)$$

where $g$ is integrable on each $X_i$ and series (4.3.19) converges. According to the theorem 4.2.2, the map $f$ is integrable on each $X_i$. According to the theorem 4.2.10,

(4.3.20) $$\left\| \int_X d\mu(x) f(x) - \int_X d\mu(x) g(x) \right\| = \left\| \int_X d\mu(x)(f(x) - g(x)) \right\|$$
$$< \frac{\epsilon}{2\mu(X)} \mu(X) = \frac{\epsilon}{2}$$

(4.3.21) $$\left\| \int_{X_i} d\mu(x) f(x) - \int_{X_i} d\mu(x) g(x) \right\| = \left\| \int_{X_i} d\mu(x)(f(x) - g(x)) \right\|$$
$$< \frac{\epsilon}{2\mu(X)} \mu(X_i)$$

follow from (4.3.18). According to statements [4]-2.1.9.3, 2.2.10.2,

(4.3.22) $$\left\| \sum_i \int_{X_i} d\mu(x) f(x) - \int_X d\mu(x) g(x) \right\|$$
$$= \left\| \sum_i \int_{X_i} d\mu(x) f(x) - \sum_i \int_{X_i} d\mu(x) g(x) \right\|$$
$$= \left\| \sum_i \left( \int_{X_i} d\mu(x) f(x) - \int_{X_i} d\mu(x) g(x) \right) \right\|$$
$$< \sum_i \left\| \int_{X_i} d\mu(x)(f(x) - g(x)) \right\|$$
$$< \frac{\epsilon}{2\mu(X)} \sum_i \mu(X_i) = \frac{\epsilon}{2}$$

follows from (4.2.9), (4.3.21). The inequality

$$(4.3.23) \quad \begin{aligned} 0 &\leq \left\| \int_X d\mu(x)f(x) - \sum_i \int_{X_i} d\mu(x)f(x) \right\| \\ &= \left\| \int_X d\mu(x)f(x) - \int_X d\mu(x)g(x) \right. \\ &\quad + \left. \int_X d\mu(x)g(x) - \sum_i \int_{X_i} d\mu(x)f(x) \right\| \\ &\leq \left\| \int_X d\mu(x)f(x) - \int_X d\mu(x)g(x) \right\| \\ &\quad + \left\| \int_X d\mu(x)g(x) - \sum_i \int_{X_i} d\mu(x)f(x) \right\| \\ &< \epsilon \end{aligned}$$

follows from (4.3.20), (4.3.22) and from the statement [4]-2.1.9.2. Since $\epsilon$ is arbitrary, then the equation (4.3.17) follows from the inequality (4.3.23). □

COROLLARY 4.3.3. *Since $\mu$-measurable map*

$$f : X \to A$$

*is integrable map on the set $X$, then map $f$ is integrable map on $\mu$-measurable set $X' \subset X$.* □

THEOREM 4.3.4. *Let*

$$f : X \to A$$
$$g : X \to A$$

*be $\mu$-measurable maps such that*

$$(4.3.24) \qquad \|f(x)\| \leq M$$

$$(4.3.25) \qquad \|g(x)\| \leq M$$

*Since $f(x) = g(x)$ almost everywhere, then*

$$(4.3.26) \qquad \int_X d\mu(x)f(x) = \int_X d\mu(x)g(x)$$

PROOF. From equations (4.3.24), (4.3.25) and the theorem 4.2.10, it follows that integrals

$$\int_X d\mu(x)f(x)$$
$$\int_X d\mu(x)g(x)$$

exist. Let

$$(4.3.27) \qquad X_1 = \{x \in X : f(x) = g(x)\}$$
$$X_2 = \{x \in X : f(x) \neq g(x)\}$$

Since $X = X_1 \cup X_2$, $X_1 \cap X_2 = \emptyset$, then, according to the theorem 4.3.2,

$$(4.3.28) \qquad \int_X d\mu(x)f(x) = \int_{X_1} d\mu(x)f(x) + \int_{X_2} d\mu(x)f(x)$$

$$(4.3.29) \quad \int_X d\mu(x)g(x) = \int_{X_1} d\mu(x)g(x) + \int_{X_2} d\mu(x)g(x)$$

According to the condition of the theorem

$$(4.3.30) \quad \mu(X_2) = 0$$

From the theorem 4.2.10, the statement [4]-2.1.9.1 and the equation (4.3.30), it follows that

$$(4.3.31) \quad 0 \leq \left\| \int_{X_2} d\mu(x)f(x) \right\| \leq M\mu(X_2) = 0$$

$$(4.3.32) \quad 0 \leq \left\| \int_{X_2} d\mu(x)g(x) \right\| \leq M\mu(X_2) = 0$$

The equation

$$(4.3.33) \quad \int_{X_2} d\mu(x)f(x) = \int_{X_2} d\mu(x)g(x) = 0$$

follows from (4.3.31), (4.3.32) and the statement [4]-2.1.9.2. The equation

$$(4.3.34) \quad \int_{X_1} d\mu(x)f(x) = \int_{X_1} d\mu(x)g(x)$$

follows from (4.3.27). The equation (4.3.26) follows from (4.3.28), (4.3.29), (4.3.33), (4.3.34). $\square$

THEOREM 4.3.5 (Chebyshev's inequality). *Since $c > 0$, then*[4.10]

$$(4.3.35) \quad \mu\{x \in X : \|f(x)\| \geq c\} \leq \frac{1}{c} \int_X d\mu(x)\|f(x)\|$$

PROOF. Let

$$(4.3.36) \quad X' = \{x \in X : \|f(x)\| \geq c\}$$

Then

$$(4.3.37) \quad \begin{aligned} \int_X d\mu(x)\|f(x)\| &= \int_{X'} d\mu(x)\|f(x)\| + \int_{X-X'} d\mu(x)\|f(x)\| \\ &\geq \int_{X'} d\mu(x)\|f(x)\| \geq c\mu(X') \end{aligned}$$

(4.3.35) follows from (4.3.36), (4.3.37). $\square$

THEOREM 4.3.6. *If*

$$(4.3.38) \quad \int_X d\mu(x)\|f(x)\| = 0$$

*then $f(x) = 0$ almost everywhere.*

PROOF. From (4.3.35), (4.3.38), it follows that

$$(4.3.39) \quad \mu\left\{x \in X : \|f(x)\| \geq \frac{1}{n}\right\} \leq n \int_X d\mu(x)\|f(x)\| = 0$$

for $n = 1, 2, \ldots$ . The equality

$$\mu\{x \in X : f(x) \geq 0\} = \lim_{n\to\infty} \left\{x \in X : \|f(x)\| \geq \frac{1}{n}\right\} = 0$$

---

[4.10]See also the theorem [1]-5, pages 299, 300.

follows from (4.3.39). □

THEOREM 4.3.7 (Normal continuity of Lebesgue integral). *If*[4.11]

$$f : X \to A$$

*is a map integrable on $X$, then for any $\epsilon \in R$, $\epsilon > 0$, there is $\delta > 0$ such that*

(4.3.40) $$\left\| \int_E d\mu(x) f(x) \right\| < \epsilon$$

*for any $E \in \mathcal{C}_X$ such that $\mu(E) < \delta$.*

PROOF. According to the theorem 4.2.4, integral

$$\int_X d\mu(x) f(x)$$

exists iff integral

$$\int_X d\mu(x) \|f(x)\|$$

exists. Let

(4.3.41) $$X_i = \{x \in X : i \leq \|f(x)\| < i+1\}$$

(4.3.42) $$Y_N = \bigcup_{i=1}^{N} X_i$$

$$Z_N = X \setminus Y_N$$

According to the theorem 4.3.2,

(4.3.43) $$\int_X d\mu(x) \|f(x)\| = \sum_{i=0}^{\infty} \int_{X_i} d\mu(x) \|f(x)\|$$

Since the series in (4.3.43) normally converges, then there exists $N$ such that

(4.3.44) $$\sum_{i=N+1}^{\infty} \int_{X_i} d\mu(x) \|f(x)\| = \int_{Z_N} d\mu(x) \|f(x)\| < \frac{\epsilon}{2}$$

Let

(4.3.45) $$0 < \delta < \frac{\epsilon}{2(N+1)}$$

Let $\mu(E) < \delta$. Then

(4.3.46) $$\left\| \int_E d\mu(x) f(x) \right\| \leq \int_E d\mu(x) \|f(x)\|$$
$$= \int_{E \cap Y_N} d\mu(x) \|f(x)\| + \int_{E \cap Z_N} d\mu(x) \|f(x)\|$$

According to the theorem 4.2.9, the inequality

(4.3.47) $$\int_{E \cap Y_N} d\mu(x) \|f(x)\| \leq (N+1)\mu(E) < (N+1)\delta = \frac{\epsilon}{2}$$

follows from (4.3.41), (4.3.42). The inequality

(4.3.48) $$\int_{E \cap Z_N} d\mu(x) \|f(x)\| \leq \int_{Z_N} d\mu(x) \|f(x)\| < \frac{\epsilon}{2}$$

---

[4.11]See also the theorem [1]-6, pages 300, 301.

follows from (4.3.46). The inequality (4.3.40) follows from inequalities (4.3.46), (4.3.47), (4.3.48). □

### 4.4. Passage to Limit in Lebesgue Integral

THEOREM 4.4.1. *Let the map*

$$g : X \to R$$

*be integrable map and for any $n \in N$ almost everywhere*

$$\|f(x)\| \leq g(x)$$

*Then the map $f$ is integrable map and*

(4.4.1) $$\left\| \int_X d\mu(x) f(x) \right\| \leq \int_X d\mu(x) g(x)$$

PROOF. According to the theorem 4.3.4,

$$\mu(X) = 0 \implies \int_X d\mu(x) f(x) = 0$$

Therefore, we may assume that condition of the theorem is true for any $x \in X$. The theorem follows from the theorems 4.2.4, [1]-3 on page 297. □

THEOREM 4.4.2 (Lebesgue's bounded convergence theorem). *Let*[4.12]

$$f_n : X \to A$$

*be sequence of $\mu$-measurable map into $\Omega$-group $A$ and*

$$f(x) = \lim_{n \to \infty} f_n(x)$$

*almost everywhere. Let the map*

$$g : X \to R$$

*be integrable map and for any $n \in N$ almost everywhere*

(4.4.2) $$\|f_n(x)\| \leq g(x)$$

*Then the map $f$ is integrable map and*

(4.4.3) $$\int_X d\mu(x) f(x) = \lim_{n \to \infty} \int_X d\mu(x) f_n(x)$$

PROOF. According to the theorem 4.3.4,

$$\mu(X) = 0 \implies \int_X d\mu(x) f(x) = 0$$

Therefore, we may assume that condition of the theorem is true for any $x \in X$. According to the theorem [5]-2 on page 56, the inequality

(4.4.4) $$\|f(x)\| \leq g(x)$$

follows from the inequality (4.4.2). According to the theorem 4.4.1, $f$ is integrable map. According to the theorem 4.3.7, for any $\epsilon \in R$, $\epsilon > 0$, there is $\delta > 0$ such that

(4.4.5) $$\int_Y d\mu(x) g(x) < \frac{\epsilon}{4}$$

---

[4.12] See also the theorem [1]-1 on pages 303, 304.

for any $Y \in \mathcal{C}_X$ such that

(4.4.6) $$\mu(Y) < \delta$$

According to the theorem 3.4.5, the set $Y$ satisfying the condition (4.4.6) can be chosen such way that the sequence $f_n$ converges uniformly on the set $Z = X \setminus Y$. According to the definition [4]-2.6.3, there exists $N$ such that $n > N$, $x \in Z$ imply

(4.4.7) $$\|f(x) - f_n(x)\| < \frac{\epsilon}{2\mu(C)}$$

The inequality

(4.4.8) $$\begin{aligned}&\left\|\int_X d\mu(x)f(x) - \int_X d\mu(x)f_n(x)\right\|\\ &\leq \left\|\int_Z d\mu(x)(f(x) - f_n(x))\right\| + \left\|\int_Y d\mu(x)f(x)\right\| - \left\|\int_Y d\mu(x)f_n(x)\right\|\\ &< \frac{\epsilon}{2} + \frac{\epsilon}{4} + \frac{\epsilon}{4} = \epsilon\end{aligned}$$

follows from the equation

$$\begin{aligned}&\int_X d\mu(x)f(x) - \int_X d\mu(x)f_n(x)\\ &= \int_Y d\mu(x)f(x) + \int_Z d\mu(x)f(x) - \int_Y d\mu(x)f_n(x) - \int_Z d\mu(x)f_n(x)\\ &= \int_Z d\mu(x)(f(x) - f_n(x)) + \int_Y d\mu(x)f(x) - \int_Y d\mu(x)f_n(x)\end{aligned}$$

from inequalities (4.4.2), (4.4.4), (4.4.7) and from the theorem 4.2.10. The equation (4.4.3) follows from the inequality (4.4.8) and from the definition [4]-2.1.17. □

COROLLARY 4.4.3. *If* $\|f_n(x)\| \leq M$ *and* $f_n \to f$, *then*

$$\int_X d\mu(x)f(x) = \lim_{n\to\infty} \int_X d\mu(x)f_n(x)$$

□

THEOREM 4.4.4 (Beppo Levi). *Let maps*[4.13]

(4.4.9) $$f_n : A \to R \quad i = 1, ...$$

*be integrable maps and*

(4.4.10) $$\int_A d\mu(x)f_n(x) \leq M \quad n = 1, ...$$

*for some constant $M$. For any $n$, let*

(4.4.11) $$f_n(x) \leq f_{n+1}(x)$$

*Then, almost everywhere on $A$, there exists finite limit*

(4.4.12) $$f(x) = \lim_{n\to\infty} f_n(x)$$

*The map $f$ is integrable map on $A$ and*

(4.4.13) $$\int_A d\mu(x)f(x) = \lim_{i\to\infty} \int_A d\mu(x)f_i(x)$$

---

[4.13]See also the theorem [1]-2 on page 305.

PROOF. We assume $f_1(x) \geq 0$, since in general we can set
$$f'_i(x) = f_i(x) - f_1(x)$$
According to this assumption and to the condition (4.4.11)

(4.4.14) $$f_n(x) \geq 0$$

From the statement (4.4.14) and from definitions 4.1.2, 4.2.1, it follows that for any measurable set $B \subseteq A$

(4.4.15) $$\int_B d\mu(x) f(x) \geq 0$$

From the theorem 4.3.2, it follows that for any measurable set $B \subseteq A$

(4.4.16) $$\int_A d\mu(x) f(x) = \int_B d\mu(x) f(x) + \int_{A \setminus B} d\mu(x) f(x)$$

The statement

(4.4.17) $$\int_B d\mu(x) f_n(x) \leq M \quad n = 1, ...$$

follows from (4.4.10), (4.4.15), (4.4.16).

Consider the set
$$\Omega = \{x \in A : f_n(x) \to \infty\}$$
Then
$$\Omega = \bigcap_r \bigcup_n \Omega_n(r)$$
where

(4.4.18) $$\Omega_n(r) = \{x \in A : f_n(x) > r\}$$

By Chebyshev's inequality (4.3.35),

(4.4.19) $$\mu(\Omega_n(r)) \leq \frac{M}{r}$$

follows from (4.4.10), (4.4.18). Since
$$\Omega_1(r) \subseteq \Omega_2(r) \subseteq ... \subseteq \Omega_n(r) \subseteq ...$$
then

(4.4.20) $$\mu\left(\bigcup_n \Omega_n(r)\right) \leq \frac{M}{r}$$

follows from (4.4.19). Since for any $r$
$$\Omega \subseteq \bigcup_n \Omega_n(r)$$
then

(4.4.21) $$\mu(\Omega) \leq \frac{M}{r}$$

Since $r$ is arbitrary, then $\mu(\Omega) = 0$ follows from (4.4.21). Therefore, monotonic sequence $f_n(x)$ has finite limit $f(x)$ almost everywhere on $A$.

Let

(4.4.22) $$A(r) = \{x \in A : r - 1 \leq f(x) < r\}$$

(4.4.23) $$\phi : A \to R \quad x \in A(r) => \phi(x) = r$$

Let

$$B(s) = \bigcup_{r=1}^{s} A(r) \qquad (4.4.24)$$

From (4.4.22), (4.4.24), it follows that maps $f_n$ and $f$ are bounded on $B(s)$. From (4.4.22), (4.4.23), it follows that

$$\phi(x) \leq f(x) + 1 \qquad (4.4.25)$$

According to theorems 2.2.6, 4.4.1,

$$\int_{B_s} d\mu(x)\phi(x) \leq \int_{B_s} d\mu(x)f(x) + \mu(A) \qquad (4.4.26)$$

follows from (4.4.25). According to corollary 4.4.3,

$$\int_{B_s} d\mu(x)\phi(x) \leq \lim_{n\to\infty} \int_{B_s} d\mu(x)f_n(x) + \mu(A) \leq M + \mu(A) \qquad (4.4.27)$$

follows from (4.4.17), (4.4.22), (4.4.24), (4.4.26). According to equations (4.4.22), (4.4.23) and the definition 3.2.1, the map $\phi$ is simple map. According to the definition 4.1.2,

$$\int_{B_s} d\mu(x)\phi(x) = \sum_{r=1}^{s} r\mu(A(r)) \qquad (4.4.28)$$

Convergence of series

$$\sum_{r=1}^{\infty} r\mu(A(r)) \leq M + \mu(A) < \infty$$

follows from (4.4.27), (4.4.28). According to the definition 3.2.1, the map $\phi$ is integrable

$$\int_{A} d\mu(x)\phi(x) = \sum_{r=1}^{\infty} r\mu(A(r)) \qquad (4.4.29)$$

Therefore, te theorem follows from the theorem 4.4.2, since

$$f_n(x) \leq f(x) \leq \phi(x)$$

follows from (4.4.11), (4.4.13), (4.4.22), (4.4.23). □

CHAPTER 5

# Fubini's Theorem

## 5.1. Product of semirings of sets

DEFINITION 5.1.1. Let $\mathcal{L}_i$, $i = 1, ..., n$, be system of subsets of set $A_i$. **Cartesian product of systems of subsets**
$$\mathcal{M} = \mathcal{L}_1 \times ... \times \mathcal{L}_n$$
is systems of subsets of set $A = A_1 \times ... \times A_n$ which can be represented as[5.1]
$$C = C_1 \times ... \times C_n \quad C_i \in \mathcal{L}_i$$
If $\mathcal{L}_1 = ... = \mathcal{L}_n = \mathcal{L}$, then the system of sets
$$\mathcal{M} = \mathcal{L}^n$$
is called **Cartesian power $n$ of systems of subsets** □

THEOREM 5.1.2. *Cartesian product*[5.2]
$$\mathcal{S} = \mathcal{S}_1 \times ...\mathcal{S}_n$$
*of semirings of sets $\mathcal{S}_1, ..., \mathcal{S}_n$ is semiring of sets.*

PROOF. We will prove the theorem in case $n = 2$. In general case, the proof is similar.

5.1.2.1: Let $A, B \in \mathcal{S}_1 \times \mathcal{S}_2$. Therefore

(5.1.1) $$A = A_1 \times A_2 \quad A_1 \in \mathcal{S}_1 \quad A_2 \in \mathcal{S}_2$$

(5.1.2) $$B = B_1 \times B_2 \quad B_1 \in \mathcal{S}_1 \quad B_2 \in \mathcal{S}_2$$

It is evident that
$$(x_1, x_2) \in (A_1 \times A_2) \cap (B_1 \times B_2)$$
iff
$$x_1 \in A_1 \cap B_1 \quad x_2 \in A_2 \cap B_2$$
Therefore

(5.1.3) $$(A_1 \times A_2) \cap (B_1 \times B_2) = (A_1 \cap B_1) \times (A_2 \cap B_2)$$

From (5.1.1), (5.1.2) and from the statement 2.1.1.2, it follows that

(5.1.4) $$A_1 \cap B_1 \in \mathcal{S}_1 \quad A_2 \cap B_2 \in \mathcal{S}_2$$

From (5.1.1), (5.1.2), (5.1.4), it follows that
$$(A_1 \times A_2) \cap (B_1 \times B_2) \in \mathcal{S}_1 \times \mathcal{S}_2$$

---

[5.1]See also the definition in [1] on page 353.
[5.2]See also the theorem [1]-1 on page 353.

Therefore, the statement 2.1.1.2 is true for $\mathcal{S}_1 \times \mathcal{S}_2$.

5.1.2.2: Let $B_1 \subset A_1$, $B_2 \subset A_2$. From the statement 2.1.1.3, it follows that

(5.1.5) $$A_1 = B_1 \cup B_{1\cdot 1} \cup ... \cup B_{1\cdot k_1} \quad B_{1\cdot i} \in \mathcal{S}_1$$

(5.1.6) $$A_2 = B_2 \cup B_{2\cdot 1} \cup ... \cup B_{2\cdot k_2} \quad B_{2\cdot i} \in \mathcal{S}_2$$

From (5.1.5), (5.1.6), it follows that

$$(x_1, x_2) \in A = A_1 \times A_2$$

iff $(x_1, x_2)$ belongs to one of the following sets: $B = B_1 \times B_2$, $B_1 \times B_{2\cdot j} \in \mathcal{S}_1 \times \mathcal{S}_2$, $B_{1\cdot i} \times B_2 \in \mathcal{S}_1 \times \mathcal{S}_2$, $B_{1\cdot i} \times B_{2\cdot j} \in \mathcal{S}_1 \times \mathcal{S}_2$. Therefore,

$$A = B \cup (B_1 \times B_{2\cdot 1}) \cup ... \cup (B_1 \times B_{2\cdot k_2})$$
$$\cup (B_{1\cdot 1} \times B_2) \cup (B_{1\cdot 1} \times B_{2\cdot 1}) \cup ... \cup (B_{1\cdot 1} \times B_{2\cdot k_2})$$
$$...$$
$$\cup (B_{k_1 \cdot 1} \times B_2) \cup (B_{k_1 \cdot 1} \times B_{2\cdot 1}) \cup ... \cup (B_{k_1 \cdot 1} \times B_{2\cdot k_2})$$

and the statement 2.1.1.3 is true for $\mathcal{S}_1 \times \mathcal{S}_2$.

From reasonings 5.1.2.1, 5.1.2.2, it follows that $\mathcal{S}_1 \times \mathcal{S}_2$ is semiring of sets. □

However, in general, Cartesian product of rings of sets is not ring of sets.

DEFINITION 5.1.3. The ring of sets

$$\mathcal{R} = \mathcal{R}_1 \otimes ... \otimes \mathcal{R}_n$$

which is generated by the semiring of sets

$$\mathcal{R}_1 \times ... \times \mathcal{R}_n$$

is called **product of rings of sets** $\mathcal{R}_1, ..., \mathcal{R}_n$. □

THEOREM 5.1.4. *Product of algebras of sets is algebra of sets.*

PROOF. For $i = 1, ..., n$, let the ring of sets $\mathcal{R}_i$ be algebra of sets. By the definition 2.1.2, algebra of sets $\mathcal{R}_i$ has unit $E_i$ such that

(5.1.7) $$A_i \cap E_i = A_i$$

for any $A_i \in \mathcal{R}_i$. From the reasoning 5.1.2.1 and from (5.1.7), it follows that

(5.1.8) $$(A_1 \times ... \times A_n) \cap (E_1 \times ... \times E_n) = A_1 \times ... \times A_n$$

From the theorem 2.1.11 and from the equation (5.1.8), it follows that for any $A \in \mathcal{R}_1 \otimes ... \otimes R_n$ the following equation is true

$$A \cap (E_1 \times ... \times E_n) = A$$

Therefore, the set $E_1 \times ... \times E_n$ is unit of the ring of sets $\mathcal{R}_1 \times ... \times \mathcal{R}_n$. By the definition 2.1.2, the ring of sets $\mathcal{R}_1 \times ... \times \mathcal{R}_n$ is the algebra of sets. □

## 5.2. Product of Measures

THEOREM 5.2.1. *Let*
$$\mu_i : \mathcal{R}_i \to R$$
*be measure on semiring* $\mathcal{R}_i$, $i = 1, ..., n$. **Cartesian product**
$$\mu = \mu_1 \times ... \times \mu_n$$
**of measures** *is defined by the formula*[5.3]

(5.2.1) $$\mu(A_1 \times ... \times A_n) = \mu_1(A_1)...\mu_n(A_n)$$

PROOF. We will prove the theorem in case $n = 2$. In general case, the proof is similar. Let

(5.2.2) $$A = A_1 \times A_2 = \bigcup_{i=1}^{t} B^i$$

where $i \neq j => B^i \cap B^j = \emptyset$ and $B^i = B_1^i \times B_2^i$. According to the lemma 2.1.10, there exist finite expansions
$$A_1 = \bigcup_{m=1}^{r} C_1^m \quad A_2 = \bigcup_{n=1}^{s} C_2^n$$
such that

(5.2.3) $$B_i^k = \bigcup_{m \in M_{ik}} C_i^m \quad i = 1, 2$$

$$M_{1k} \subseteq \{1, ..., r\} \quad M_{2k} \subseteq \{1, ..., s\}$$

According to the statement 2.2.1.2 with respect to measures $\mu_1, \mu_2$, , from (5.2.3), it follows that

$$\mu(A) = \mu_1(A_1)\mu_2(A_2) = \sum_{m=1}^{r} \mu(C_1^m) \sum_{n=1}^{s} \mu(C_2^n)$$

(5.2.4) $$= \sum_{k=1}^{t} \sum_{m \in M_{1k}} C_1^m \sum_{n \in M_{2k}} C_2^n$$

$$= \sum_{k=1}^{t} \mu_1(B_1^k)\mu_2(B_2^k) = \sum_{k=1}^{t} \mu_1(B^k)$$

From (5.2.2), (5.2.4), it follows that the statement 2.2.1.2 is true for $\mu$. □

THEOREM 5.2.2. *Let* $\mu_i$, $i = 1, ..., n$, *be* $\sigma$-*additive measure*[5.4] *defined on* $\sigma$-*algebra* $\mathcal{C}_i$. *Cartesian product of measures* $\mu_1 \times ... \times \mu_n$ *is* $\sigma$-*additive measure on semiring* $\mathcal{C}_1 \times ... \times \mathcal{C}_n$.

PROOF. We will prove the theorem in case $n = 2$. In general case, the proof is similar. Let

(5.2.5) $$C = \bigcup_{n=1}^{\infty} C_n \quad C, C_n \in \mathcal{C}_1 \times \mathcal{C}_2$$

---

[5.3]See [1], pages 354, 355, the proof, that the map defined by the formula (5.2.1), is measure.
[5.4]See also the theorem [1]-2 on page 355.

where

(5.2.6) $$n \neq m => C_n \cap C_m = \emptyset$$

By (5.2.5)

$$C = A \times B \quad A \in \mathcal{C}_1 \quad B \in \mathcal{C}_2$$

$$C_n = A_n \times B_n \quad A_n \in \mathcal{C}_1 \quad B_n \in \mathcal{C}_2$$

Consider the set of maps

$$f_n : X \to R$$

defined by the rule

(5.2.7) $$f_n(x) = \begin{cases} \mu_2(B_n) & x \in A_n \\ 0 & x \notin A_n \end{cases}$$

Let $x \in A$. Let

$$N_x = \{n : \exists y \in B, (x, y) \in A_n \times B_n\}$$

By (5.2.6), for any $y \in B$, there exists unique $n \in N_x$ such that

$$(x, y) \in A_n \times B_n \quad y \in B_n$$

and

$$n \in N_x, m \in N_x, n \neq m => B_n \cap B_m = \emptyset$$

$$\bigcup_{n \in N_x} B_n = B$$

According to the statement 2.2.10.2 with respect to measure $\mu_2$

(5.2.8) $$\sum_n f_n(x) = \mu_2(B)$$

According to theorems 4.2.3, 4.2.7 and to the corollary 4.4.3, the equation

(5.2.9) $$\sum_n \int_A d\mu_1(x) f_n(x) = \int_A d\mu_1(x) \mu_2(B) = \mu_1(A)\mu_2(B) = \mu(C)$$

follows from the equation (5.2.8). According to the theorem 4.2.7, the equation

(5.2.10) $$\int_A d\mu_1(x) f_n(x) = \mu_1(A_n)\mu_2(B_n) = \mu(C_n)$$

follows from the equation (5.2.7). The equation

$$\mu(C) = \sum_n \mu(C_n)$$

follows from equations (5.2.9) (5.2.10). Therefore, the statement 2.2.10.2 is true for the measure $\mu$. □

Lebesgue extension of measure $\mu_1 \times ... \times \mu_n$ defined on $\sigma$-algebra $\mathcal{C} \supset \mathcal{C}_1 \times ... \times \mathcal{C}_n$ is called **product of measures**

$$\mu_1 \otimes ... \otimes \mu_n = \bigotimes \mu_i$$

When $\mu_1 = ... = \mu_n = \mu$, then the product of measures is called **power of measure** $\mu$

$$\mu^n = \mu_1 \otimes ... \otimes \mu_n$$

THEOREM 5.2.3. Let $\mu_i$, $i = 1, 2$, be $\sigma$-additive measure[5.5] defined on $\sigma$-algebra $\mathcal{C}_i$ of subsets of the set $X_i$. Let $\mu = \mu_1 \otimes \mu_2$. Let $A \in \mathcal{C}_1 \otimes \mathcal{C}_2$.

Let for any $x_1 \in X_1$

$$(5.2.11) \qquad A_{2x_1} = \{x_2 \in X_2 : (x_1, x_2) \in A\} \in \mathcal{C}_2$$

If the map
$$x_1 \to \mu_2(A_{2x_1})$$
is $\mu_1$-integrable map, then

$$(5.2.12) \qquad \mu(A) = \int_{X_1} d\mu_1(x_1) \mu_2(A_{2x_1})$$

Let for any $x_2 \in X_2$
$$A_{1x_2} = \{x_1 \in X_1 : (x_1, x_2) \in A\} \in \mathcal{C}_1$$

If the map
$$x_2 \to \mu_1(A_{1x_2})$$
is $\mu_2$-integrable map, then

$$(5.2.13) \qquad \mu(A) = \int_{X_2} d\mu_2(x_2) \mu_1(A_{1x_2})$$

PROOF. We will prove the equation (5.2.12). A proof of the equation (5.2.13) is similar.

LEMMA 5.2.4. The equation (5.2.12) holds for the set of the form

$$(5.2.14) \qquad A = A_1 \times A_2 \quad A_1 \in \mathcal{C}_1 \quad A_2 \in \mathcal{C}_2$$

PROOF. According to the definition 5.2.1,

$$(5.2.15) \qquad \mu(A) = \mu_1(A_1)\mu_2(A_2)$$

By (5.2.11)

$$(5.2.16) \qquad A_{2x_1} = \begin{cases} A_2 & x_1 \in A_1 \\ \emptyset & x_1 \notin A_1 \end{cases}$$

From (5.2.16) and from the theorem 2.2.2, it follows that

$$(5.2.17) \qquad \mu_2(A_{2x_1}) = \begin{cases} \mu_2(A_2) & x_1 \in A_1 \\ 0 & x_1 \notin A_1 \end{cases}$$

From (5.2.14), (5.2.17) and from the definition 3.2.1, it follows that the map
$$x_1 \to \mu_2(A_{2x_1})$$
is simple map. According to the definition 4.1.2,

$$(5.2.18) \qquad \int_{X_1} d\mu_1(x_1) \mu_2(A_{2x_1}) = \mu_1(A_1)\mu_2(A_2)$$

The lemma follows from equations (5.2.15), (5.2.18). ⊙

---

[5.5]See also the theorem [1]-3 on pages 356, 357. I do not assume that $X_i$ is set of real numbers.

LEMMA 5.2.5. The equation (5.2.12) holds for the set of the form

$$(5.2.19) \qquad A = \bigcup_{i=1}^{n} A_{1i} \times A_{2i} \quad A_{1i} \in \mathcal{C}_1 \quad A_{2i} \in \mathcal{C}_2$$

PROOF. From the lemma 2.1.10, there exists finite system of sets $B_{11}$, ..., $B_{1t_1} \in \mathcal{C}_1$ such that

$$(5.2.20) \qquad i \neq j \implies B_{1i} \cap B_{1j} = \emptyset$$

$$(5.2.21) \qquad A_{1i} = \bigcup_{j \in M_{1i}} B_{1j}$$

where $M_{1i} \subset \{1, ..., t_1\}$. From the lemma 2.1.10, there exists finite system of sets $B_{21}$, ..., $B_{2t_2} \in \mathcal{C}_2$ such that

$$(5.2.22) \qquad i \neq j \implies B_{2i} \cap B_{2j} = \emptyset$$

$$(5.2.23) \qquad A_{2i} = \bigcup_{j \in M_{2i}} B_{2j}$$

where $M_{2i} \subset \{1, ..., t_2\}$. From (5.2.20), (5.2.22), it follows that

$$(5.2.24) \qquad i \neq j \land k \neq l \implies (B_{1i} \times B_{2k}) \cap (B_{1j} \times B_{2l}) = \emptyset$$

From equations (5.2.21), (5.2.23), it follows that for any $i$

$$(5.2.25) \qquad A_{1i} \times A_{2i} = \bigcup_{k \in M_{1i}} \bigcup_{l \in M_{2i}} B_{1k} \times B_{2l}$$

Let

$$(5.2.26) \qquad M_i = M_{1i} \times M_{2i}$$

The equation

$$(5.2.27) \qquad A_{1i} \times A_{2i} = \bigcup_{(k,l) \in M_i} B_{1k} \times B_{2l}$$

follows from (5.2.25), (5.2.26). Let

$$(5.2.28) \qquad M = \bigcup_{i} M_i$$

The equation

$$(5.2.29) \qquad A = \bigcup_{(k,l) \in M} B_{1k} \times B_{2l}$$

follows from (5.2.19), (5.2.27), (5.2.28). According to the statement 2.2.1.2 and to the definition 5.2.1,

$$(5.2.30) \qquad \mu(A) = \sum_{(k,l) \in M} \mu(B_{1k} \times B_{2l}) = \sum_{(k,l) \in M} \mu_1(B_{1k})\mu_2(B_{2l})$$

follows from (5.2.24), (5.2.29).

By (5.2.11) (5.2.29),

$$(5.2.31) \qquad A_{2x_1} = \begin{cases} \bigcup_{(k,l)\in M} B_{2l} & x_1 \in B_{1k} \\ \emptyset & x_1 \notin \bigcup_k B_{1k} \end{cases}$$

From the equation (5.2.31), from the statement 2.2.1.2 and from the theorem 2.2.2, it follows that

$$(5.2.32) \qquad \mu_2(A_{2x_1}) = \begin{cases} \sum_{(k,l)\in M} \mu_2(B_{2l}) & x_1 \in B_{1k} \\ 0 & x_1 \notin \bigcup_k B_{1k} \end{cases}$$

From (5.2.29), (5.2.32) and from the definition 3.2.1, it follows that the map

$$x_1 \to \mu_2(A_{2x_1})$$

is simple map. According to the definition 4.1.2,

$$(5.2.33) \qquad \int_{X_1} d\mu_1(x_1)\mu_2(A_{2x_1}) = \sum_{(k,l)\in M} \mu_1(B_{1k})\mu_2(B_{2l})$$

The lemma follows from equations (5.2.30), (5.2.33). ⊙

According to the theorem 2.3.13, for any set $A \in \mathcal{C}_{\mu_1 \otimes \mu_2}$, there exist set $B$ such that

$$(5.2.34) \qquad A \subseteq B$$

$$(5.2.35) \qquad \mu(A) = \mu(B)$$

$$(5.2.36) \qquad B = \bigcap_n B_n$$

$$(5.2.37) \qquad B_1 \supseteq B_2 \supseteq ... \supseteq B_n \supseteq ...$$

$$(5.2.38) \qquad B_n = \bigcup_k B_{nk}$$

$$(5.2.39) \qquad B_{nk} \in \mathcal{R}(\mathcal{C}_{\mu_1} \times \mathcal{C}_{\mu_2})$$

$$(5.2.40) \qquad \mu(B_{nk}) < \mu(A) + \frac{1}{n}$$

$$(5.2.41) \qquad B_{n1} \subseteq B_{n2} \subseteq ... \subseteq B_{nk} \subseteq ...$$

According to lemmas 5.2.4, 5.2.5, from the statement (5.2.39), it follows that the equation (5.2.12) holds for the sets $B_{nk}$

$$(5.2.42) \qquad \mu(B_{nk}) = \int_{X_1} d\mu_1(x_1)\mu_2(B_{nk \cdot 2x_1})$$

where

$$(5.2.43) \qquad B_{nk \cdot 2x_1} = \{x_2 \in X_2 : (x_1, x_2) \in B_{nk}\}$$

The statement

$$(5.2.44) \qquad B_{n1 \cdot 2x_1} \subseteq B_{n2 \cdot 2x_1} \subseteq ... \subseteq B_{nk \cdot 2x_1} \subseteq ...$$

follows from (5.2.41), (5.2.43). Consider the set of maps
$$f_{nk} : X_1 \to R \quad f_{nk}(x_1) = \mu_2(B_{nk \cdot 2x_1})$$
According to the theorem 2.2.6, the statement
$$(5.2.45) \qquad f_{nk}(x_1) \leq f_{nk+1}(x_1)$$
follows from the statement (5.2.44). The statement
$$(5.2.46) \qquad \int_{X_1} d\mu_1(x_1) f_{nk}(x_1) = \int_{X_1} d\mu_1(x_1) \mu_2(B_{nk \cdot 2x_1}) < \mu(A) + \frac{1}{n}$$
follows from statements (5.2.40), (5.2.42). According to the theorem 4.4.4,

- almost everywhere on $X_1$, there exists finite limit
$$(5.2.47) \qquad f_n(x_1) = \lim_{k \to \infty} f_{nk}(x_1) = \lim_{k \to \infty} \mu_2(B_{nk \cdot 2x_1})$$
- the map $f_n$ is integrable map on $X_1$ and
$$(5.2.48) \qquad \begin{aligned} \int_{X_1} d\mu_1(x_1) f_n(x_1) &= \lim_{k \to \infty} \int_{X_1} d\mu_1(x_1) f_{nk}(x_1) \\ &= \lim_{k \to \infty} \int_{X_1} d\mu_1(x_1) \mu_2(B_{nk \cdot 2x_1}) \end{aligned}$$

The statement
$$(5.2.49) \qquad \lim_{k \to \infty} \int_{X_1} d\mu_1(x_1) \mu_2(B_{nk \cdot 2x_1}) = \int_{X_1} d\mu_1(x_1) \lim_{k \to \infty} \mu_2(B_{nk \cdot 2x_1})$$
follows from statements (5.2.47), (5.2.48). The statement
$$(5.2.50) \qquad B_{n \cdot 2x_1} = \bigcup_k B_{nk \cdot 2x_1}$$
where
$$B_{n \cdot 2x_1} = \{x_2 \in X_2 : (x_1, x_2) \in B_n\}$$
follows from the chain of statements
$$x_2 \in B_{n \cdot 2x_1} <=> (x_1, x_2) \in B_n <=> \exists k, (x_1, x_2) \in B_{nk}$$
$$<=> \exists k, x_2 \in B_{nk \cdot 2x_1} <=> x_2 \in \bigcup_k B_{nk \cdot 2x_1}$$
According to the theorem 2.2.13,

- equations
$$(5.2.51) \qquad \lim_{k \to \infty} \mu_2(B_{nk \cdot 2x_1}) = \mu_2(B_{n \cdot 2x_1})$$

$$(5.2.52) \qquad \lim_{k \to \infty} \int_{X_1} d\mu_1(x_1) \mu_2(B_{nk \cdot 2x_1}) = \int_{X_1} d\mu_1(x_1) \mu_2(B_{n \cdot 2x_1})$$
  follow from statements (5.2.44), (5.2.49), (5.2.50).
- the equation
$$(5.2.53) \qquad \lim_{k \to \infty} \mu(B_{nk}) = \mu(B_n)$$
  follows from statements (5.2.38), (5.2.41).

From equations (5.2.42), (5.2.52), (5.2.53), it follows that the equation (5.2.12) holds for the sets $B_n$

$$(5.2.54) \qquad \mu(B_n) = \int_{X_1} d\mu_1(x_1) \mu_2(B_{n \cdot 2x_1})$$

The statement

$$(5.2.55) \qquad B_{1 \cdot 2x_1} \supseteq B_{2 \cdot 2x_1} \supseteq ... \supseteq B_{k \cdot 2x_1} \supseteq ...$$

follows from (5.2.37), (5.2.43). The equation

$$f_n : X_1 \to R \qquad f_n(x_1) = \mu_2(B_{n \cdot 2x_1})$$

follows from equations (5.2.47), (5.2.51). According to the theorem 2.2.6, the statement

$$(5.2.56) \qquad f_n(x_1) \geq f_{n+1}(x_1)$$

follows from the statement (5.2.55). According to the theorem 4.2.10, the statement

$$(5.2.57) \qquad \int_{X_1} d\mu_1(x_1) f_n(x_1) = \int_{X_1} d\mu_1(x_1) \mu_2(B_{n \cdot 2x_1}) > 0$$

follows from the statement 2.2.1.1. According to the theorem 4.4.4,

- almost everywhere on $X_1$, there exists finite limit

$$(5.2.58) \qquad f(x_1) = \lim_{n \to \infty} f_n(x_1) = \lim_{n \to \infty} \mu_2(B_{n \cdot 2x_1})$$

- the map $f$ is integrable map on $X_1$ and

$$(5.2.59) \qquad \begin{aligned} \int_{X_1} d\mu_1(x_1) f(x_1) &= \lim_{n \to \infty} \int_{X_1} d\mu_1(x_1) f_n(x_1) \\ &= \lim_{n \to \infty} \int_{X_1} d\mu_1(x_1) \mu_2(B_{n \cdot 2x_1}) \end{aligned}$$

The statement

$$(5.2.60) \qquad \lim_{n \to \infty} \int_{X_1} d\mu_1(x_1) \mu_2(B_{n \cdot 2x_1}) = \int_{X_1} d\mu_1(x_1) \lim_{n \to \infty} \mu_2(B_{n \cdot 2x_1})$$

follows from statements (5.2.58), (5.2.59). The statement

$$(5.2.61) \qquad B_{2x_1} = \bigcap_k B_{n \cdot 2x_1}$$

where

$$B_{2x_1} = \{x_2 \in X_2 : (x_1, x_2) \in B\}$$

follows from the chain of statements

$$x_2 \in B_{2x_1} <=> (x_1, x_2) \in B <=> \forall n, (x_1, x_2) \in B_n$$
$$<=> \forall n, x_2 \in B_{n \cdot 2x_1} <=> x_2 \in \bigcap_n B_{n \cdot 2x_1}$$

According to the theorem 2.2.12,

- the equation

$$(5.2.62) \qquad \lim_{n \to \infty} \int_{X_1} d\mu_1(x_1) \mu_2(B_{n \cdot 2x_1}) = \int_{X_1} d\mu_1(x_1) \mu_2(B_{2x_1})$$

follows from statements (5.2.55), (5.2.60), (5.2.61).

- the equation

$$(5.2.63) \qquad \lim_{n \to \infty} \mu(B_n) = \mu(B)$$

follows from statements (5.2.36), (5.2.37).

From equations (5.2.54), (5.2.62), (5.2.63), it follows that the equation (5.2.12) holds for the set $B$

$$(5.2.64) \qquad \mu(B) = \int_{X_1} d\mu_1(x_1) \mu_2(B_{2x_1})$$

LEMMA 5.2.6. Since

$$(5.2.65) \qquad \mu(A) = 0$$

then

$$(5.2.66) \qquad \int_{X_1} d\mu_1(x_1) \mu_2(A_{2x_1}) = 0$$

PROOF. The equation

$$(5.2.67) \qquad \int_{X_1} d\mu_1(x_1) \mu_2(B_{2x_1}) = 0$$

follows from equations (5.2.35), (5.2.65). According to the theorem 4.3.6 and to the statement 2.2.1.1, the equation

$$(5.2.68) \qquad \mu_2(B_{2x_1}) = 0$$

almost everywhere follows from the equation (5.2.67). The statement

$$(5.2.69) \qquad A_{2x_1} \subseteq B_{2x_1}$$

follows from the statement (5.2.34). According to the theorem 2.2.8, the equation

$$(5.2.70) \qquad \mu_2(A_{2x_1}) = 0$$

almost everywhere follows from (5.2.68), (5.2.69). According to the definition 4.1.2, the equation (5.2.66) follows from the equation (5.2.70). ⊙

Let $C = B \setminus A$. Then

$$(5.2.71) \qquad B = A \cup C \quad A \cap C = \emptyset$$

$$(5.2.72) \qquad B_{2x_1} = A_{2x_1} \cup C_{2x_1} \quad A_{2x_1} \cap C_{2x_1} = \emptyset$$

Equations

$$(5.2.73) \qquad \mu(B) = \mu(A) + \mu(C)$$

$$(5.2.74) \qquad \mu_2(B_{2x_1}) = \mu_2(A_{2x_1}) + \mu_2(C_{2x_1})$$

follow from equations (5.2.71), (5.2.72). The equation

$$(5.2.75) \qquad \mu(A) + \mu(C) = \int_{X_1} d\mu_1(x_1) \mu_2(A_{2x_1}) + \int_{X_1} d\mu_1(x_1) \mu_2(C_{2x_1})$$

follows from equations (5.2.64), (5.2.73), (5.2.74) and the theorem 4.2.3. The equation

$$(5.2.76) \qquad \mu(C) = 0$$

follows from equations (5.2.35), (5.2.73). According to the lemma 5.2.6, the equation

$$\int_{X_1} d\mu_1(x_1)\mu_2(C_{2x_1}) = 0 \tag{5.2.77}$$

follows from the equation (5.2.76). The equation (5.2.12) follows from equations (5.2.75), (5.2.76), (5.2.77). □

## 5.3. Fubini's Theorem

THEOREM 5.3.1 (Fubini). *Let $\mu_1$ and $\mu_2$ be $\sigma$-additive complete measures defined on $\sigma$-algebras.*[5.6] *Let $\mu = \mu_1 \otimes \mu_2$. Let $A$ be complete Abelian $\Omega$-group. Let the map*

$$f : X_1 \times X_2 \to A$$

*with compact range be $\mu$-integrable on the set $B \subseteq X_1 \times X_2$. Then*

$$\int_B d\mu(x_1, x_2) f(x_1, x_2) = \int_{X_1} d\mu(x_1) \int_{B_{2x_1}} d\mu(x_2) f(x_1, x_2) \tag{5.3.1}$$

$$\int_B d\mu(x_1, x_2) f(x_1, x_2) = \int_{X_2} d\mu(x_2) \int_{B_{1x_2}} d\mu(x_1) f(x_1, x_2) \tag{5.3.2}$$

PROOF. We will prove the equation (5.3.1). A proof of the equation (5.3.2) is similar.

LEMMA 5.3.2. *The equation (5.3.1) holds for the simple map $f$.*

PROOF. Let $y_1, y_2, ...$ be the range of the map $f$. Let

$$B_i = \{(x_1, x_2) \in B : f(x_1, x_2) = y_i\} \tag{5.3.3}$$

Since

$$B = \bigcup_i B_i \quad i \neq j \implies B_i \cap B_j = \emptyset \tag{5.3.4}$$

then the equation

$$\int_B d\mu(x_1, x_2) f(x_1, x_2) = \sum_n \mu(B_n) y_n \tag{5.3.5}$$

follows from the definition 4.1.2. The equation

$$\mu(B_n) = \int_{X_1} d\mu_1(x_1)\mu_2(B_{n \cdot 2x_1}) \tag{5.3.6}$$

where

$$B_{n \cdot 2x_1} = \{x_2 \in X_2 : (x_1, x_2) \in B_n\} \tag{5.3.7}$$

follows from the theorem 5.2.3. The equation

$$\int_B d\mu(x_1, x_2) f(x_1, x_2) = \sum_n \left( \int_{X_1} d\mu_1(x_1)\mu_2(B_{n \cdot 2x_1}) \right) y_n \tag{5.3.8}$$

---

[5.6]See also the theorem [1]-4 on page 359.

follows from equations (5.3.5), (5.3.6). According to the definition 4.1.2, the equation

$$\int_{B_{2x_1}} d\mu_2(x_2) f(x_1, x_2) = \sum_n \mu_2(B_{n \cdot 2x_1}) y_n \tag{5.3.9}$$

follows from the equation (5.3.3), (5.3.7). According to theorems 4.2.3, 4.2.8, the equation

$$\begin{aligned}
\int_{X_1} d\mu_1(x_1) \int_{B_{2x_1}} d\mu(x_2) f(x_1, x_2) &= \int_{X_1} d\mu_1(x_1) \left( \sum_n \mu_2(B_{n \cdot 2x_1}) y_n \right) \\
&= \sum_n \int_{X_1} d\mu_1(x_1)(\mu_2(B_{n \cdot 2x_1}) y_n) \\
&= \sum_n \left( \int_{X_1} d\mu_1(x_1) \mu_2(B_{n \cdot 2x_1}) \right) y_n
\end{aligned} \tag{5.3.10}$$

follows from the equation (5.3.9). The equation (5.3.1) follows from equations (5.3.8), (5.3.10). $\odot$

According to the theorem 3.3.2, there exists a sequence of simple integrable over the set $X$ maps

$$f_n : X_1 \times X_2 \to A$$

converging uniformly to $f$

$$f(x) = \lim_{n \to \infty} f_n(x) \tag{5.3.11}$$

According to the lemma 5.3.2, the equation

$$\int_B d\mu(x_1, x_2) f_n(x_1, x_2) = \int_{X_1} d\mu(x_1) \int_{B_{2x_1}} d\mu(x_2) f_n(x_1, x_2) \tag{5.3.12}$$

holds for every map $f_n$. Equations

$$\int_B d\mu(x_1, x_2) f(x_1, x_2) = \lim_{n \to \infty} \int_B d\mu(x_1, x_2) f_n(x_1, x_2) \tag{5.3.13}$$

$$\int_{B_{2x_1}} d\mu_2(x_2) f(x_1, x_2) = \lim_{n \to \infty} \int_{B_{2x_1}} d\mu_2(x_2) f_n(x_1, x_2) \tag{5.3.14}$$

follow from definition 4.2.1. Since the map $f$ has compact range, then there exists $M > 0$ such that

$$\|f(x_1, x_2)\| \le M \tag{5.3.15}$$

According to the construction in the proof of the theorem 3.3.2, the statement

$$\|f_n(x_1, x_2)\| < M + \frac{1}{n} \le M + 1 \tag{5.3.16}$$

follows from the statement (5.3.15). According to the theorem 4.2.10, the statement

$$\left| \int_{B_{2x_1}} d\mu_2(x_2) f_n(x_1, x_2) \right| \le \mu_2(B_{2x_1})(M + 1) \tag{5.3.17}$$

follows from the statement (5.3.16). According to theorems 4.2.8, 5.2.3, the map

$$x_1 \to \mu_2(B_{2x_1})(M + 1)$$

is $\mu_1$-integralable. According to the theorem 4.4.2, from (5.3.11), (5.3.14), it follows that the map
$$x_1 \to \int_{B_{2x_1}} d\mu_2(x_2) f(x_1, x_2)$$
is $\mu_1$-integralable and

(5.3.18)
$$\begin{aligned}
&\int_{X_1} d\mu_1(x_1) \int_{B_{2x_1}} d\mu_2(x_2) f(x_1, x_2) \\
&= \int_{X_1} d\mu_1(x_1) \lim_{n\to\infty} \int_{B_{2x_1}} d\mu_2(x_2) f_n(x_1, x_2) \\
&= \lim_{n\to\infty} \int_{X_1} d\mu_1(x_1) \int_{B_{2x_1}} d\mu_2(x_2) f_n(x_1, x_2)
\end{aligned}$$

follows from (5.3.14). The equation (5.3.1). follows from equations (5.3.13), (5.3.18).
$\square$

CHAPTER 6

# CHAPTER 7

# Index

# CHAPTER 8

# Special Symbols and Notations

# Интеграл отображения в абелеву Ω-группу

Александр Клейн

Aleks_Kleyn@MailAPS.org
http://AleksKleyn.dyndns-home.com:4080/
http://sites.google.com/site/AleksKleyn/
http://arxiv.org/a/kleyn_a_1
http://AleksKleyn.blogspot.com/

Аннотация. Кольцо, модуль и алгебра имеют то общее, что они являются абелевыми группами относительно сложения. Этого свойства достаточно для изучения операции интегрирования. Рассмотрен интеграл измеримого отображения в нормированную абелеву $\Omega$-группу. Теория интегрирования отображений в $\Omega$-группу имеет много общего с теорией интегрирования функций действительного переменного. Однако многие утверждения необходимо изменить, так как они неявно предполагают компактность области значений либо отношение полного порядка в $\Omega$-группе.

# Оглавление

Глава 1

# Предисловие

### 1.1. Предисловие к изданию 1

Теория измеримых функций и интеграла, построенная в середине XX века, решает немало математических проблем. Однако эта теория рассматривает функции действительного переменного.

Когда я начал исследовать алгебры с непрерывным базисом, я понял, что мне необходима аналогичная теория интегрирования в нормированном векторном пространстве. Я рассматривал интеграл отображений в кольца, модули и алгебры. Эти алгебраические структуры имеют то общее, что они являются абелевыми группами относительно сложения. Этого свойства достаточно для изучения операции интегрирования.

Естественно было рассмотреть проблему в общем, и я вспомнил, что однажды такая задача передо мной стояла. Так возникло решение изучить нормированные $\Omega$-группы и теорию интегрирования отображений в $\Omega$-группу.

Из многочисленных ссылок на аналогичные определения и теоремы, читатель увидит, что теория интегрирования отображений в $\Omega$-группу имеет много общего с теорией интегрирования функций действительного переменного. Однако многие утверждения необходимо изменить, так как они неявно предполагают компактность области значений либо отношение полного порядка в $\Omega$-группе.

Я уделил большое внимание топологии нормированной $\Omega$-группы. Для определения интеграла измеримого отображения $f$, мне надо рассмотреть последовательность простых отображений, равномерно сходящихся к отображению $f$.

Октябрь, 2013

### 1.2. Предисловие к изданию 2

Когда я писал текст издания 1, я был вынужден остановиться на теореме Фубини. Доказательство теоремы существенно зависит от свойств действительных чисел. Это неприемлемо в случае произвольной $\Omega$-группы. Однако без теоремы Фубини невозможно выделить матрицу линейного преобразования алгебры с непрерывным базисом. Поэтому я вернулся к статье с целью доказать теорему Фубини.

Сейчас для меня очень важен раздел 4.3. Из теорем раздела 4.3 следует, что свойства интеграла похожи на свойства меры. Заманчиво рассмотреть более общее определение меры. Текст статьи готов к обобщению. Однако должна быть решена следующая задача. Пусть $\mu$ - мера со значением в $\Omega$-группе. Пусть $\mu(A) \neq 0$, $\mu(B) \neq 0$. Мы требуем, чтобы $\mu(A \cup B) \neq \emptyset$. Это требование может быть частью определения меры.

Так как я рассматриваю интеграл отображения в $\Omega$-группу $A$, я должен был определить взаимодействие действительного числа и $A$-числа. Я рассматриваю представление поля действительных чисел в $\Omega$-группу $A$. Эта модель удобна при изучении алгебры с непрерывным базисом.

<div align="right">Март, 2014</div>

### 1.3. Соглашения

Соглашение 1.3.1. *Элемент $\Omega$-группы $A$ называется $A$-**числом**. Например, комплексное число также называется $C$-число, а кватернион называется $H$-число.* □

Соглашение 1.3.2. *Пусть $A$ - $\Omega_1$-алгебра. Пусть $B$ - $\Omega_2$-алгебра. Запись*
$$A \xrightarrow{\ \ *\ \ } B$$
*означает, что определено представление $\Omega_1$-алгебры $A$ в $\Omega_2$-алгебре $B$.* □

Без сомнения, у читателя могут быть вопросы, замечания, возражения. Я буду признателен любому отзыву.

Глава 2

# Мера

## 2.1. Алгебра множеств

ОПРЕДЕЛЕНИЕ 2.1.1. Непустая система множеств $\mathcal{S}$ называется **полукольцом множеств**,[2.1] если

2.1.1.1: $\emptyset \in \mathcal{S}$
2.1.1.2: Если $A, B \in \mathcal{S}$, то $A \cap B \in \mathcal{S}$
2.1.1.3: Если $A, A_1 \in \mathcal{S}, A_1 \subset A$, то множество $A$ может быть представлено в виде

$$(2.1.1) \qquad A = \bigcup_{i=1}^{n} A_i \quad A_i \in \mathcal{S}$$

где $i \neq j => A_i \cap A_j = \emptyset$

Представление (2.1.1) множества $A$ называется **конечным разложением множества** $A$. □

ОПРЕДЕЛЕНИЕ 2.1.2. Непустая система множеств $\mathcal{R}$ называется **кольцом множеств**,[2.2] если условие $A, B \in \mathcal{R}$ влечёт $A \triangle B, A \cap B \in \mathcal{R}$. Множество $E \in \mathcal{R}$ называется **единицей кольца множеств**, если

$$A \cap E = A$$

Кольцо множеств с единицей называется **алгеброй множеств**. □

ЗАМЕЧАНИЕ 2.1.3. Для любых $A, B$

$$A \cup B = (A \triangle B) \triangle (A \cap B)$$
$$A \setminus B = A \triangle (A \cap B)$$

Следовательно, если $A, B \in \mathcal{R}$, то $A \cup B \in \mathcal{R}, A \setminus B \in \mathcal{R}$. □

ТЕОРЕМА 2.1.4. *Кольцо множеств $\mathcal{R}$ является полукольцом.*

ДОКАЗАТЕЛЬСТВО. Пусть $A, A_1 \in \mathcal{R}, A_1 \subset A$. Тогда

$$A = A_1 \cup A_2$$

где

$$A_2 = A \setminus A_1 \in \mathcal{R}$$

□

ТЕОРЕМА 2.1.5. *Пересечение $\mathcal{R} = \bigcap \mathcal{R}_i$ любого множества колец является кольцом.*[2.3]

---

[2.1]Смотри так же определение [1]-2, страница 43.
[2.2]Смотри так же определение [1]-1, страница 41.
[2.3]Смотри также теорему [1]-1 на странице 42.

Доказательство. Пусть $A, B \in \mathcal{R}$. Тогда для любого $i$, $A, B \in \mathcal{R}_i$. Согласно определению 2.1.2, для любого $i$, $A \triangle B, A \cap B \in \mathcal{R}_i$. Следовательно, $A \triangle B, A \cap B \in \mathcal{R}$. Согласно определению 2.1.2, множество $\mathcal{R}$ является кольцом множеств. □

Теорема 2.1.6. *Для любой непустой системы множеств $\mathcal{C}$, существует одно и только одно кольцо множеств[2.4] $\mathcal{R}(\mathcal{C})$, содержащее $\mathcal{C}$ и содержащееся в любом кольце множеств $\mathcal{R}$ таком, что $\mathcal{C} \subset \mathcal{R}$.*

Доказательство. Пусть кольца множеств $\mathcal{R}_1, \mathcal{R}_2, \mathcal{R}_1 \neq \mathcal{R}_2$, удовлетворяют условию теоремы. Согласно теореме 2.1.5, множество $\mathcal{R}_1 \cap \mathcal{R}_2$, является кольцом множеств, удовлетворяющим условию теоремы. Следовательно, если кольцо множеств $\mathcal{R}(\mathcal{C})$ существует, то оно единственно.

Пусть
$$R = \bigcup_{X \in \mathcal{C}} X$$
Множество $\mathcal{B}(R)$ всех подмножеств множества $R$ является кольцом множеств и $\mathcal{C} \subseteq \mathcal{B}(R)$. Пусть $\Sigma$ - совокупность колец множеств $\mathcal{R}$ таких, что $\mathcal{C} \subseteq \mathcal{R} \subseteq \mathcal{B}(R)$. Тогда, согласно теореме 2.1.5, множество $\mathcal{P} = \bigcap_{\mathcal{R} \in \Sigma} \mathcal{R}$ является кольцом множеств, которое удовлетворяет теореме. □

Теорема 2.1.7. *Пусть $\mathcal{C}$ - непустая система множеств. Пусть*
$$R = \bigcup_{X \in \mathcal{C}} X \tag{2.1.2}$$

*Если*
$$R \in \mathcal{C} \tag{2.1.3}$$

*то кольцо множеств $\mathcal{R}(\mathcal{C})$ является алгеброй множеств.*

Доказательство. Множество $\mathcal{B}(R)$ всех подмножеств множества $R$ является кольцом множеств и $\mathcal{C} \subseteq \mathcal{B}(R)$. Пусть $\Sigma$ - совокупность колец множеств $\mathcal{R}$ таких, что
$$\mathcal{C} \subseteq \mathcal{R} \subseteq \mathcal{B}(R) \tag{2.1.4}$$
Из (2.1.3), (2.1.4) следует, что $R \in \mathcal{R}$ для любого $\mathcal{R} \in \Sigma$. Тогда, согласно теореме 2.1.5, множество $\mathcal{P} = \bigcap_{\mathcal{R} \in \Sigma} \mathcal{R}$ является наименьшим кольцом множеств таким, что $R \in \mathcal{P}, \mathcal{C} \subseteq \mathcal{P}$. Из (2.1.2) следует, что $R \cap A = A$ для любого множества $A \in \mathcal{P}$. Согласно определению 2.1.2, кольцо множеств $\mathcal{P}$ является алгеброй множеств. □

Теорема 2.1.8. *Для любой непустой системы множеств $\mathcal{C}$, существует одна и только одна алгебра множеств $\mathcal{A}(\mathcal{C})$, содержащая $\mathcal{C}$ и содержащаяся в любой алгебре множеств $\mathcal{R}$ такой, что $\mathcal{C} \subset \mathcal{R}$.*

Доказательство. Пусть
$$R = \bigcup_{X \in \mathcal{C}} X$$

---
[2.4] Смотри также теорему [1]-2 на странице 42.

Теорема следует из теорем 2.1.6, 2.1.7, если мы положим
$$\mathcal{A}(\mathcal{C}) = \mathcal{R}(\{R\} \cup \mathcal{C})$$

$\square$

Лемма 2.1.9. Пусть[2.5] $\mathcal{S}$ - полукольцо. Пусть $A$, $A_1, ..., A_n \in \mathcal{S}$, $i \neq j \Rightarrow A_i \cap A_j = \emptyset$. Тогда существует конечное разложение множества $A$

$$A = \bigcup_{i=1}^{s} A_i \quad s \geq n$$

Доказательство. Для $n = 1$ лемма следует из утверждения 2.1.1.3.

Пусть лемма верна для $n = m$. Пусть множества $A_1, ..., A_{m+1}$ удовлетворяют условию леммы. Согласно предположению

(2.1.5) $$A = A_1 \cup ... \cup A_m \cup B_1 \cup ... \cup B_p$$

где $B_i \in \mathcal{S}, i = 1, ..., p,\ i \neq j \Rightarrow A_i \cap A_j = \emptyset,\ A_i \cap B_j = \emptyset,\ i \neq j \Rightarrow B_i \cap B_j = \emptyset$. Согласно утверждению 2.1.1.2

(2.1.6) $$B_{i \cdot 1} = A_{m+1} \cap B_i \in \mathcal{S}$$

Согласно утверждению 2.1.1.3

(2.1.7) $$B_i = B_{i \cdot 1} \cup ... \cup B_{i \cdot r_s} \quad B_{i \cdot j} \in \mathcal{S}$$

Так как $A_{m+1} \subset B_1 \cup ... \cup B_p$, то

$$A = A_1 \cup ... \cup A_m \cup A_{m+1} \cup \bigcup_{i=1}^{p} \bigcup_{j=2}^{r_i} B_{i \cdot j}$$

является следствием (2.1.5), (2.1.6), (2.1.7). Следовательно, лемма верна для $n = m + 1$.

Согласно математической индукции, лемма верна для любого $n$. $\square$

Лемма 2.1.10. Для любой конечной системы множеств $A_1, ..., A_n \in \mathcal{S}$ существует конечная система множеств $B_1, ..., B_t \in \mathcal{S}, i \neq j \Rightarrow B_i \cap B_j = \emptyset$, такая, что

(2.1.8) $$A_i = \bigcup_{j \in M_i} B_j$$

где $M_i \subset \{1, ..., t\}$.

Доказательство. Для $n = 1$ лемма очевидна, так как достаточно положить $t = 1$, $B_1 = A_1$.

Пусть лемма верна для $n = m$. Рассмотрим систему множеств $A_1, ..., A_{m+1}$. Пусть $B_i \in \mathcal{S}, i = 1, ..., p$, - множества, удовлетворяющие условию леммы по отношению к $A_1, ..., A_m$. Согласно утверждению 2.1.1.2

(2.1.9) $$B_{i \cdot 1} = A_{m+1} \cap B_i \in \mathcal{S}$$

Согласно утверждению 2.1.1.3

(2.1.10) $$B_i = B_{i \cdot 1} \cup ... \cup B_{i \cdot r_s} \quad B_{i \cdot j} \in \mathcal{S}$$

---

[2.5]Леммы 2.1.9, 2.1.10 и теорема 2.1.11 аналогичны леммам 1, 2 и теореме 3, [1], страницы 43, 45.

Разложение

(2.1.11) $$A_{m+1} = \bigcup_{s=1}^{t} B_{s \cdot 1} \cup \bigcup_{p=1}^{q} B'_p \quad B'_p \in \mathcal{S}$$

следует из леммы 2.1.9. Разложение

$$A_i = \bigcup_{j \in M_i} \bigcup_{p=1}^{r_j} B_{j \cdot p}$$

следует из равенств (2.1.8), (2.1.10). Из равенств (2.1.9), (2.1.11) следует, что $B_i \cap B'_p = \emptyset$. Следовательно, из равенства (2.1.10) следует, что $B_{i \cdot j} \cap B'_p = \emptyset$. Следовательно, множества $B_{i \cdot j}$, $B'_p$ удовлетворяют условиям леммы по отношению к $A_1, ..., A_{m+1}$. Следовательно, лемма верна для $n = m+1$.

Согласно математической индукции, лемма верна для любого $n$. $\square$

Теорема 2.1.11. *Пусть $\mathcal{S}$ - полукольцо множеств. Система $\mathcal{R}$ множеств $A$, которые имеют конечное разложение*

$$A = \bigcup_{i=1}^{n} A_i \quad A_i \in \mathcal{S}$$

*является* **кольцом множеств, порождённым полукольцом множеств $\mathcal{S}$**.

Доказательство. Пусть $A, B \in \mathcal{R}$ Тогда

(2.1.12) $$A = \bigcup_{i=1}^{n} A_i \quad A_i \in \mathcal{S}$$

(2.1.13) $$B = \bigcup_{i=1}^{m} B_i \quad B_i \in \mathcal{S}$$

Из (2.1.12), (2.1.13) и утверждения 2.1.1.2 следует, что

$$C_{ij} = A_i \cap B_j \in \mathcal{S}$$

Согласно лемме 2.1.9

(2.1.14) $$A_i = \bigcup_{j=1}^{m} C_{ij} \cup \bigcup_{k=1}^{r_i} D_{ik} \quad D_{ik} \in \mathcal{S}$$
$$B_j = \bigcup_{i=1}^{n} C_{ij} \cup \bigcup_{k=1}^{s_j} E_{kj} \quad E_{kj} \in \mathcal{S}$$

Из (2.1.14) следует, что

$$A \cup B = \bigcup_{ij} C_{ij} \in \mathcal{R}$$

$$A \Delta B = \bigcup_{ik} D_{ik} \cup \bigcup_{jl} E_{jl} \in \mathcal{R}$$

Следовательно, $\mathcal{R}$ является кольцом множеств. $\square$

Теорема 2.1.12. *Пусть полукольцо множеств $\mathcal{C}$ содержит единицу. Тогда кольцо множеств $\mathcal{R}(\mathcal{C})$ является алгеброй множеств.*

Доказательство. Теорема является следствием теорем 2.1.8, 2.1.11. $\square$

Определение 2.1.13. Кольцо множеств $\mathcal{R}$ называется $\sigma$-**кольцом множеств**,[2.6] если условие $A_i \in \mathcal{R}$, $i = 1, ..., n, ...,$ влечёт
$$\bigcup_n A_n \in \mathcal{R}$$

$\sigma$-Кольцо множеств с единицей называется $\sigma$-**алгеброй множеств**. □

Теорема 2.1.14. *Пересечение* $\mathcal{R} = \bigcap \mathcal{R}_i$ *любого множества $\sigma$-колец является $\sigma$-кольцом.*[2.7]

Доказательство. Если $A_i \in \mathcal{R}$, $i = 1, ..., n, ...,$ . то для любых $i, j$, $A_i \in \mathcal{R}_j$. Следовательно, для любого $j$,
$$\bigcup_i A_i \in \mathcal{R}_j$$
Следовательно,
$$\bigcup_i A_i \in \mathcal{R}$$
Согласно определению 2.1.13, множество $\mathcal{R}$ является $\sigma$-кольцом множеств. □

Теорема 2.1.15. *Для любой непустой системы множеств $\mathcal{C}$, существует одно и только одно $\sigma$-кольцо множеств*[2.8] $\mathcal{R}_\sigma(\mathcal{C})$, *содержащее $\mathcal{C}$ и содержащееся в любом $\sigma$-кольце множеств $\mathcal{R}$ таком, что $\mathcal{C} \subset \mathcal{R}$.*

Доказательство. Пусть $\sigma$-кольца множеств $\mathcal{R}_1$, $\mathcal{R}_2$, $\mathcal{R}_1 \neq \mathcal{R}_2$, удовлетворяют условию теоремы. Согласно теореме 2.1.14, множество $\mathcal{R}_1 \cap \mathcal{R}_2$, является $\sigma$-кольцом множеств, удовлетворяющим условию теоремы. Следовательно, если $\sigma$-кольцо множеств $\mathcal{R}_\sigma(\mathcal{C})$ существует, то оно единственно.

Пусть
$$R = \bigcup_{X \in \mathcal{C}} X$$
Множество $\mathcal{B}(R)$ всех подмножеств множества $R$ является $\sigma$-кольцом множеств и $\mathcal{C} \subseteq \mathcal{B}(R)$. Пусть $\Sigma$ - совокупность $\sigma$-колец множеств $\mathcal{R}$ таких, что $\mathcal{C} \subseteq \mathcal{R} \subseteq \mathcal{B}(R)$. Тогда, согласно теореме 2.1.5, множество $\mathcal{P} = \bigcap_{\mathcal{R} \in \Sigma} \mathcal{R}$ является $\sigma$-кольцом множеств, которое удовлетворяет теореме. □

Теорема 2.1.16. *Для любой непустой системы множеств $\mathcal{C}$, существует одна и только одна $\sigma$-алгебра множеств $\mathcal{A}_\sigma(\mathcal{C})$, содержащая $\mathcal{C}$ и содержащаяся в любой $\sigma$-алгебре множеств $\mathcal{R}$ такой, что $\mathcal{C} \subset \mathcal{R}$.*

Доказательство. Пусть
$$R = \bigcup_{X \in \mathcal{C}} X$$
Теорема следует из теорем 2.1.15, 2.1.7, если мы положим
$$\mathcal{A}_\sigma(\mathcal{C}) = \mathcal{R}_\sigma(\{R\} \cup \mathcal{C})$$

□

---

[2.6]Смотри аналогичное определение в [1], с. 45, определение 3.
[2.7]Смотри также теорему [1]-1 на странице 42.
[2.8]Смотри также теорему [1]-2 на странице 42.

## 2.2. Мера

Определение 2.2.1. Пусть $\mathcal{C}_m$ - полукольцо множеств.[2.9] Отображение
$$m : \mathcal{C}_m \to R$$
называется **мерой**, если

2.2.1.1: $m(A) \geq 0$

2.2.1.2: Отображение $m$ аддитивно. Если множество $A \in \mathcal{C}_m$ имеет конечное разлоение
$$A = \bigcup_{i=1}^n A_i \quad A_i \in \mathcal{C}_m$$
где $i \neq j => A_i \cap A_j = \emptyset$, то
$$m(A) = \sum_{i=1}^n m(A_i)$$

□

Теорема 2.2.2. $m(\emptyset) = 0$.

Доказательство. Так как $\emptyset = \emptyset \cup \emptyset$, $\emptyset \cap \emptyset = \emptyset$, то теорема является следствием утверждений 2.1.1.1, 2.2.1.2. □

Определение 2.2.3. Мера $\mu$ называется[2.10] **продолжением меры** $m$, если $\mathcal{C}_m \subseteq \mathcal{C}_\mu$ и $\mu(A) = m(A)$, $A \in \mathcal{C}_m$. □

Теорема 2.2.4. *Пусть*[2.11] $\mathcal{R}(\mathcal{C}_m)$ - *кольцо множеств, порождённое полукольцом множеств* $\mathcal{C}_m$. *Для меры* $m$, *заданной на полукольце множеств* $\mathcal{C}_m$, *существует одно и только одно продолжение* $\mu$, *заданное на кольце множеств* $\mathcal{R}(\mathcal{C}_m)$.

Доказательство. Согласно теореме 2.1.11, для каждого множества $A \in \mathcal{R}(\mathcal{C}_m)$ существует конечное разложение

$$A = \bigcup_{i=1}^n A_i \quad A_i \in \mathcal{C}_m \tag{2.2.1}$$

где $i \neq j => A_i \cap A_j = \emptyset$. Положим

$$\mu(A) = \sum_{i=1}^n m(A_i) \tag{2.2.2}$$

Если существует два конечных разложения
$$A = \bigcup_{i=1}^n A_i = \bigcup_{j=1}^m B_j \quad A_i, B_j \in \mathcal{C}_m$$
Так как $A_i \cap B_j \in \mathcal{C}_m$ согласно утверждению 2.1.1.2, то
$$\sum_{i=1}^n m(A_i) = \sum_{i=1}^n \sum_{j=1}^m m(A_i \cap B_j) = \sum_{j=1}^m m(B_j)$$

---

[2.9]Смотри также определение [1]-1 на странице 265.

[2.10]Смотри также определение [1]-2 на странице 266.

[2.11]Смотри также теорему [1]-1 на странице 266.

следует из утверждения 2.2.1.2. Следовательно, величина $\mu(A)$, определённая равенством (2.2.2), не зависит от выбора конечного разложения (2.2.1).

Следовательно, мы построили отображение

$$\mu : \mathcal{R}(\mathcal{C}_m) \to R$$

которое удовлетворяет определению 2.2.1.

Для любого продолжения $\mu'$ меры $m$ и для конечного разложения (2.2.1), из утверждения 2.2.1.2 и определения 2.2.3 следует, что

$$\mu'(A) = \sum_{i=1}^{n} \mu'(A_i) = \sum_{i=1}^{n} \mu(A_i) = \mu(A)$$

Следовательно, мера $\mu'$ совпадает с мерой $\mu$, определённой равенством (2.2.2). □

ОПРЕДЕЛЕНИЕ 2.2.5. Мера $\mu$ называется **полной**, если из условий $B \subset A$, $\mu(A) = 0$ следует, что множество $B$ измеримо. □

ТЕОРЕМА 2.2.6. *Пусть на множестве $X$ определена мера $\mu$. Пусть*[2.12]

(2.2.3) $$A \subset B \quad A, B \in \mathcal{C}_\mu$$

*Тогда*

(2.2.4) $$\mu(A) \leq \mu(B)$$

ДОКАЗАТЕЛЬСТВО. Равенство

(2.2.5) $$B = A \cup (B \setminus A)$$

следует из утверждения (2.2.3). Согласно замечанию 2.1.3, из утверждения (2.2.3) следует, что

$$B \setminus A \in \mathcal{C}_\mu$$

Равенство

(2.2.6) $$\mu(B) = \mu(A) + \mu(B \setminus A)$$

следует из равенства (2.2.5). Согласно утверждению 2.2.10.1,

(2.2.7) $$\mu(B \setminus A) \geq 0$$

Утверждение (2.2.4) является следствием утверждений 2.2.10.2, (2.2.6). □

ТЕОРЕМА 2.2.7. *Пусть на множестве $X$ определена мера $\mu$. Пусть*

$$A \subset \bigcup_i A_i \quad A, A_i \in \mathcal{C}_\mu$$

*где $\{A_i\}$ конечная или счётная система множеств. Тогда*

$$\mu(A) \leq \mu(B)$$

ДОКАЗАТЕЛЬСТВО. Мы докажем теорему для случая 2 множеств. Общий случай легко доказать методом математической индукции. Если

$$A \subseteq A_1 \cup A_2$$

то теорема следует из теоремы 2.2.6, утверждения

$$\mu(A_1 \cup A_2) = \mu(A_1) + \mu(A_2) - \mu(A_1 \cap A_2) \leq \mu(A_1) + \mu(A_2)$$

---

[2.12] Смотри также теорему [1]-2, с. 267.

и утверждения 2.2.10.1. □

Теорема 2.2.8. *Пусть на множестве $X$ определена полная мера $\mu$. Пусть*
$$A \subset B \quad A, B \in \mathcal{C}_\mu$$
*Если $\mu(B) = 0$, то $\mu(A) = 0$.*

Доказательство. Теорема является следствием определения 2.2.5 и теоремы 2.2.6. □

Замечание 2.2.9. Мы будем полагать, что рассматриваемая мера является полной мерой. □

Определение 2.2.10. Пусть $\mathcal{C}_\mu$ - $\sigma$-алгебра множеств множества $F$.[2.13] Отображение
$$\mu : \mathcal{C}_\mu \to R$$
в поле действительных чисел $R$ называется **$\sigma$-аддитивной мерой**, если для любого множества $X \in \mathcal{C}_\mu$ выполнены следующие условия.

2.2.10.1: $\mu(X) \geq 0$

2.2.10.2: Пусть
$$X = \bigcup_i X_i \quad i \neq j => X_i \cap X_j = \emptyset$$
конечное или счётное объединение множеств $X_n \in \mathcal{C}_\mu$. Тогда
$$\mu(X) = \sum_i \mu(X_i)$$
где ряд в правой части сходится абсолютно.

□

Теорема 2.2.11. *Пусть $m$ - $\sigma$-аддитивная мера, определённая на полукольце множеств $\mathcal{C}_m$. Тогда продолжение $\mu$ меры $m$, заданное на кольце множеств $\mathcal{R}(\mathcal{C}_m)$, является $\sigma$-аддитивной мерой.*

Доказательство. Пусть $A, B_n \in \mathcal{R}(\mathcal{C}_m)$, $n = 1, 2, ...$, $i \neq j => B_i \cap B_j = \emptyset$. Пусть
$$A = \bigcup_{n=1}^\infty B_n$$
Согласно теореме 2.1.11, существуют конечные разложения
$$(2.2.8) \quad A = \bigcup_j A_j \quad B_n = \bigcup_j B_{nj}$$
где
$$A_k \cap A_l = \emptyset \quad B_{nk} \cap B_{nl} = \emptyset \quad k \neq l$$
Пусть $C_{nil} = B_{ni} \cap A_l$. Очевидно, что множества $C_{nil}$ попарно не пересекаются и
$$(2.2.9) \quad A_j = \bigcup_{n=1}^\infty \bigcup_i C_{nij} \quad B_{ni} = \bigcup_j C_{nij}$$

---

[2.13]Смотри аналогичные определения в [1], определение 1 на странице 265 и определение 3 на странице 268.

Из равенства (2.2.9) и утверждения 2.2.10.2 следует, что

$$m(A_j) = \sum_{n=1}^{\infty} \sum_{i} m(C_{nij})$$
(2.2.10)
$$m(B_{ni}) = \sum_{j} m(C_{nij})$$

Из равенства (2.2.8) и теоремы 2.2.4 следует, что

$$\mu(A) = \sum_{j} m(A_j)$$
(2.2.11)
$$\mu(B_n) = \sum_{i} m(B_{ni})$$

Равенство

$$\mu(A) = \sum_{n=1}^{\infty} \mu(B_n)$$

является следствием (2.2.10), (2.2.11). □

Теорема 2.2.12 (Непрерывность $\sigma$-аддитивной меры). *Пусть*[2.14]

(2.2.12) $$A_1 \supset A_2 \supset ...$$

*последовательность $\mu$-измеримых множеств. Тогда*

$$\mu(A) = \lim_{n \to \infty} \mu(A_n)$$

*где* $A = \bigcap_{n} A_n$ .

Доказательство. Мы рассмотрим случай $A = \emptyset$. Общий случай сводится к этому заменой $A_n$ на $A_n \setminus A$. Согласно утверждению (2.2.12)

(2.2.13) $$A_n = (A_n \setminus A_{n+1}) \cup (A_{n+1} \setminus A_{n+2}) \cup ... \quad n = 1, ...$$

Согласно утверждению 2.2.10.2, равенство

(2.2.14) $$\mu(A_n) = \sum_{k=n}^{\infty} \mu(A_k \setminus A_{k+1}) \quad n = 1, ...$$

следует из (2.2.13). Поскольку ряд (2.2.14) для $n = 1$ сходится, его остаток (2.2.14) стремится к 0, когда $n \to \infty$. Следовательно,

$$\lim_{n \to \infty} \mu(A_n) = 0$$

□

Теорема 2.2.13 (Непрерывность $\sigma$-аддитивной меры). *Пусть*[2.15]

$$A_1 \subset A_2 \subset ...$$

*последовательность $\mu$-измеримых множеств. Тогда*

$$\mu(A) = \lim_{n \to \infty} \mu(A_n)$$

*где* $A = \bigcup_{n} A_n$ .

---

[2.14]Смотри также теорему [1]-9 на странице 261.

[2.15]Смотри также следствие теоремы [1]-9 на страницах 261, 262.

Доказательство. Теорема следует из теоремы 2.2.12, если мы рассмотрим множества $X \setminus A_n$. □

## 2.3. Лебегово продолжение меры

Определение 2.3.1. Пусть $m$ - $\sigma$-аддитивная мера на полукольце $\mathcal{S}$ с единицей $E$.[2.16] **Внешняя мера** множества $A \subseteq E$ определена равенством

$$\mu^*(A) = \inf_{A \subseteq \bigcup_k B_k} \sum_k m(B_k)$$

где точная верхняя грань берётся по всем покрытиям множества $A$ конечными или счётными системами множеств $B_n \in \mathcal{S}$. □

Теорема 2.3.2 (счётная полуаддитивность). *Если*[2.17]

$$(2.3.1) \qquad A \subseteq \bigcup_n A_n$$

*где $\{A_n\}$ - конечная или счётная система множеств, то*

$$(2.3.2) \qquad \mu^*(A) \leq \sum_n \mu^*(A_n)$$

Доказательство. Согласно определению 2.3.1, для каждого $A_n$ существует конечная или счётная система множеств $\{P_{nk}\}$ такая, что

$$(2.3.3) \qquad A_n \subseteq \bigcup_k P_{nk}$$

$$(2.3.4) \qquad \sum_k m(P_{nk}) \leq \mu^*(A_n) + \frac{\epsilon}{2^n}$$

Из (2.3.1), (2.3.3) следует, что

$$(2.3.5) \qquad A \subseteq \bigcup_n \bigcup_k P_{nk}$$

Согласно утверждению (2.3.5) и определению 2.3.1, неравенство

$$(2.3.6) \qquad \mu^*(A) \leq \sum_n \sum_k m(P_{nk}) \leq \sum_n \mu^*(A_n) + \epsilon$$

следует из неравенства (2.3.4). (2.3.2) является следствием (2.3.6), так как $\epsilon$ произвольно. □

Теорема 2.3.3. *Для любых множеств $A$, $B$,*

$$(2.3.7) \qquad |\mu^*(A) - \mu^*(B)| \leq \mu^*(A \Delta B)$$

Доказательство. Если $\mu^*(A) \geq \mu^*(B)$, то утверждение

$$(2.3.8) \qquad \mu^*(A) \leq \mu^*(B) + \mu(A \Delta B)$$

является следствием утверждения

$$A \subseteq B \cup (A \Delta B)$$

---

[2.16]Смотри определение [1]-1, страница 272.

[2.17]Смотри также теоремы [1]-3 на страницах 256, 257, [1]-1 на странице 272.

и теоремы 2.3.2. Если $\mu^*(B) \geq \mu^*(A)$, то утверждение
$$(2.3.9) \qquad \mu^*(B) \leq \mu^*(A) + \mu(A\Delta B)$$
является следствием утверждения
$$B \subseteq A \cup (A\Delta B)$$
и теоремы 2.3.2. Утверждение (2.3.7) является следствием утверждений (2.3.8), (2.3.9). □

ОПРЕДЕЛЕНИЕ 2.3.4. Множество $A$ называется **измеримым по Лебегу**, для любого $\epsilon \in R, \epsilon > 0$, существует $B \in \mathcal{R}(\mathcal{S})$ такое, что[2.18]
$$\mu^*(A\Delta B) < \epsilon$$
□

Пусть $\mathcal{C}_\mu$ - система измеримых по Лебегу множеств.

ТЕОРЕМА 2.3.5. *Пусть $m$ - $\sigma$-аддитивная мера на полукольце $\mathcal{S}$ с единицей $E$.[2.19] Если множество $A$ измеримо по Лебегу, то множество $E \setminus A$ также измеримо по Лебегу.*

ДОКАЗАТЕЛЬСТВО. Пусть $A$ - множество, измеримое по Лебегу. Согласно определению 2.3.4, для любого $\epsilon \in R, \epsilon > 0$, существует $B \in \mathcal{R}(\mathcal{S})$ такое, что
$$(2.3.10) \qquad \mu^*(A\Delta B) < \epsilon$$
Согласно замечанию 2.1.3, $E \setminus B \in \mathcal{R}(\mathcal{S})$. Утверждение
$$\mu^*((E \setminus A)\Delta(E \setminus B)) < \epsilon$$
следует из утверждения (2.3.10) и равенства
$$A\Delta B = (E \setminus A)\Delta(E \setminus B)$$
□

ТЕОРЕМА 2.3.6. *Пусть $A_1, A_2 \in \mathcal{C}_\mu$. Тогда $A = A_1 \setminus A_2 \in \mathcal{C}_\mu$.*

ДОКАЗАТЕЛЬСТВО. Согласно определению 2.3.4, для любого $\epsilon \in R, \epsilon > 0$, существуют $B_1, B_2 \in \mathcal{R}(\mathcal{S})$ такое, что
$$(2.3.11) \qquad \begin{aligned} \mu^*(A_1\Delta B_1) &< \frac{\epsilon}{2} \\ \mu^*(A_2\Delta B_2) &< \frac{\epsilon}{2} \end{aligned}$$
Согласно замечанию 2.1.3, $B = B_1 \setminus B_2 \in \mathcal{R}(\mathcal{S})$. Утверждение
$$\mu^*(A\Delta B) < \epsilon$$
следует из утверждения (2.3.11), утверждения
$$(A_1 \setminus A_2)\Delta(B_1 \setminus B_2) \subseteq (A_1\Delta B_1) \cup (A_2\Delta B_2)$$
и теоремы 2.3.2. □

---

[2.18]Смотри определение [1]-2, страница 272.

[2.19]Смотри также замечание в [1] после определения 2 на странице 272.

Теорема 2.3.7. *Пусть $m$ - $\sigma$-аддитивная мера на полукольце $\mathcal{S}$ с единицей $E$.[2.20] Пусть $\mu$ - продолжение меры $m$ на кольцо множеств $\mathcal{R}(\mathcal{S})$. Любое множество $A \in \mathcal{R}(\mathcal{S})$ измеримо по Лебегу и*

$$\mu^*(A) = \mu(A) \tag{2.3.12}$$

Доказательство. Пусть $A \in \mathcal{R}(\mathcal{S})$. Согласно определению 2.1.1 и теореме 2.1.11, множество $A$ может быть представлено в виде

$$A = \bigcup_{i=1}^{n} A_i \quad A_i \in \mathcal{S}$$

где $i \neq j => A_i \cap A_j = \emptyset$. Согласно теореме 2.2.4,

$$\mu(A) = \sum_{i=1}^{n} m(A_i) \tag{2.3.13}$$

Согласно определению 2.3.1, так как множества $A_i$ покрывают множество $A$, то

$$\mu^*(A) \leq \sum_{i=1}^{n} m(A_i) = \mu(A) \tag{2.3.14}$$

следует из (2.3.13).

Пусть $\{Q_i\}, Q_i \in \mathcal{S}$, - конечная или счётная система множеств, покрывающих множество $A$. Согласно теоремам 2.2.4, 2.2.7,

$$\mu(A) \leq \sum_j m(Q_j) \tag{2.3.15}$$

Согласно определению 2.3.1,

$$\mu(A) \leq \mu^*(A) \tag{2.3.16}$$

следует из (2.3.15).

(2.3.12) следует из (2.3.14), (2.3.16). $\square$

Теорема 2.3.8. *Система $\mathcal{C}_\mu$ измеримых по Лебегу множеств является алгеброй множеств.*

Доказательство. Из теоремы 2.3.6, определения 2.1.2 и равенств

$$A_1 \cap A_2 = A_1 \setminus (A_1 \setminus A_2)$$
$$A_1 \cup A_2 = E \setminus ((E \setminus A_1) \cap (E \setminus A_2))$$
$$A_1 \Delta A_2 = (A_1 \cup A_2) \setminus (A_1 \cap A_2)$$

следует, что $\mathcal{C}_\mu$ является кольцом множеств. Кольцо множеств $\mathcal{C}_\mu$ является алгеброй множеств, так как $E \in \mathcal{C}_\mu$ является единицей кольца множеств $\mathcal{C}_\mu$. $\square$

Теорема 2.3.9. *Отображение $\mu^*(A)$ аддитивно на алгебре множеств $\mathcal{C}_\mu$.[2.21]*

---

[2.20]Смотри также замечания в [1] на страницах 256 и 272.

[2.21]Смотри также теоремы [1]-6 на странице 258 и [1]-3 на странице 273.

Доказательство. Для доказательства теоремы достаточно рассмотреть случай двух множеств. Пусть $A_1, A_2 \in \mathcal{C}_\mu$. Согласно определению 2.3.4, для любого $\epsilon \in R, \epsilon > 0,$ существуют $B_1, B_2 \in \mathcal{R}(\mathcal{S})$ такое, что

$$(2.3.17) \quad \begin{aligned} \mu^*(A_1 \Delta B_1) &< \frac{\epsilon}{2} \\ \mu^*(A_2 \Delta B_2) &< \frac{\epsilon}{2} \end{aligned}$$

Положим $A = A_1 \cup A_2$, $B = B_1 \cup B_2$. Согласно теореме 2.3.8, $A \in \mathcal{C}_\mu$. Так как $A_1 \cap A_2 = \emptyset$, то

$$(2.3.18) \quad B_1 \cap B_2 \subseteq (A_1 \Delta B_1) \cup (A_2 \Delta B_2)$$

Утверждение

$$(2.3.19) \quad \mu(B_1 \cap B_2) \leq \epsilon$$

следует из утверждения (2.3.18) и теорем 2.3.2, 2.3.7. Утверждение

$$(2.3.20) \quad \begin{aligned} |\mu(B_1) - \mu^*(A_1)| &\leq \frac{\epsilon}{2} \\ |\mu(B_2) - \mu^*(A_2)| &\leq \frac{\epsilon}{2} \end{aligned}$$

следует из утверждения (2.3.17) и теорем 2.3.3, 2.3.7. Так как мера аддитивна на алгебре множеств $\mathcal{R}(\mathcal{S})$, то утверждение

$$(2.3.21) \quad \mu(B) = \mu(B_1) + \mu(B_2) - \mu(B_1 \cap B_2) \geq \mu^*(A_1) + \mu^*(A_2) - 2\epsilon$$

следует из утверждений (2.3.19), (2.3.20). Утверждение

$$(2.3.22) \quad \mu^*(A) \geq \mu(B) - \mu^*(A \Delta B) \geq \mu(B) - 2\epsilon \geq \mu^*(A_1) + \mu^*(A_2) - 3\epsilon$$

следует из утверждения

$$A \Delta B \subseteq (A_1 \Delta B_1) \cup (A_2 \Delta B_2)$$

Утверждение

$$(2.3.23) \quad \mu^*(A_1) + \mu^*(A_2) \geq \mu^*(A) \geq \mu^*(A_1) + \mu^*(A_2) - 3\epsilon$$

следует из утверждения (2.3.22) и теоремы 2.3.2. Так как $\epsilon$ может быть выбрано произвольно малым, то

$$\mu^*(A) = \mu^*(A_1) + \mu^*(A_2)$$

следует из утверждения (2.3.23). $\square$

Определение 2.3.10. Если множество $A$ измеримо по Лебегу,[2.22] то величина $\mu(A) = \mu^*(A)$ называется **мерой Лебега**. Отображение $\mu$, определённое на алгебре множеств $\mathcal{C}_\mu$, называется **лебеговым продолжением меры** $m$. $\square$

Теорема 2.3.11. *Алгебра множеств $\mathcal{C}_\mu$ является $\sigma$-алгеброй[2.23] с единицей $E$.*

---

[2.22]Смотри определение [1]-2, страница 272.

[2.23]Смотри также теоремы [1]-7 на странице 259 и [1]-5 на страницах 273, 274.

Доказательство. Пусть $A_1, ...,$ - счётная система множеств, измеримых по Лебегу. Пусть

$$A = \bigcup_{n=1}^{\infty} A_n \tag{2.3.24}$$

Положим

$$A'_n = A_n \setminus \bigcup_{k=1}^{n-1} A_k \tag{2.3.25}$$

Из (2.3.24), (2.3.25) следует, что

$$A = \bigcup_{n=1}^{\infty} A'_n \tag{2.3.26}$$

где $i \neq j => A'_i \cap A'_j = \emptyset$. Согласно теореме 2.3.8, определению 2.1.2 и замечанию 2.1.3, все множества $A'_n$ измеримы по Лебегу. Согласно теореме 2.3.8 и определению 2.3.1, для любого $n$

$$\sum_{k=1}^{n} \mu(A'_k) = \mu\left(\bigcup_{k=1}^{n} A'_k\right) \leq \mu^*(A)$$

Следовательно, ряд $\sum_{n=1}^{\infty} \mu(A'_n)$ сходится. Следовательно, для любого $\epsilon \in R$, $\epsilon > 0$, существует $N$ такое, что

$$\sum_{n>N} \mu(A'_n) < \frac{\epsilon}{2} \tag{2.3.27}$$

Согласно теореме 2.3.8,

$$C = \sum_{n=1}^{N} \mu(A'_n) \in \mathcal{C}_\mu \tag{2.3.28}$$

Согласно определению 2.3.4, из утверждения (2.3.28) следует, что существует $B \in \mathcal{R}(\mathcal{S})$ такое, что

$$\mu^*(C \Delta B) < \frac{\epsilon}{2} \tag{2.3.29}$$

Так как

$$A \Delta B \subseteq (C \Delta B) \cup \left(\bigcup_{n>N} A'_n\right)$$

то из (2.3.27), (2.3.29) и теоремы 2.3.2 следует, что

$$\mu^*(A \Delta B) < \epsilon$$

Согласно определению 2.3.4, $A \in \mathcal{C}_\mu$.

Так как

$$\bigcap_n A_n = E \setminus \bigcup_n (E \setminus A_n)$$

то теорема является следствием теоремы 2.3.5. $\qquad\square$

Теорема 2.3.12. *Отображение $\mu(A)$ $\sigma$-аддитивно[2.24] на алгебре множеств $\mathcal{C}_\mu$.*

Доказательство. Пусть
$$A = \bigcup_{i=1}^{\infty} A_i \quad A_i \in \mathcal{C}_\mu$$

где $i \neq j => A_i \cap A_j = \emptyset$. Согласно теореме 2.3.11, $A \in \mathcal{C}_\mu$. Согласно теореме 2.3.2,

(2.3.30) $$\mu(A) \leq \sum_i \mu(A_i)$$

Согласно теореме 2.3.9, для любого $N$
$$\mu(A) \geq \mu\left(\bigcup_{i=1}^{N} A_i\right) = \sum_{i=1}^{N} \mu(A_i)$$

и следовательно

(2.3.31) $$\mu(A) \geq \sum_i \mu(A_i)$$

Равенство
$$\mu(A) = \sum_i \mu(A_i)$$

следует из (2.3.30), (2.3.31). $\square$

Теорема 2.3.13. *Пусть[2.25] $A \in \mathcal{C}_\mu$. Тогда существует множество $B$ такое, что*

(2.3.32) $$A \subseteq B$$

(2.3.33) $$\mu(A) = \mu(B)$$

(2.3.34) $$B = \bigcap_n B_n$$

(2.3.35) $$B_1 \supseteq B_2 \supseteq ... \supseteq B_n \supseteq ...$$

(2.3.36) $$B_n = \bigcup_k B_{nk}$$

(2.3.37) $$B_{nk} \in \mathcal{R}(\mathcal{S})$$

(2.3.38) $$\mu(B_{nk}) < \mu(A) + \frac{1}{n}$$

(2.3.39) $$B_{n1} \subseteq B_{n2} \subseteq ... \subseteq B_{nk} \subseteq ...$$

---

[2.24]Смотри также теорему [1]-4 на странице 273.
[2.25]Смотри также лемму в [1] на странице 315.

Доказательство. Согласно определению 2.3.4 и теореме 2.1.11, для любого $n$ существует множество $C_n$ такое, что

$$A \subseteq C_n \tag{2.3.40}$$

$$\mu(C_n) < \mu(A) + \frac{1}{n} \tag{2.3.41}$$

$$C_n = \bigcup_r \Delta_{nr} \quad \Delta_{nr} \in \mathcal{S} \tag{2.3.42}$$

Положим

$$B_n = \bigcap_{k=1}^{n} C_k \tag{2.3.43}$$

Утверждение (2.3.35) следует из (2.3.43). Утверждение

$$A \subseteq B_n \tag{2.3.44}$$

следует из (2.3.40), (2.3.43). Утверждение (2.3.32) следует из (2.3.34), (2.3.44). Из (2.3.42), (2.3.43) следует, что

$$B_n = \bigcup_r \delta_{nr} \quad \delta_{nr} \in \mathcal{S} \tag{2.3.45}$$

Пусть

$$B_{nk} = \bigcup_{r=1}^{k} \delta_{nr} \tag{2.3.46}$$

Утверждение (2.3.39) следует из (2.3.46). Утверждение (2.3.36) следует из (2.3.45), (2.3.46). Утверждение (2.3.37) следует из (2.3.45), (2.3.46) и теоремы 2.1.11.

Из утверждения (2.3.32) и теоремы 2.2.6 следует, что

$$\mu(A) \leq \mu(B) \tag{2.3.47}$$

Из утверждений (2.3.34), (2.3.41), (2.3.43) и теоремы 2.2.6 следует, что

$$\mu(B) \leq \mu(B_n) \leq \mu(C_n) < \mu(A) + \frac{1}{n} \tag{2.3.48}$$

Так как $n$ произвольно, то утверждение

$$\mu(B) \leq \mu(A) \tag{2.3.49}$$

является следствием утверждения (2.3.48). Утверждение (2.3.33) является следствием утверждений (2.3.47), (2.3.49). Утверждение (2.3.38) является следствием утверждений (2.3.36), (2.3.48) и теоремы 2.2.6. $\square$

# Измеримое отображение в абелеву $\Omega$-группу

## 3.1. Измеримое отображение

Определение 3.1.1. Минимальная $\sigma$-алгебра $\mathcal{B}(A)$ над совокупностью всех открытых шаров нормированной $\Omega$-группы $A$, называется **алгеброй Бореля**.[3.1] Множество, принадлежащее алгебре Бореля, называется **борелевским множеством** или $B$-**множеством**. □

Определение 3.1.2. Пусть $\mathcal{C}_X$ - $\sigma$-алгебра множеств множества $X$. Пусть $\mathcal{C}_Y$ - $\sigma$-алгебра множеств множества $Y$. Отображение

$$f : X \to Y$$

называется $(\mathcal{C}_X, \mathcal{C}_Y)$-**измеримым**,[3.2] если для всякого множества $C \in \mathcal{C}_Y$

$$f^{-1}(C) \in \mathcal{C}_X$$

□

Пример 3.1.3. Пусть на множестве $X$ определена $\sigma$-аддитивная мера $\mu$. Пусть $\mathcal{C}_\mu$ - $\sigma$-алгебра измеримых относительно меры $\mu$ множеств. Пусть $\mathcal{B}(A)$ - алгебра Бореля нормированной $\Omega$-группы $A$. Отображение

$$f : X \to A$$

называется $\mu$-**измеримым**,[3.3] если для всякого множества $C \in \mathcal{B}(A)$

$$f^{-1}(C) \in \mathcal{C}_\mu$$

□

Пример 3.1.4. Пусть $\mathcal{B}(A)$ - алгебра Бореля нормированной $\Omega_1$-группы $A$. Пусть $\mathcal{B}(B)$ - алгебра Бореля нормированной $\Omega_2$-группы $B$. Отображение

$$f : A \to B$$

называется **борелевским**,[3.4] если для всякого множества $C \in \mathcal{B}(B)$

$$f^{-1}(C) \in \mathcal{B}(A)$$

□

---

[3.1]Смотри определение в [1], с. 46. Согласно замечанию 2.1.3, алгебра Бореля может быть также порождена множеством замкнутых шаров.

[3.2]Смотри аналогичное определение в [1], с. 282.

[3.3]Смотри аналогичное определение в [1], с. 282, определение 1. Если мера $\mu$ на множестве $X$ определена в контексте, мы также будем называть отображение

$$f : X \to A$$

**измеримым**.

[3.4]Смотри аналогичное определение в [1], с. 282, определение 1.

Теорема 3.1.5. *Пусть $\mathcal{C}_X$ - $\sigma$-алгебра множеств множества $X$. Пусть $\mathcal{C}_Y$ - $\sigma$-алгебра множеств множества $Y$. Пусть $\mathcal{C}_Z$ - $\sigma$-алгебра множеств множества $Z$.*

3.1.5.1: *Пусть отображение*

$$f : X \to Y$$

*$(\mathcal{C}_X, \mathcal{C}_Y)$-измеримо.*

3.1.5.2: *Пусть отображение*

$$g : Y \to Z$$

*$(\mathcal{C}_Y, \mathcal{C}_Z)$-измеримо.*

*Тогда отображение*

$$g \circ f : X \to Z$$

*$(\mathcal{C}_X, \mathcal{C}_Z)$-измеримо.*[3.5]

Доказательство. Пусть $A \in \mathcal{C}_Z$. Согласно определению 3.1.2 и утверждению 3.1.5.2,

$$g^{-1}(A) \in \mathcal{C}_Y$$

Согласно определению 3.1.2 и утверждению 3.1.5.1,

$$f^{-1}(g^{-1}(A)) = (gf)^{-1}(A) \in \mathcal{C}_X$$

Следовательно, отображение

$$g \circ f : X \to Z$$

$(\mathcal{C}_X, \mathcal{C}_Z)$-измеримо. □

### 3.2. Простое отображение

Определение 3.2.1. Пусть на множестве $X$ определена $\sigma$-аддитивная мера $\mu$. Отображение

$$f : X \to A$$

в нормированную $\Omega$-группу $A$ называется **простым отображением**, если это отображение $\mu$-измеримо и принимает не более, чем счётное множество значений. □

Теорема 3.2.2. *Пусть отображение*

$$f : X \to A$$

*принимает не более, чем счётное множество значений $y_1, y_2, \ldots$ . Отображение $f$ $\mu$-измеримо тогда и только тогда, когда все множества*

$$A_n = \{x : f(x) = y_n\}$$

*$\mu$-измеримы.*[3.6]

---

[3.5]Смотри аналогичную теорему в [1], страницы 282, 283, теорема 1.

[3.6]Смотри аналогичную теорему в [1], страница 292, теорема 1.

Доказательство. Каждое множество $\{y_n\}$ является борелевским множеством. Так как $A_n$ является прообразом множества $\{y_n\}$, то $A_n$ является $\mu$-измеримым, если отображение $f$ является $\mu$-измеримым. Следовательно, условие теоремы необходимо.

Пусть все множества $A_n$ $\mu$-измеримы. Прообраз $f^{-1}(B)$ борелевского множества $B \subset A$ $\mu$-измерим, так как он является объединением

$$\bigcup_{y_n \in B} A_n$$

не более чем счётного множества $\mu$-измеримых множеств $A_n$. Следовательно, отображение $f$ $\mu$-измеримо и условие теоремы достаточно. $\square$

Теорема 3.2.3. *Пусть на множестве $X$ определена $\sigma$-аддитивная мера $\mu$. Пусть*

$$f : X \to A$$

*простое отображения в нормированную $\Omega$-группу $A$. Пусть отображение $f$ $\mu$-измеримо на множестве $X_i \subset X$, $i = 1, 2, ...$ . Тогда отображение $f$ $\mu$-измеримо на множестве $\bigcup_i X_i$ .*

Доказательство. Пусть $y_1, y_2, ...$ - область значений отображения $f$. Согласно теореме 3.2.2, множество

$$Y_{i \cdot k} = \{x \in X_i : f(x) = y_k\}$$

$\mu$-измеримо. Согласно определениям 2.1.13, 2.2.10, множество

$$Y_k = \bigcup_i Y_{i \cdot k} = \{x \in \bigcup_i X_i : f(x) = y_k\}$$

$\mu$-измеримо. Согласно теореме 3.2.2, отображение $f$ $\mu$-измеримо на множестве $\bigcup_i X_i$ . $\square$

Теорема 3.2.4. *Пусть на множестве $X$ определена $\sigma$-аддитивная мера $\mu$. Пусть*

$$f : X \to A$$
$$g : X \to A$$

*простые отображения в нормированную $\Omega$-группу $A$. Тогда отображение*

$$h = f + g$$

*является простым отображением.*[3.7]

Доказательство. Согласно определению 3.2.1, простые отображения $f$ и $g$ имеют конечные или счётные области значений. Пусть $y_1, y_2, ...$ - область значений отображения $f$. Пусть $z_1, z_2, ...$ - область значений отображения $g$. Тогда область значений отображения $h$ состоит из значений

$$c_{ij} = y_i + z_j$$

и является конечным или счётным множеством. Для каждого $c_{ij}$ множество

$$\{x : h(x) = c_{ij}\} = \bigcup_{y_i + z_j = c_{ij}} \{x : f(x) = y_i\} \cap \{x : g(x) = z_j\}$$

---

[3.7]Смотри аналогичную теорему в [1], страница 283, теорема 3.

$\mu$-измеримо. Следовательно, отображение $h$ является простым отображением. $\square$

Теорема 3.2.5. *Пусть на множестве $X$ определена $\sigma$-аддитивная мера $\mu$. Пусть $\omega \in \Omega$ - $n$-арная операция. Пусть*

$$f_i : X \to A \quad i = 1, ..., n$$

*простые отображения в нормированную $\Omega$-группу $A$. Тогда отображение*

$$h = f_1 ... f_n \omega$$

*является простым отображением.*

Доказательство. Согласно определению 3.2.1, простые отображения $f_i$, $i = 1, ..., n,$ имеют конечные или счётные области значений. Пусть $y_{i \cdot 1}$, $y_{i \cdot 2}$, ... - область значений отображения $f_i$. Тогда область значений отображения $h$ состоит из значений

$$y_{i_1 ... i_n} = y_{1 \cdot i_1} ... y_{n \cdot i_n} \omega$$

и является конечным или счётным множеством. Для каждого $y_{i_1 ... i_n}$ множество

$$\{x : h(x) = y_{i_1 ... i_n}\} = \bigcup_{y_{1 \cdot i_1} ... y_{n \cdot i_n} \omega = y_{i_1 ... i_n}} \bigcap_{j=1}^{n} \{x : f_j(x) = y_{j \cdot j_i}\}$$

$\mu$-измеримо. Следовательно, отображение $h$ является простым отображением. $\square$

Теорема 3.2.6. *Пусть*

$$f : A_1 \relbar\joinrel\ast\joinrel\rightarrow A_2$$

*представление $\Omega_1$-группы $A_1$ с нормой $\|x\|_1$ в $\Omega_2$-группе $A_2$ с нормой $\|x\|_2$. Пусть на множестве $X$ определена $\sigma$-аддитивная мера $\mu$. Пусть*

$$g_i : X \to A_i \quad i = 1, 2$$

*простое отображение. Тогда отображение*

$$h = f_X(g_1)(g_2)$$

*является простым отображением.*

Доказательство. Согласно определению 3.2.1, простые отображения $g_1$ и $g_2$ имеют конечные или счётные области значений. Пусть $y_{1 \cdot 1}$, $y_{1 \cdot 2}$, ... - область значений отображения $g_1$. Пусть $y_{2 \cdot 1}$, $y_{2 \cdot 2}$, ... - область значений отображения $g_2$. Тогда область значений отображения $h$ состоит из значений

$$c_{ij} = f(y_{1 \cdot i})(y_{2 \cdot j})$$

и является конечным или счётным множеством. Для каждого $c_{ij}$ множество

$$\{x : h(x) = c_{ij}\} = \bigcup_{f(y_{1 \cdot i})(y_{2 \cdot j}) = c_{ij}} \{x : g_1(x) = y_{1 \cdot i}\} \cap \{x : g_2(x) = y_{2 \cdot j}\}$$

$\mu$-измеримо. Следовательно, отображение $h$ является простым отображением. $\square$

### 3.3. Действия над измеримыми отображениями

Теорема 3.3.1. *Пусть $A$ - нормированная $\Omega$-группа. Пусть $\{f_n\}$ - последовательность $\mu$-измеримых отображений*

$$f_n : X \to A$$

*Пусть*

$$f : X \to A$$

*такое отображение, что*

$$(3.3.1) \qquad f(x) = \lim_{n \to \infty} f_n(x)$$

*для каждого $x$. Тогда отображение $f$ является $\mu$-измеримым отображением.*

Доказательство. Докажем следующее равенство

$$(3.3.2) \qquad \{x : f(x) \in B_c(a, R)\} = \bigcup_k \bigcup_n \bigcap_{m>n} \{x : f_m(x) \in B_o(a, R + 1/k)\}$$

- Пусть $f(x) \in B_c(a, R)$. Согласно определению [4]-2.1.15

$$(3.3.3) \qquad \|f(x) - a\| \leq R$$

Из равенства (3.3.1) следует, что для любого $k$ существует такое $n$, что для $m > n$

$$(3.3.4) \qquad \|f(x) - f_m(x)\| < 1/k$$

Из равенств (3.3.3), (3.3.4) следует, что

$$(3.3.5) \qquad \begin{aligned} \|f_m(x) - a\| &= \|f_m(x) - f(x) + f(x) - a\| \\ &\leq \|f_m(x) - f(x)\| + \|f(x) - a\| \\ &< R + 1/k \end{aligned}$$

Из определения [4]-2.1.14 и равенства (3.3.5) следует, что для любого $k$ существует такое $n$, что для $m > n$

$$f_m(x) \in B_o(a, R + 1/k)$$

Следовательно, мы доказали, что

$$(3.3.6) \qquad \{x : f(x) \in B_c(a, R)\} \subseteq \bigcup_k \bigcup_n \bigcap_{m>n} \{x : f_m(x) \in B_o(a, R + 1/k)\}$$

- Пусть

$$x \in \bigcup_k \bigcup_n \bigcap_{m>n} \{x : f_m(x) \in B_o(a, R + 1/k)\}$$

Тогда существуют $k$, $n$ такие, что для $m > n$

$$(3.3.7) \qquad f_m(x) \in B_o(a, R + 1/k)$$

Из определения [4]-2.1.14 и равенства (3.3.7) следует, что

$$(3.3.8) \qquad \|f_m(x) - a\| < R + 1/k$$

Из теоремы [5]-2, страница 56, равенства (3.3.1) и неравенства (3.3.8) следует, что

$$(3.3.9) \qquad \|f(x) - a\| \leq R$$

Из определения [4]-2.1.15 и неравенства (3.3.9) следует, что $f(x) \in B_c(a, R)$. Следовательно, мы доказали, что

$$\bigcup_k \bigcup_n \bigcap_{m>n} \{x : f_m(x) \in B_o(a, R + 1/k)\} \subseteq \{x : f(x) \in B_c(a, R)\} \quad (3.3.10)$$

- Равенство (3.3.2) следует из равенств (3.3.6), (3.3.10).

Если отображения $f_n$ измеримы, то, согласно примеру 3.1.3, множества

$$\{x : f_m(x) \in B_o(a, R + 1/k)\}$$

измеримы. Так как совокупность измеримых множеств является $\sigma$-алгеброй, то из равенства (3.3.2) следует, что множество

$$\{x : f(x) \in B_c(a, R)\}$$

измеримо. Согласно примеру 3.1.3, отображение $f$ измеримо. $\square$

Теорема 3.3.2. *Пусть $A$ - нормированная $\Omega$-группа. Пусть множество значений отображения*

$$f : X \to A$$

*компактно. Отображение $f$ $\mu$-измеримо тогда и только тогда, когда оно может быть представлено как предел равномерно сходящейся последовательности простых отображений.*[3.8]

Доказательство. Пусть отображение $f$ является пределом равномерно сходящейся последовательности простых отображений. Согласно определению 3.2.1, простое отображение $\mu$-измеримо. Согласно теореме 3.3.1, предел последовательности простых отображений $\mu$-измерим. Следовательно, отображение $f$ $\mu$-измеримо.

Пусть отображение $f$ $\mu$-измеримо. Для заданного $n$, рассмотрим множество открытых шаров

$$B_f = \{B_o(y, 1/n) : \exists x, y = f(x)\}$$

Множество $B_f$ является открытым покрытием множества значений отображения $f$. Следовательно, существует конечное множество $B'_f$, $B'_f \subset B_f$, являющееся открытым покрытием множества значений отображения $f$. Для $x \in X$, выберем открытый шар

$$B_o(y, 1/n) \in B_f \quad f(x) \in B_o(y, 1/n) \quad (3.3.11)$$

и мы положим

$$f_n(x) = y \quad (3.3.12)$$

Отображение $f_n$ - простое. Из равенств (3.3.11), (3.3.12) и определения [4]-2.1.14 следует, что

$$\|f_n(x) - f(x)\| < \frac{1}{n} \quad (3.3.13)$$

Из равенства (3.3.13) следует, что последовательность простых отображений $f_n$ равномерно сходится к отображению $f$. $\square$

---

[3.8] Смотри аналогичную теорему в [1], страница 292, теорема 2.

Теорема 3.3.3. *Пусть на множестве $X$ определена $\sigma$-аддитивная мера $\mu$. Пусть $A$ - нормированная $\Omega$-группа. Пусть множество значений отображения*
$$f : X \to A$$
*компактно. Пусть отображение $f$ $\mu$-измеримо на множестве $X_i \subset X$, $i = 1$, ..., $n$,*
$$i \neq j \implies X_i \cap X_j = \emptyset$$
*Тогда отображение $f$ $\mu$-измеримо на множестве $\bigcup_i X_i$.*

Доказательство. Согласно теореме 3.3.2, для каждого $i$ существует последовательность простых отображений
$$f_{i \cdot k} : X_i \to A$$
равномерно сходящейся к отображению $f$ на множестве $X_i$. Для каждого $k$, рассмотрим отображение
$$f_k : \bigcup_i X_i \to A$$
определённое правилом
$$x \in X_i \implies f_k(x) = f_{i \cdot k}(x)$$
Так как отображение $f_k$ принимает не более, чем счётное множество значений, на множестве $X_i$, то отображение $f_k$ принимает не более, чем счётное множество значений, на множестве $\bigcup_i X_i$. Согласно теореме 3.2.3, отображение $f_k$ измеримо на множестве $\bigcup_i X_i$.

Так как последовательность отображений $f_k$ равномерно сходится к отображению $f$ на множестве $X_i$, то, согласно определению [4]-2.6.3, для заданного $\epsilon \in R$, $\epsilon > 0$, существует $K_i$ такое, что из условия $k > K_i$ следует
$$\|f_k(x) - f(x)\| < \epsilon$$
для любого $x \in X_i$. Пусть
$$K = \max(K_1, ..., K_n)$$
Тогда для заданного $\epsilon \in R$, $\epsilon > 0$, из условия $k > K$ следует
$$\|f_k(x) - f(x)\| < \epsilon$$
для любого $x \in \bigcup_i X_i$. Согласно определению [4]-2.6.3, последовательность отображений $f_k$ равномерно сходится к отображению $f$ на множестве $\bigcup_i X_i$. Согласно теореме 3.3.2, отображение $f$ $\mu$-измеримо на множестве $\bigcup_i X_i$. □

Теорема 3.3.4. *Пусть на множестве $X$ определена $\sigma$-аддитивная мера $\mu$. Пусть $A$ - полная $\Omega$-группа. Пусть*
$$f : X \to A$$
$$g : X \to A$$

$\mu$-измеримые отображения с компактным множеством значений. Тогда отображение

(3.3.14) $$h = f + g$$

является $\mu$-измеримым отображением с компактным множеством значений.[3.9]

Доказательство. Согласно теореме 3.3.2, существует последовательность простых отображений $f_n$ равномерно сходящихся к отображению $f$ и существует последовательность простых отображений $g_n$ равномерно сходящихся к отображению $g$. Для каждого $n$, отображение

$$h_n = f_n + g_n$$

является простым отображением согласно теореме 3.2.4. Согласно теореме [4]-2.6.7 последовательность отображений $h_n$ равномерно сходится к отображению $h$. Согласно теореме 3.3.1, отображению $h$ является $\mu$-измеримым отображением.

Рассмотри открытое покрытие $O_h$ множества значений отображения $h$. Согласно определению [4]-2.2.2, открытое покрытие $O_h$ содержит множество открытых шаров

$$B_h = \{B_o(y, 2\epsilon_x) : \exists x, y = h(x)\}$$

Рассмотрим множество открытых шаров

$$B_f = \{B_o(y, \epsilon_x) : \exists x, y = f(x)\}$$

Множество $B_f$ является открытым покрытием множества значений отображения $f$. Следовательно, существует конечное множество $B'_f$, $B'_f \subset B_f$, являющееся открытым покрытием множества значений отображения $f$. Положим

$$I_f = \{x \in X : y = f(x), B_o(y, \epsilon_x) \in B'_f\}$$

Рассмотрим множество открытых шаров

$$B_g = \{B_o(y, \epsilon_x) : \exists x, y = g(x)\}$$

Множество $B_g$ является открытым покрытием множества значений отображения $g$. Следовательно, существует конечное множество $B'_g$, $B'_g \subset B_g$, являющееся открытым покрытием множества значений отображения $g$. Положим

$$I_g = \{x \in X : y = g(x), B_o(y, \epsilon_x) \in B'_g\}$$

Положим

$$I_h = I_f \cup I_g$$

Для любого $x \in X$, существует $x' \in I_h$ такой, что

(3.3.15) $$f(x) \in B_o(f(x'), \epsilon_{x'}) \quad g(x) \in B_o(g(x'), \epsilon_{x'})$$

Согласно теореме [4]-2.4.2, утверждение

$$h(x) \in B_o(h(x'), 2\epsilon_{x'})$$

следует из (3.3.14), (3.3.15). Следовательно, множество открытых шаров

$$B'_h = \{B_o(y, 2\epsilon_x) : \exists x \in I_h, y = h(x)\} \subseteq B_h \subseteq O_h$$

---

[3.9]Смотри аналогичную теорему в [1], страница 283, теорема 3.

является конечным открытым покрытием множества значений отображения $h$. Согласно определению [4]-2.2.15, множество значений отображения $h$ компактно. □

Теорема 3.3.5. *Пусть на множестве $X$ определена $\sigma$-аддитивная мера $\mu$. Пусть $\omega \in \Omega$ - $n$-арная операция. Пусть*

$$f_i : X \to A \quad i = 1, ..., n$$

*$\mu$-измеримые отображения в полную $\Omega$-группу $A$. Тогда отображение*

$$h = f_1...f_n\omega$$

*является $\mu$-измеримым отображением.*

Доказательство. Согласно теореме 3.3.2, существует последовательность простых отображений $f_{i\cdot m}$ равномерно сходящихся к отображению $f_i$. Для каждого $m$, отображение

$$h_m = f_{1\cdot m}...f_{n\cdot m}\omega$$

является простым отображением согласно теореме 3.2.5. Согласно теореме [4]-2.6.8 последовательность отображений $h_n$ равномерно сходится к отображению $h$. □

Теорема 3.3.6. *Пусть на множестве $X$ определена $\sigma$-аддитивная мера $\mu$. Пусть*

$$f : A_1 \relbar\joinrel\ast\joinrel\rightarrow A_2$$

*представление $\Omega_1$-группы $A_1$ с нормой $\|x\|_1$ в $\Omega_2$-группе $A_2$ с нормой $\|x\|_2$. Пусть*

$$g_i : X \to A_i \quad i = 1, 2$$

*$\mu$-измеримое отображение. Тогда отображение*

$$h = f_X(g_1)(g_2)$$

*является $\mu$-измеримым отображением.*

Доказательство. Согласно теореме 3.3.2, существует последовательность простых отображений $f_{i\cdot m}$ равномерно сходящихся к отображению $f_i$. Для каждого $m$, отображение $f_X(g_{1\cdot n})(g_{2\cdot n})$ является простым отображением согласно теореме 3.2.6. Согласно теореме [4]-3.2.2 последовательность отображений $h_n$ равномерно сходится к отображению $h$. □

### 3.4. Сходимость почти всюду

Определение 3.4.1. Пусть на множестве $X$ определена $\sigma$-аддитивная мера $\mu$. Последовательность[3.10]

$$f_n : X \to A$$

$\mu$-измеримых отображений в $\Omega$-группу $A$ **сходится почти всюду**, если

(3.4.1) $$f(x) = \lim_{n \to \infty} f_n(x)$$

для почти всех $x \in X$, т. е. множество $x$, в которых (3.4.1) не выполняется, имеет меру нуль. □

---

[3.10]Смотри так же определение [1]-2 на странице 286.

Теорема 3.4.2. *Пусть на множестве $X$ определена $\sigma$-аддитивная мера $\mu$. Если последовательность*
$$f_n : X \to A$$
*$\mu$-измеримых отображений в $\Omega$-группу $A$ сходится к отображению*
$$f : X \to A$$
*почти всюду, то отображение $f$ также $\mu$-измеримо.*

Доказательство. Пусть
$$A = \{x \in X : \lim_{n \to \infty} f_n(x) = f(x)\} \quad B = X \setminus A$$
Согласно определению 3.4.1, $\mu(B) = 0$.

Лемма 3.4.3. *Отображение $f$ $\mu$-измеримо на множестве $B$.*

Доказательство. Утверждение является следствием определения 3.1.2 и теоремы 2.2.8. $\odot$

Лемма 3.4.4. *Отображение $f$ $\mu$-измеримо на множестве $A$.*

Доказательство. Согласно замечанию 2.1.3, множество $A$ $\mu$-измеримо. Согласно теореме 3.3.1, отображение $f$ $\mu$-измеримо на множестве $A$. $\odot$

Из теоремы 3.3.3 и из лемм 3.4.3, 3.4.4 следует, что отображение $f$ $\mu$-измеримо. $\square$

Теорема 3.4.5 (Дмитрий Фёдорович Егоров). *Пусть[3.11] последовательность $\mu$-измеримых отображений*
$$f_n : X \to A$$
*сходится на измеримом множестве $E$ почти всюду к отображению $f$. Тогда, для любого $\delta > 0$, существует измеримое множество $E_\delta \subset E$ такое, что*
3.4.5.1: $\mu(E_\delta) > \mu(E) - \delta$
3.4.5.2: *Последовательность $f_n$ равномерно сходится на множестве $E_\delta$.*

Доказательство. Согласно теореме 3.4.2, отображение $f$ $\mu$-измеримо на множестве $E$. Пусть

$$(3.4.2) \qquad E_n^m = \bigcap_{i \geq n} \left\{ x \in E : \|f_i(x) - f(x)\| < \frac{1}{m} \right\}$$

$$(3.4.3) \qquad E^m = \bigcup_{n=1}^{\infty} E_n^m$$

Утверждение

$$(3.4.4) \qquad E_1^m \subset E_2^m \subset ...$$

следует из (3.4.2). Из утверждений (3.4.3), (3.4.4) и теоремы 2.2.13 следует, что для любого $m$ и для любого $\delta > 0$, существует $n_0(m)$ такое, что

$$(3.4.5) \qquad \mu(E^m - E_{n_0(m)}^m) < \frac{\delta}{2^m}$$

---

[3.11]Смотри так же теорему [1]-6 на странице 287.

Положим

$$E_\delta = \bigcap_{m=1}^{\infty} E_{n_0(m)}^m \tag{3.4.6}$$

Если $x \in E_\delta$, то, для любого $m = 1, 2, ...$ и для любого $i > n_0(m)$, неравенство

$$\|f_i(x) - f(x)\| < \frac{1}{m} \tag{3.4.7}$$

следует из (3.4.2), (3.4.6). Утверждение 3.4.5.2 следует из неравенства (3.4.7) и определения [4]-2.6.3.

Если $x \in E \setminus E^m$, то, из (3.4.2), (3.4.3) следует, что существуют сколь угодно большие значения $i$, для которых

$$\|f_i(x) - f(x)\| > \frac{1}{m} \tag{3.4.8}$$

Следовательно, последовательность $f_n(x)$ в точке $x \in E \setminus E^m$ не сходится к $f(x)$. Так как $f_n$ сходится к $f$ почти всюду, то из теоремы 2.2.8 следует, что

$$\mu(E \setminus E^m) = 0 \tag{3.4.9}$$

Неравенство

$$\mu(E - E_{n_0(m)}^m) = \mu(E^m - E_{n_0(m)}^m) < \frac{\delta}{2^m} \tag{3.4.10}$$

следует из (3.4.5), (3.4.9). Утверждение 3.4.5.1 является следствием неравенства

$$\mu(E \setminus E_\delta) = \mu\left(E \setminus \bigcap_{m=1}^{\infty} E_{n_0(m)}^m\right) = \mu\left(\bigcup_{m=1}^{\infty} (E \setminus E_{n_0(m)}^m)\right)$$
$$\leq \sum_{m=1}^{\infty} \mu(E \setminus E_{n_0(m)}^m) < \sum_{m=1}^{\infty} \frac{\delta}{2^m} = \delta$$

которое следует из неравенства (3.4.10). $\square$

Глава 4

# Интеграл отображения в абелеву $\Omega$-группу

Пусть на множестве $X$ определена $\sigma$-аддитивная мера $\mu$. Пусть определено эффективное представление поля действительных чисел $R$ в полной абелевой $\Omega$-группе $A$.[4.1] В этой главе мы рассмотрим определение интеграл Лебега $\mu$-измеримого отображения

$$f : X \to A$$

## 4.1. Интеграл простого отображения

Определение 4.1.1. Пусть $a_i$ - последовательность $A$-чисел. Если

$$\sum_{i=1}^{\infty} \|a^i\| < \infty$$

то мы будем говорить, что **ряд**

$$\sum_{i=1}^{\infty} a^i$$

**сходится нормально**.[4.2] $\square$

Определение 4.1.2. Для простого отображения

$$f : X \to A$$

рассмотрим ряд

(4.1.1) $$\sum_n \mu(F_n) f_n$$

где

- Множество $\{f_1, f_2, ...\}$ является областью определения отображения $f$
- Если $n \neq m$, то $f_n \neq f_m$
- $F_n = \{x \in X : f(x) = f_n\}$

Простое отображение

$$f : X \to A$$

называется **интегрируемым** по множеству $X$, если ряд (4.1.1) сходится нормально.[4.3] Если отображение $f$ интегрируемо, то сумма ряда (4.1.1) называется **интегралом отображения** $f$ по множеству $X$

(4.1.2) $$\int_X d\mu(x) f(x) = \sum_n \mu(F_n) f_n$$

---

[4.1]Другими словами, $\Omega$-группа $A$ является $R$-векторным пространством.

[4.2]Смотри также определение нормальной сходимости ряда на странице [6]-12.

[4.3]Смотри аналогичное определение в [1], определение 2, с. 293.

$\square$

Теорема 4.1.3. *Пусть*
$$f : X \to A$$
*простое отображение. Пусть $f$ принимает значение[4.4] $f_n$ на множестве $F_n \subset X$. Пусть $X = \bigcup\limits_n F_n$, $F_n \cap F_m = \emptyset$. Отображение $f$ интегрируемо тогда и только тогда, когда ряд*

$$\sum_n \mu(F_n) f_n \tag{4.1.3}$$

*сходится по норме. Тогда*

$$\int_X d\mu(x) f(x) = \sum_n \mu(F_n) f_n \tag{4.1.4}$$

Доказательство. Пусть
$$X_i = \{x \in X : f(x) = f_i\}$$
Тогда
$$X_i = \bigcup_{f(X_i) = f(F_n)} F_n \tag{4.1.5}$$

Из (4.1.5) следует, что

$$\mu(X_i) = \sum_{f(X_i) = f(F_n)} \mu(F_n) \tag{4.1.6}$$

Так как ряд (4.1.3) сходится по норме, то

$$\sum_n \mu(F_n) f_n = \sum_i \left( \sum_{f(X_i) = f(F_n)} \mu(F_n) \right) f(X_i) = \sum_i \mu(X_i) f(X_i) \tag{4.1.7}$$

следует из (4.1.6). (4.1.4) следует из (4.1.2), (4.1.7). $\square$

Теорема 4.1.4. *Пусть*
$$f : X \to A$$
$$g : X \to A$$
*простые отображения.[4.5] Если существуют интегралы*
$$\int_X d\mu(x) f(x)$$
$$\int_X d\mu(x) g(x)$$
*то существует интеграл*
$$\int_X d\mu(x)(f(x) + g(x))$$
*и*

$$\int_X d\mu(x)(f(x) + g(x)) = \int_X d\mu(x) f(x) + \int_X d\mu(x) g(x) \tag{4.1.8}$$

---

[4.4]Мы не требуем $f_n \neq f_m$ при условии $n \neq m$. Однако мы требуем $F_n \cap F_m = \emptyset$. Смотри также лемма в [1], с. 293.

[4.5]Смотри аналогичное утверждение в [1], свойство A, с. 294.

Доказательство. Пусть $f$ принимает значение $f_n$ на множестве $F_n \subset X$. Пусть

$$F_n \cap F_m = \emptyset \qquad (4.1.9)$$

Пусть $g$ принимает значение $g_k$ на множестве $G_k \subset X$. Пусть

$$G_k \cap G_l = \emptyset \qquad (4.1.10)$$

Равенство

$$\mu(F_n) = \sum_k \mu(F_n \cap G_k) \qquad (4.1.11)$$

следует из равенства

$$F_n = \bigcup_k F_n \cap G_k$$

и условия (4.1.10). Равенство

$$\mu(G_k) = \sum_n \mu(F_n \cap G_k) \qquad (4.1.12)$$

следует из равенства

$$G_k = \bigcup_n F_n \cap G_k$$

и условия (4.1.9). Условие

$$(F_n \cap G_k) \cap (F_m \cap G_l) = \emptyset$$

следует из условий (4.1.9), (4.1.10).

4.1.4.1: Так как отображение $f$ интегрируемо, то равенство

$$\begin{aligned}\int_X d\mu(x) f(x) &= \sum_n \mu(F_n) f_n = \sum_n \left( \sum_k \mu(F_n \cap G_k) \right) f_n \\ &= \sum_n \sum_k \mu(F_n \cap G_k) f_n\end{aligned} \qquad (4.1.13)$$

следует из (4.1.4), (4.1.11) и ряд

$$\sum_n \sum_k \mu(F_n \cap G_k) f_n$$

сходится по норме.

4.1.4.2: Так как отображение $g$ интегрируемо, то равенство

$$\begin{aligned}\int_X d\mu(x) g(x) &= \sum_k \mu(G_k) g_k = \sum_k \left( \sum_n \mu(G_k \cap F_n) \right) g_k \\ &= \sum_k \sum_n \mu(G_k \cap F_n) g_k \\ &= \sum_n \sum_k \mu(F_n \cap G_k) g_k\end{aligned} \qquad (4.1.14)$$

следует из (4.1.4), (4.1.12) и ряд

$$\sum_k \sum_n \mu(F_n \cap G_k) g_k = \sum_n \sum_k \mu(F_n \cap G_k) g_k$$

сходится по норме.

Из утверждений 4.1.4.1, 4.1.4.2 следует, что

$$
\begin{aligned}
&\int_X d\mu(x)f(x) + \int_X d\mu(x)g(x) \\
&= \sum_n \sum_k \mu(F_n \cap G_k)f_n + \sum_n \sum_k \mu(F_n \cap G_k)g_k \\
&= \sum_n \sum_k \mu(F_n \cap G_k)(f_n + g_k)
\end{aligned}
\tag{4.1.15}
$$

и ряд

$$\sum_n \sum_k \mu(F_n \cap G_k)(f_n + g_k)$$

сходится по норме. Следовательно, согласно теореме 4.1.3, равенство (4.1.8) следует из (4.1.15). □

Теорема 4.1.5. *Пусть*

$$f : X \to A$$

*простое отображение. Интеграл*

$$\int_X d\mu(x)f(x)$$

*существует тогда и только тогда, когда интеграл*

$$\int_X d\mu(x)\|f(x)\|$$

*существует. Тогда*

$$\left\| \int_X d\mu(x)f(x) \right\| \leq \int_X d\mu(x)\|f(x)\| \tag{4.1.16}$$

Доказательство. Равенство

$$\sum_n \|\mu(F_n)f_n\| = \sum_n \mu(F_n)\|f_n\| \tag{4.1.17}$$

является следствием равенства

$$\|\mu(F_n)f_n\| = \mu(F_n)\|f_n\|$$

Следовательно, ряд

$$\sum_n \mu(F_n)f_n$$

сходится по норме тогда и только тогда, когда сходится ряд

$$\sum_n \mu(F_n)\|f_n\|$$

Неравенство (4.1.16) следует из неравенства

$$\left\| \sum_n \mu(F_n)f_n \right\| \leq \sum_n \|\mu(F_n)f_n\|$$

равенства (4.1.17) и определения интеграла. □

Теорема 4.1.6. *Пусть $\omega \in \Omega$ - n-арная операция в абелевой $\Omega$-группе $A$. Пусть простое отображение*

$$f_i : X \to A \quad i = 1, ..., n$$

*интегрируемо. Тогда отображение*

$$h = f_1...f_n\omega$$

*интегрируемо и*

$$\left\| \int_X d\mu(x) h(x) \right\| \leq \int_X d\mu(x)(\|\omega\| \|f_1(x)\|...\|f_n(x)\|)$$

Доказательство. Теорема является следствием теоремы 4.1.5 и неравенства [4]-(2.1.6). □

Теорема 4.1.7. *Пусть на множестве $X$ определена $\sigma$-аддитивная мера $\mu$. Пусть*

$$f : A_1 \relbar\joinrel\mathrel{*\!}\joinrel\rightarrow A_2$$

*представление $\Omega_1$-группы $A_1$ с нормой $\|x\|_1$ в $\Omega_2$-группе $A_2$ с нормой $\|x\|_2$. Пусть*

$$g_i : X \to A_i \quad i = 1, 2$$

*простое интегрируемое отображение. Тогда отображение*

$$h = f_X(g_1)(g_2)$$

*интегрируемо и*

$$\left\| \int_X d\mu(x) h(x) \right\|_2 \leq \int_X d\mu(x)(\|f\| \|g_1(x)\|_1 \|g_2(x)\|_2)$$

Доказательство. Теорема является следствием теоремы 4.1.5 и неравенства [4]-(3.1.2). □

Теорема 4.1.8. *Пусть на множестве $X$ определена $\sigma$-аддитивная мера $\mu$. Пусть*

$$f : A_1 \relbar\joinrel\mathrel{*\!}\joinrel\rightarrow A_2$$

*эффективное представление $\Omega_1$-группы $A_1$ с нормой $\|x\|_1$ в $\Omega_2$-группе $A_2$ с нормой $\|x\|_2$. Пусть преобразования представления $f$ являются автоморфизмами представления*

$$R \relbar\joinrel\mathrel{*\!}\joinrel\rightarrow A_2$$

*Пусть*

$$g_2 : X \to A_2$$

*простое интегрируемое отображение в нормированную $\Omega$-группу $A_2$. Тогда отображение*

$$h = a_1 g_2$$

*интегрируемо и*

$$(4.1.18) \qquad \int_X d\mu(x)(a_1 g_2(x)) = a_1 \int_X d\mu(x) g_2(x)$$

Доказательство. Пусть $g_2$ принимает значение $g_{2\cdot k}$ на множестве $G_k \subset X$. Пусть
$$G_k \cap G_l = \emptyset$$
Так как отображение $g_2$ интегрируемо, то равенство

(4.1.19) $$\int_X d\mu(x) g_2(x) = \sum_k \mu(G_k) g_{2\cdot k}$$

следует из (4.1.4) и ряд в (4.1.19) сходится по норме. Согласно теореме [4]-3.1.2,
$$\|a_1 g_{2\cdot j}\|_2 \le \|f\| \|a_1\|_1 \|g_{2\cdot j}\|_2$$
Следовательно, ряд
$$\sum_k \mu(G_k)(a_1 g_{2\cdot k})$$
сходится по норме. Так как преобразование, порождённое $a_1 \in A_1$, является автоморфизмом представления[4.6]
$$R \longrightarrow_* A_2$$
равенство (4.1.18) следует из равенства
$$\sum_k \mu(G_k)(a g_{2\cdot k}) = \sum_k a(\mu(G_k) g_{2\cdot k}) = a \sum_k \mu(G_k) g_{2\cdot k}$$
и теоремы 4.1.3. $\square$

Теорема 4.1.9. *Пусть на множестве $X$ определена $\sigma$-аддитивная мера $\mu$. Пусть*
$$f : A_1 \longrightarrow_* A_2$$
*эффективное представление $\Omega_1$-группы $A_1$ с нормой $\|x\|_1$ в $\Omega_2$-группе $A_2$ с нормой $\|x\|_2$. Пусть преобразования представления*
$$R \longrightarrow_* A_1$$
*являются автоморфизмами представления $f$. Пусть*
$$g_1 : X \to A_1$$
*простое интегрируемое отображение в нормированную $\Omega$-группу $A_1$. Тогда отображение*
$$h = g_1 a_2$$
*интегрируемо и*

(4.1.20) $$\int_X d\mu(x)(g_1(x) a_2) = \left( \int_X d\mu(x) g_1(x) \right) a_2$$

---

[4.6]Для нас здесь важно, что для любых $p \in R$, $a_1 \in A_1$, $a_2, b_2 \in A_2$,
$$p(a_1 a_2) = a_1(p a_2)$$
$$a_1(a_2 + b_2) = a_1 a_2 + a_1 b_2$$
Смотри также определение [2]-2.2.2.

ДОКАЗАТЕЛЬСТВО. Пусть $g_1$ принимает значение $g_{1\cdot k}$ на множестве $G_k \subset X$. Пусть
$$G_k \cap G_l = \emptyset$$
Так как отображение $g_1$ интегрируемо, то равенство
$$\int_X d\mu(x) g_1(x) = \sum_k \mu(G_k) g_{1\cdot k} \qquad (4.1.21)$$
следует из (4.1.4) и ряд в (4.1.21) сходится по норме. Согласно теореме [4]-3.1.2,
$$\|g_{1\cdot j} a_2\|_1 \leq \|f\| \|g_{1\cdot j}\|_1 \|a_2\|_2$$
Следовательно, ряд
$$\sum_k \mu(G_k)(g_{1\cdot k}) a_2$$
сходится по норме. Так как преобразование представления
$$R \longrightarrow_* A_1$$
является автоморфизмом представления $f$, равенство[4.7] (4.1.20) следует из равенства
$$\sum_k \mu(G_k)(g_{1\cdot k} a) = \sum_k (\mu(G_k) g_{1\cdot k}) a = \left( \sum_k \mu(G_k) g_{1\cdot k} \right) a$$
и теоремы 4.1.3. $\square$

ТЕОРЕМА 4.1.10. *Пусть простое отображение*
$$f : X \to A$$
*удовлетворяет условию*
$$\|f(x)\| \leq M \qquad (4.1.22)$$
*Если мера множества $X$ конечна, то*
$$\int_X d\mu(x) \|f(x)\| \leq M \mu(X) \qquad (4.1.23)$$

ДОКАЗАТЕЛЬСТВО. Пусть $f$ принимает значение $f_n$ на множестве $F_n \subset X$. Пусть $X = \bigcup_n F_n$, $F_n \cap F_m = \emptyset$. Согласно теореме 4.1.3, неравенство
$$\int_X d\mu(x) \|f(x)\| = \sum_n \mu(F_n) \|f_n\| \leq M \sum_n \mu(F_n) = M\mu(X) \qquad (4.1.24)$$
следует из равенства (4.1.4) и неравенства (4.1.22). Неравенство (4.1.23) следует из неравенства (4.1.24). $\square$

---

[4.7]Для нас здесь важно, что для любых $p \in R$, $a_1, b_1 \in A_1$, $a_2 \in A_2$,
$$p(a_1 a_2) = pf(a_1)(a_2) = f(pa_1)(a_2) = (pa_1) a_2$$
$$(a_1 + b_1) a_2 = a_1 a_2 + b_1 a_2$$
Смотри также определение [2]-2.2.2.

Теорема 4.1.11. *Пусть простое отображение*
$$f : X \to A$$
*удовлетворяет условию*
$$\|f(x)\| \le M$$
*Если мера множества $X$ конечна, то отображение $f$ интегрируемо и*

(4.1.25) $$\left\| \int_X d\mu(x) f(x) \right\| \le M\mu(X)$$

Доказательство. Интегрируемость отображения $f$ следует из теорем 4.1.5, 4.1.10. Неравенство (4.1.25) следует из неравенств (4.1.16), (4.1.23). □

## 4.2. Интеграл измеримого отображения на множестве конечной меры

Определение 4.2.1. $\mu$-измеримое отображение
$$f : X \to A$$
называется **интегрируемым** по множеству $X$,[4.8] если существует последовательность простых отображений
$$f_n : X \to A$$
сходящаяся равномерно к $f$. Если отображение $f$ интегрируемо, то предел

(4.2.1) $$\int_X d\mu(x) f(x) = \lim_{n \to \infty} \int_X d\mu(x) f_n(x)$$

называется **интегралом отображения** $f$ по множеству $X$. □

Теорема 4.2.2. *Пусть*
$$f : X \to A$$
*$\mu$-измеримое отображение.*[4.9] *Пусть мера множества $X$ конечна.*

4.2.2.1: *Для любой равномерно сходящейся последовательности $f_n$ простых интегрируемых отображений*
$$f_n : X \to A$$
*предел (4.2.1) существует.*

4.2.2.2: *Предел (4.2.1) не зависит от выбора последовательности $f_n$.*

4.2.2.3: *Для простого отображения, определение 4.2.1 сводится к определению 4.1.2.*

Доказательство. Согласно теореме [4]-2.6.5, так как последовательность $f_n$ сходится равномерно к отображению $f$, то для любого $\epsilon \in R$, $\epsilon > 0$, существует $N$ такое, что

(4.2.2) $$\|f_n(x) - f_m(x)\| < \frac{\epsilon}{\mu(X)}$$

---

[4.8]Смотри также определение [1]-3, страницы 294, 295.

[4.9]Смотри также анализ определения [1]-3, страница 295.

для любых $n, m > N$. Согласно теоремам 4.1.4, 4.1.11,

$$(4.2.3) \quad \left\|\int_X d\mu(x) f_m(x) - \int_X d\mu(x) f_n(x)\right\| = \left\|\int_X d\mu(x)(f_m(x) - f_n(x))\right\| \leq \frac{\epsilon}{\mu(X)} \mu(X) = \epsilon$$

следует из неравенства (4.2.2) для любых $n, m > N$. Согласно определению [4]-2.1.19 и неравенству (4.2.3), последовательность интегралов

$$\int_X d\mu(x) f_n(x)$$

является фундаментальной. Следовательно, существует предел (4.2.1) и утверждение 4.2.2.1 верно.

Пусть $f_{1 \cdot n}$, $f_{2 \cdot n}$ - фундаментальные последовательности простых отображений, равномерно сходящиеся к $f$. Согласно определению [4]-2.6.3, так как последовательность $f_{1 \cdot n}$ сходится равномерно к отображению $f$, то для любого $\epsilon \in R$, $\epsilon > 0$, существует $N_1$ такое, что

$$(4.2.4) \quad \|f_{1 \cdot n}(x) - f(x)\| < \frac{\epsilon}{2\mu(X)}$$

для любых $n > N_1$. Согласно определению [4]-2.6.3, так как последовательность $f_{2 \cdot n}$ сходится равномерно к отображению $f$, то для любого $\epsilon \in R$, $\epsilon > 0$, существует $N_2$ такое, что

$$(4.2.5) \quad \|f_{2 \cdot n}(x) - f(x)\| < \frac{\epsilon}{2\mu(X)}$$

для любых $n > N_2$. Пусть

$$N = \max(N_1, N_2)$$

Из неравенств (4.2.4), (4.2.5) следует, что для заданного $\epsilon \in R$, $\epsilon > 0$, существует, зависящее от $\epsilon$, натуральное число $N$ такое, что

$$(4.2.6) \quad \begin{aligned} \|f_{1 \cdot n}(x) - f_{2 \cdot n}(x)\| &= \|f_{1 \cdot n}(x) - f(x) + f(x) - f_{2 \cdot n}(x)\| \\ &\leq \|f_{1 \cdot n}(x) - f(x)\| + \|f_{2 \cdot n}(x) - f(x)\| \\ &< \frac{\epsilon}{\mu(X)} \end{aligned}$$

для любого $n > N$. Согласно теоремам 4.1.4, 4.1.11,

$$(4.2.7) \quad \left\|\int_X d\mu(x) f_{1 \cdot n}(x) - \int_X d\mu(x) f_{2 \cdot n}(x)\right\| = \left\|\int_X d\mu(x)(f_{1 \cdot n}(x) - f_{2 \cdot n}(x))\right\| \leq \frac{\epsilon}{\mu(X)} \mu(X) = \epsilon$$

следует из неравенства (4.2.6) для любого $n > N$. Согласно теореме [4]-2.1.22 и определению [4]-2.1.17,

$$(4.2.8) \quad \lim_{n \to \infty} \int_X d\mu(x) f_{1 \cdot n}(x) = \lim_{n \to \infty} \int_X d\mu(x) f_{2 \cdot n}(x)$$

следует из неравенства (4.2.7) для любого $n > N$. Из равенства (4.2.8) следует, что утверждение 4.2.2.2 верно.

Пусть $f$ - простое отображение. Для доказательства справедливости утверждения 4.2.2.3 достаточно рассмотреть последовательность, в которой $f_n = f$ для любого $n$. $\square$

Теорема 4.2.3. *Пусть*
$$f : X \to A$$
$$g : X \to A$$
*$\mu$-измеримые отображения с компактным множеством значений. Если существуют интегралы*
$$\int_X d\mu(x) f(x)$$
$$\int_X d\mu(x) g(x)$$
*то существует интеграл*
$$\int_X d\mu(x)(f(x) + g(x))$$
*и*
$$(4.2.9) \qquad \int_X d\mu(x)(f(x) + g(x)) = \int_X d\mu(x) f(x) + \int_X d\mu(x) g(x)$$

Доказательство. Согласно теореме 3.3.2, существует последовательность простых отображений
$$f_n : X \to A$$
сходящаяся равномерно к $f$. Согласно теореме 4.2.2, для любого $\epsilon \in R$, $\epsilon > 0$, существует $N_1$ такое, что
$$(4.2.10) \qquad \left\| \int_X d\mu(x) f(x) - \int_X d\mu(x) f_n(x) \right\| < \frac{\epsilon}{2}$$
для любого $n > N_1$. Согласно теореме 3.3.2, существует последовательность простых отображений
$$g_n : X \to A$$
сходящаяся равномерно к $g$. Согласно теореме 4.2.2, для любого $\epsilon \in R$, $\epsilon > 0$, существует $N_2$ такое, что
$$(4.2.11) \qquad \left\| \int_X d\mu(x) g(x) - \int_X d\mu(x) g_n(x) \right\| < \frac{\epsilon}{2}$$
для любого $n > N_2$. Пусть
$$N = \max(N_1, N_2)$$
Согласно теореме 4.1.4, для любого $n > N$ существует интеграл
$$(4.2.12) \qquad \int_X d\mu(x)(f_n(x) + g_n(x)) = \int_X d\mu(x) f_n(x) + \int_X d\mu(x) g_n(x)$$

Согласно теореме 3.3.4, отображение
$$(4.2.13) \qquad h = f + g$$
является $\mu$-измеримым отображением и
$$(4.2.14) \qquad h(x) = \lim_{n \to \infty} f_n(x) + g_n(x)$$

Если $k, n > N$, то

$$
\begin{aligned}
&\left\| \int_X d\mu(x)(f_n(x) + g_n(x)) - \int_X d\mu(x)(f_k(x) + g_k(x)) \right\| \\
&= \left\| \int_X d\mu(x) f_n(x) + \int_X d\mu(x) g_n(x) - \int_X d\mu(x) f_k(x) - \int_X d\mu(x) g_k(x) \right\| \\
&\leq \left\| \int_X d\mu(x) f_n(x) - \int_X d\mu(x) f(x) \right\| + \left\| \int_X d\mu(x) f_k(x) - \int_X d\mu(x) f(x) \right\| \\
&\quad + \left\| \int_X d\mu(x) g_n(x) - \int_X d\mu(x) g(x) \right\| + \left\| \int_X d\mu(x) g_k(x) - \int_X d\mu(x) g(x) \right\| \\
&< 2\epsilon
\end{aligned}
\tag{4.2.15}
$$

следует из (4.2.10), (4.2.11), (4.2.12). Согласно определению [4]-2.1.19, из неравенства (4.2.15) следует, что последовательность интегралов (4.2.12) является фундаментальной. Следовательно, согласно утверждениям (4.2.13), (4.2.14) и теореме 4.2.2, существует интеграл

$$
\int_X d\mu(x)(f(x) + g(x)) = \lim_{n \to \infty} \int_X d\mu(x)(f_n(x) + g_n(x)) \tag{4.2.16}
$$

Если $n > N$, то

$$
\begin{aligned}
&\left\| \int_X d\mu(x) f(x) + \int_X d\mu(x) g(x) - \int_X d\mu(x)(f_n(x) + g_n(x)) \right\| \\
&= \left\| \int_X d\mu(x) f(x) + \int_X d\mu(x) g(x) - \int_X d\mu(x) f_n(x) - \int_X d\mu(x) g_n(x) \right\| \\
&\leq \left\| \int_X d\mu(x) f_n(x) - \int_X d\mu(x) f(x) \right\| + \left\| \int_X d\mu(x) g_n(x) - \int_X d\mu(x) g(x) \right\| \\
&< \epsilon
\end{aligned}
\tag{4.2.17}
$$

следует из (4.2.10), (4.2.11), (4.2.12). Из неравенства (4.2.17) следует, что последовательность интегралов (4.2.12) является фундаментальной. Согласно определению [4]-2.1.17,

$$
\int_X d\mu(x) f(x) + \int_X d\mu(x) g(x) = \lim_{n \to \infty} \int_X d\mu(x)(f_n(x) + g_n(x)) \tag{4.2.18}
$$

Равенство (4.2.9) следует из равенств (4.2.16), (4.2.18). $\square$

Теорема 4.2.4. *Пусть*

$$f : X \to A$$

*измеримое отображение. Интеграл*

$$\int_X d\mu(x) f(x)$$

*существует тогда и только тогда, когда интеграл*

$$\int_X d\mu(x) \|f(x)\|$$

*существует. Тогда*

$$
\left\| \int_X d\mu(x) f(x) \right\| \leq \int_X d\mu(x) \|f(x)\| \tag{4.2.19}
$$

Доказательство. Теорема является следствием теоремы 4.1.5 и определения 4.2.1. $\square$

Теорема 4.2.5. *Пусть $\omega \in \Omega$ - n-арная операция в абелевой $\Omega$-группе $A$. Пусть*

$$f_i : X \to A \quad i = 1, ..., n$$

*$\mu$-измеримое отображение с компактным множеством значений. Если отображение $f_i$, $i = 1, ..., n$, интегрируемо, то отображение*

$$h = f_1 ... f_n \omega$$

*интегрируемо и*

$$\left\| \int_X d\mu(x) h(x) \right\| \leq \int_X d\mu(x)(\|\omega\| \|f_1(x)\| ... \|f_n(x)\|) \quad (4.2.20)$$

Доказательство. Так как множество значений отображения $f_i$ компактно, то, согласно теореме [4]-2.3.11, определена следующая величина

$$F_i = \sup \|f_i(x)\| \quad (4.2.21)$$

Из равенства (4.2.21) и утверждения [4]-2.1.9.1 следует, что

$$F_i \geq 0 \quad (4.2.22)$$

Согласно теореме 3.3.2, для $i = 1, ..., n$, существует последовательность простых отображений

$$f_{i \cdot m} : X \to A$$

сходящаяся равномерно к $f_i$. Из равенства

$$f_i(x) = \lim_{m \to \infty} f_{i \cdot m}(x)$$

и теоремы [4]-2.6.5 следует, что для заданного

$$\delta_1 \in R, \; \delta_1 > 0 \quad (4.2.23)$$

существует $M_i$ такое, что из условия $m > M_i$ следует

$$\|f_{i \cdot m}(x) - f_i(x)\| < \delta_1 \quad (4.2.24)$$

Пусть

$$M = \max(M_1, ..., M_n)$$

Для $k, m > M$, неравенство

$$\|f_{i \cdot k}(x) - f_{i \cdot m}(x)\| = \|f_{i \cdot k}(x) - f_i(x) + f_i(x) - f_{i \cdot m}(x)\| \\ \leq \|f_{i \cdot k}(x) - f_i(x)\| + \|f_i(x) - f_{i \cdot m}(x)\| < 2\delta_1 \quad (4.2.25)$$

следует из неравенства (4.2.24) и утверждения [4]-2.5.14.3. Для $m > M$, неравенство

$$\|f_{i \cdot m}(x)\| < F_i + \delta_1 \quad (4.2.26)$$

следует из [4]-(2.1.1), (4.2.21), (4.2.24) и неравенства

$$\|f_{i \cdot m}(x)\| - \|f_i(x)\| \leq \|f_{i \cdot m}(x) - f_i(x)\| < \delta_1$$

Согласно теореме 4.1.6, для любого $m > M$ существует интеграл

$$\int_X d\mu(x)(f_{1 \cdot m}(x) ... f_{n \cdot m}(x)\omega) \quad (4.2.27)$$

Согласно утверждению [4]-2.1.9.3,

$$\left\| \int_X d\mu(x)(f_{1\cdot k}(x)...f_{n\cdot k}(x)\omega) - \int_X d\mu(x)(f_{1\cdot m}(x)...f_{n\cdot m}(x)\omega) \right\|$$
$$= \left\| \int_X d\mu(x)(f_{1\cdot k}(x)f_{2\cdot k}(x)...f_{n\cdot k}(x)\omega) \right.$$
$$- \int_X d\mu(x)(f_{1\cdot m}(x)f_{2\cdot k}(x)...f_{n\cdot k}(x)\omega)$$
$$+ \int_X d\mu(x)(f_{1\cdot m}(x)f_{2\cdot k}(x)...f_{n\cdot k}(x)\omega)$$
$$(4.2.28) \qquad - ... - \left. \int_X d\mu(x)(f_{1\cdot m}(x)...f_{n\cdot m}(x)\omega) \right\|$$
$$\leq \left\| \int_X d\mu(x)(f_{1\cdot k}(x)f_{2\cdot k}(x)...f_{n\cdot k}(x)\omega) \right.$$
$$- \left. \int_X d\mu(x)(f_{1\cdot m}(x)f_{2\cdot k}(x)...f_{n\cdot k}(x)\omega) \right\|$$
$$+ ... + \left\| \int_X d\mu(x)(f_{1\cdot m}(x)...f_{n-1\cdot m}(x)f_{n\cdot k}(x)\omega) \right.$$
$$- \left. \int_X d\mu(x)(f_{1\cdot m}(x)...f_{n-1\cdot m}(x)f_{n\cdot m}(x)\omega) \right\|$$

Согласно определению [4]-2.1.3,

$$\left\| \int_X d\mu(x)(f_{1\cdot k}(x)...f_{n\cdot k}(x)\omega) - \int_X d\mu(x)(f_{1\cdot m}(x)...f_{n\cdot m}(x)\omega) \right\|$$
$$\leq \left\| \int_X d\mu(x)(f_{1\cdot k}(x)f_{2\cdot k}(x)...f_{n\cdot k}(x)\omega - f_{1\cdot m}(x)f_{2\cdot k}(x)...f_{n\cdot k}(x)\omega) \right\| + ...$$
$$+ \left\| \int_X d\mu(x)(f_{1\cdot m}(x)...f_{n-1\cdot m}(x)f_{n\cdot k}(x)\omega \right.$$
$$(4.2.29) \qquad - \left. f_{1\cdot m}(x)...f_{n-1\cdot m}(x)f_{n\cdot m}(x)\omega) \right\|$$
$$= \left\| \int_X d\mu(x)((f_{1\cdot k}(x) - f_{1\cdot m}(x))f_{2\cdot k}(x)...f_{n\cdot k}(x)\omega) \right\| + ...$$
$$+ \left\| \int_X d\mu(x)(f_{1\cdot m}(x)...f_{n-1\cdot m}(x)(f_{n\cdot k}(x) - f_{n\cdot m}(x))\omega) \right\|$$

следует из (4.2.9), (4.2.28). Согласно теореме 4.1.6, для любого $m > M$ неравенство

$$\left\| \int_X d\mu(x)(f_{1\cdot k}(x)...f_{n\cdot k}(x)\omega) - \int_X d\mu(x)(f_{1\cdot m}(x)...f_{n\cdot m}(x)\omega) \right\|$$
$$(4.2.30) \quad \leq \int_X d\mu(x)(\|\omega\| \|f_{1\cdot k}(x) - f_{1\cdot m}(x)\| \|f_{2\cdot k}(x)\|...\|f_{n\cdot k}(x)\|) + ...$$
$$+ \int_X d\mu(x)(\|\omega\| \|f_{1\cdot m}(x)\|...\|f_{n-1\cdot m}(x)\| \|f_{n\cdot k}(x) - f_{n\cdot m}(x)\|)$$
$$\leq 2\mu(x)\|\omega\|\delta_1((F_2 + \delta_1)...(F_n + \delta_1) + ... + (F_1 + \delta_1)...(F_{n-1} + \delta_1))$$

следует из неравенств (4.2.25), (4.2.26), (4.2.29). Пусть

$$(4.2.31) \quad \epsilon_1 = 2\mu(x)\|\omega\|\delta_1((F_2+\delta_1)...(F_n+\delta_1) + ... + (F_1+\delta_1)...(F_{n-1}+\delta_1))$$

Из утверждений (4.2.22), (4.2.23) следует, что

$$(4.2.32) \qquad\qquad \epsilon_1 > 0 \quad \frac{d\epsilon_1}{d\delta_1} > 0$$

Из равенства (4.2.31) и утверждения (4.2.32) следует, что $\epsilon_1$ является полиномиальной строго монотонно возрастающей функцией $\delta_1 > 0$. При этом

$$\delta_1 = 0 \Rightarrow \epsilon_1 = 0$$

Согласно теореме [4]-2.3.5, отображение (4.2.31) отображает интервал $[0,\delta_1)$ в интервал $[0,\epsilon_1)$. Согласно теореме [4]-2.3.3, для заданного $\epsilon > 0$ существует такое $\delta > 0$, что

$$\epsilon_1(\delta) < \epsilon$$

Согласно построению, значение $M$ зависит от значения $\delta_1$. Мы выберем значение $M$, соответствующее $\delta_1 = \delta$. Следовательно, для заданного $\epsilon \in R$, $\epsilon > 0$, существует $M$ такое, что из условия $m > M$ следует

$$(4.2.33) \quad \left\|\int_X d\mu(x)(f_{1\cdot k}(x)...f_{n\cdot k}(x)\omega) - \int_X d\mu(x)(f_{1\cdot m}(x)...f_{n\cdot m}(x)\omega)\right\| < \epsilon$$

Из равенства (4.2.33) следует, что последовательность интегралов (4.2.27) является фундаментальной последовательностью. Следовательно, последовательность интегралов (4.2.27) имеет предел

$$(4.2.34) \qquad \int_X d\mu(x)(f_1(x)...f_n(x)\omega) = \lim_{m\to\infty} \int_X d\mu(x)(f_{1\cdot m}(x)...f_{n\cdot m}(x)\omega)$$

Согласно теореме 4.1.6, для любого $m > M$

$$(4.2.35) \quad \left\|\int_X d\mu(x)(f_{1\cdot m}(x)...f_{n\cdot m}(x)\omega)\right\| \le \int_X d\mu(x)(\|\omega\|\|f_{1\cdot m}(x)\|...\|f_{n\cdot m}(x)\|)$$

Неравенство (4.2.20) следует из неравенства (4.2.35) при предельном переходе $m \to \infty$. $\square$

Теорема 4.2.6. *Пусть на множестве $X$ определена $\sigma$-аддитивная мера $\mu$. Пусть*

$$f : A_1 \dashrightarrow A_2$$

*представление $\Omega_1$-группы $A_1$ с нормой $\|x\|_1$ в $\Omega_2$-группе $A_2$ с нормой $\|x\|_2$. Пусть*

$$g_i : X \to A_i \quad i = 1,2$$

*интегрируемое отображение с компактным множеством значений. Тогда отображение*

$$h = f_X(g_1)(g_2)$$

*интегрируемо и*

$$(4.2.36) \qquad \left\|\int_X d\mu(x)h(x)\right\|_2 \le \int_X d\mu(x)(\|f\|\|g_1(x)\|_1\|g_2(x)\|_2)$$

Доказательство. Так как множество значений отображения $g_i$ компактно, то, согласно теореме [4]-2.3.11, определена следующая величина

$$(4.2.37) \qquad G_i = \sup \|g_i(x)\|_i$$

Из равенства (4.2.37) и утверждения [4]-2.1.9.1 следует, что

$$(4.2.38) \qquad G_i \geq 0$$

Согласно теореме 3.3.2, для $i = 1, 2$, существует последовательность простых отображений

$$g_{i \cdot n} : X \to A_i$$

сходящаяся равномерно к $g_i$. Из равенства

$$g_i(x) = \lim_{n \to \infty} g_{i \cdot n}(x)$$

и теоремы [4]-2.6.5 следует, что для заданного

$$(4.2.39) \qquad \delta_1 \in R, \ \delta_1 > 0$$

существует $N_i$ такое, что из условия $n > N_i$ следует

$$(4.2.40) \qquad \|g_{i \cdot n}(x) - g_i(x)\|_i < \delta_1$$

Пусть

$$N = \max(N_1, N_2)$$

Для $k, n > N$, неравенство

$$(4.2.41) \qquad \begin{aligned} \|g_{i \cdot k}(x) - g_{i \cdot n}(x)\|_i &= \|g_{i \cdot k}(x) - g_i(x) + g_i(x) - g_{i \cdot n}(x)\|_i \\ &\leq \|g_{i \cdot k}(x) - g_i(x)\|_i + \|g_i(x) - g_{i \cdot n}(x)\|_i < 2\delta_1 \end{aligned}$$

следует из неравенства (4.2.40) и утверждения [4]-2.5.14.3. Для $n > N$, неравенство

$$(4.2.42) \qquad \|g_{i \cdot n}(x)\|_i < G_i + \delta_1$$

следует из [4]-(2.1.1), (4.2.37), (4.2.40) и неравенства

$$\|g_{i \cdot n}(x)\|_i - \|g_i(x)\|_i \leq \|g_{i \cdot n}(x) - g_i(x)\|_i < \delta_1$$

Согласно теореме 4.1.7, для любого $n > N$ существует интеграл

$$(4.2.43) \qquad \int_X d\mu(x)(f(g_{1 \cdot n}(x))(g_{2 \cdot n}(x)))$$

Согласно утверждению [4]-2.1.9.3,

$$(4.2.44) \qquad \begin{aligned} & \left\| \int_X d\mu(x)(f(g_{1 \cdot k}(x))(g_{2 \cdot k}(x))) - \int_X d\mu(x)(f(g_{1 \cdot n}(x))(g_{2 \cdot n}(x))) \right\|_2 \\ &= \left\| \int_X d\mu(x)(f(g_{1 \cdot k}(x))(g_{2 \cdot k}(x))) - \int_X d\mu(x)(f(g_{1 \cdot n}(x))(g_{2 \cdot k}(x))) \right. \\ &\quad \left. + \int_X d\mu(x)(f(g_{1 \cdot n}(x))(g_{2 \cdot k}(x))) - \int_X d\mu(x)(f(g_{1 \cdot n}(x))(g_{2 \cdot n}(x))) \right\|_2 \\ &\leq \left\| \int_X d\mu(x)(f(g_{1 \cdot k}(x))(g_{2 \cdot k}(x))) - \int_X d\mu(x)(f(g_{1 \cdot n}(x))(g_{2 \cdot k}(x))) \right\|_2 \\ &\quad + \left\| \int_X d\mu(x)(f(g_{1 \cdot n}(x))(g_{2 \cdot k}(x))) - \int_X d\mu(x)(f(g_{1 \cdot n}(x))(g_{2 \cdot n}(x))) \right\|_2 \end{aligned}$$

Согласно определению [2]-2.1.4, отображение $f$ является гомоморфизмом абелевой группы $A_1$ и для любого $a_1 \in A_1$ отображение $f(a_1)$ является гомоморфизмом абелевой группы $A_2$. Следовательно,

$$\left\| \int_X d\mu(x)(f(g_{1\cdot k}(x))(g_{2\cdot k}(x))) - \int_X d\mu(x)(f(g_{1\cdot n}(x))(g_{2\cdot n}(x))) \right\|_2$$
$$\leq \left\| \int_X d\mu(x)(f(g_{1\cdot k}(x))(g_{2\cdot k}(x)) - f(g_{1\cdot n}(x))(g_{2\cdot k}(x))) \right\|_2$$
(4.2.45)
$$+ \left\| \int_X d\mu(x)(f(g_{1\cdot n}(x))(g_{2\cdot k}(x))) - f(g_{1\cdot n}(x))(g_{2\cdot n}(x))) \right\|_2$$
$$= \left\| \int_X d\mu(x)(f(g_{1\cdot k}(x) - g_{1\cdot n}(x))g_{2\cdot k}(x))) \right\|_2$$
$$+ \left\| \int_X d\mu(x)(f(g_{1\cdot n}(x))(g_{2\cdot k}(x) - g_{2\cdot n}(x))) \right\|_2$$

следует из (4.2.9), (4.2.44). Согласно теореме 4.1.7, для любого $n > N$ неравенство

$$\left\| \int_X d\mu(x)(f(g_{1\cdot k}(x))(g_{2\cdot k}(x))) - \int_X d\mu(x)(f(g_{1\cdot n}(x))(g_{2\cdot n}(x))) \right\|_2$$
(4.2.46)
$$\leq \int_X d\mu(x)(\|f\| \|g_{1\cdot k}(x) - g_{1\cdot n}(x)\|_1 \|g_{2\cdot k}(x)\|_2)$$
$$+ \int_X d\mu(x)(\|f\| \|g_{1\cdot n}(x)\|_1 \|g_{2\cdot k}(x) - f_{2\cdot n}(x)\|\|_2)$$
$$\leq 2\mu(x)\|f\|\delta_1((G_2 + \delta_1) + (G_1 + \delta_1))$$

следует из неравенств (4.2.41), (4.2.42), (4.2.45). Пусть
(4.2.47) $\quad\quad\quad\quad \epsilon_1 = 2\mu(x)\|f\|\delta_1((G_2 + \delta_1) + (G_1 + \delta_1))$

Из утверждений (4.2.38), (4.2.39) следует, что
(4.2.48) $\quad\quad\quad\quad \epsilon_1 > 0 \quad \dfrac{d\epsilon_1}{d\delta_1} > 0$

Из равенства (4.2.47) и утверждения (4.2.48) следует, что $\epsilon_1$ является полиномиальной строго монотонно возрастающей функцией $\delta_1 > 0$. При этом

$$\delta_1 = 0 \Rightarrow \epsilon_1 = 0$$

Согласно теореме [4]-2.3.5, отображение (4.2.47) отображает интервал $[0, \delta_1)$ в интервал $[0, \epsilon_1)$. Согласно теореме [4]-2.3.3, для заданного $\epsilon > 0$ существует такое $\delta > 0$, что

$$\epsilon_1(\delta) < \epsilon$$

Согласно построению, значение $N$ зависит от значения $\delta_1$. Мы выберем значение $N$, соответствующее $\delta_1 = \delta$. Следовательно, для заданного $\epsilon \in R$, $\epsilon > 0$, существует $N$ такое, что из условия $n > N$ следует

(4.2.49) $\quad \left\| \int_X d\mu(x)(f(g_{1\cdot k}(x))(g_{2\cdot k}(x))) - \int_X d\mu(x)(f(g_{1\cdot n}(x))(g_{2\cdot n}(x))) \right\|_2 < \epsilon$

Из равенства (4.2.49) следует, что последовательность интегралов (4.2.43) является фундаментальной последовательностью. Следовательно, последовательность интегралов (4.2.43) имеет предел

$$(4.2.50) \qquad \int_X d\mu(x)(f(g_1(x))(g_2(x))) = \lim_{n\to\infty} \int_X d\mu(x)(f(g_{1\cdot n}(x))(g_{2\cdot n}(x)))$$

Согласно теореме 4.1.7, для любого $n > N$

$$(4.2.51) \qquad \left\| \int_X d\mu(x)(f(g_{1\cdot n}(x))(g_{2\cdot n}(x))) \right\|_2 \le \int_X d\mu(x)(\|f\| \|g_{1\cdot n}(x)\|_1 \|g_{2\cdot n}(x)\|_2)$$

Неравенство (4.2.36) следует из неравенства (4.2.51) при предельном переходе $n \to \infty$. $\square$

Теорема 4.2.7. *Пусть на множестве $X$ определена $\sigma$-аддитивная мера $\mu$. Пусть*

$$f : A_1 \relbar\joinrel\twoheadrightarrow A_2$$

*эффективное представление $\Omega_1$-группы $A_1$ с нормой $\|x\|_1$ в $\Omega_2$-группе $A_2$ с нормой $\|x\|_2$. Пусть преобразования представления $f$ являются автоморфизмами представления*

$$R \relbar\joinrel\twoheadrightarrow A_2$$

*Пусть*

$$g_2 : X \to A_2$$

*интегрируемое отображение в нормированную $\Omega$-группу $A_2$. Тогда отображение*

$$h = a_1 g_2$$

*интегрируемо и*

$$\int_X d\mu(x)(a_1 g_2(x)) = a_1 \int_X d\mu(x) g_2(x)$$

Доказательство. Теорема является следствием теоремы 4.1.8 и определения 4.2.1. $\square$

Теорема 4.2.8. *Пусть на множестве $X$ определена $\sigma$-аддитивная мера $\mu$. Пусть*

$$f : A_1 \relbar\joinrel\twoheadrightarrow A_2$$

*эффективное представление $\Omega_1$-группы $A_1$ с нормой $\|x\|_1$ в $\Omega_2$-группе $A_2$ с нормой $\|x\|_2$. Пусть преобразования представления*

$$R \relbar\joinrel\twoheadrightarrow A_1$$

*являются автоморфизмами представления $f$. Пусть*

$$g_1 : X \to A_1$$

*интегрируемое отображение в нормированную $\Omega$-группу $A_1$. Тогда отображение*

$$h = g_1 a_2$$

*интегрируемо и*

$$\int_X d\mu(x)(g_1(x) a_2) = \left( \int_X d\mu(x) g_1(x) \right) a_2$$

Доказательство. Теорема является следствием теоремы 4.1.9 и определения 4.2.1. □

Теорема 4.2.9. *Пусть $\mu$-измеримое отображение*
$$f : X \to A$$
*удовлетворяет условию*
$$\|f(x)\| \leq M$$
*Если мера множества $X$ конечна, то*
$$\text{(4.2.52)} \qquad \int_X d\mu(x)\|f(x)\| \leq M\mu(X)$$

Доказательство. Теорема является следствием теоремы 4.1.10 и определения 4.2.1. □

Теорема 4.2.10. *Пусть $\mu$-измеримое отображение*
$$f : X \to A$$
*удовлетворяет условию*
$$\|f(x)\| \leq M$$
*Если мера множества $X$ конечна, то отображение $f$ интегрируемо и*
$$\text{(4.2.53)} \qquad \left\|\int_X d\mu(x)f(x)\right\| \leq M\mu(X)$$

Доказательство. Интегрируемость отображения $f$ следует из теорем 4.2.4, 4.2.9. Неравенство (4.2.53) следует из неравенств (4.2.19), (4.2.52). □

## 4.3. Интеграл Лебега как отображение множества

Рассмотрим $\mu$-измеримое отображение
$$f : X \to A$$
в Ω-группу $A$. Пусть $\mathcal{C}_X$ - $\sigma$-алгебра измеримых множеств множества $X$. Если $Y \in \mathcal{C}_X$, то мы можем рассматривать выражение
$$\text{(4.3.1)} \qquad F(Y) = \int_Y d\mu(x)f(x)$$
как отображение
$$F : \mathcal{C}_X \to A$$

Лемма 4.3.1. *Пусть*
$$f : X \to A$$
*простое отображение. Пусть*
$$\text{(4.3.2)} \qquad X = \bigcup_i X_i \quad i \neq j => X_i \cap X_j = \emptyset$$
*конечное или счётное объединение множеств $X_i$. Отображение $f$ интегрируемо по множеству $X$ тогда и только тогда, когда отображение $f$ интегрируемо по каждому $X_i$*
$$\text{(4.3.3)} \qquad \int_X d\mu(x)f(x) = \sum_i \int_{X_i} d\mu(x)f(x)$$
*где ряд в правой части сходится по норме.*

Доказательство. Пусть множество $\{f_1, f_2, ...\}$ является областью определения отображения $f$. Пусть

(4.3.4) $$F_n = \{x \in X : f(x) = f_n\}$$

(4.3.5) $$F_{in} = \{x \in X_i : f(x) = f_n\}$$

Равенства

(4.3.6) $$F_{in} = F_n \cap X_i$$

(4.3.7) $$F_n = \bigcup_i F_{in}$$

следуют из равенств (4.3.2), (4.3.4), (4.3.5). Равенство

(4.3.8) $$F_{in} \cap F_{im} = \emptyset$$

следует из равенств (4.3.2), (4.3.6). Равенство

(4.3.9) $$\mu(F_n) = \sum_n \mu(F_{in})$$

следует из равенств (4.3.7), (4.3.8) и утверждения 2.2.10.2. Равенство

(4.3.10) $$\int_X d\mu(x) f(x) = \sum_n \mu(F_n) f_n = \sum_n \left(\sum_i \mu(F_{in})\right) f_n = \sum_n \sum_i \mu(F_{in}) f_n$$

следует из равенств (4.1.1), (4.3.9). Неравенство

(4.3.11) $$\|\mu(F_{in}) f_n\| \leq \|\mu(F_n) f_n\| \leq \sum_i \|\mu(F_{in}) f_n\|$$

следует из равенств (4.3.9), (4.3.6) и утверждения [4]-2.1.9.3. Неравенство

(4.3.12) $$\sum_n \|\mu(F_{in}) f_i\| \leq \sum_n \|\mu(F_n) f_n\| \leq \sum_n \sum_i \|\mu(F_{in}) f_n\|$$

следует из неравенства (4.3.11). Сходимость по норме ряда $\sum_n \sum_i \|\mu(F_{in}) f_n\|$ означает

(4.3.13) $$\sum_n \sum_i \|\mu(F_{in}) f_n\| = \sum_i \sum_n \|\mu(F_{in}) f_n\|$$

(4.3.14) $$\sum_n \sum_i \mu(F_{in}) f_n = \sum_i \sum_n \mu(F_{in}) f_n$$

Из (4.3.12), (4.3.13) следует, что ряд $\sum_n \mu(F_n) f_n$ сходится нормально тогда и только тогда, когда для любого $i$ ряд $\sum_n \mu(F_{in}) f_n$ сходится нормально. Согласно определению 4.1.2 для каждого $i$, определён интеграл

(4.3.15) $$\int_{X_i} d\mu(x) f(x) = \sum_n \mu(F_{in}) f_n$$

Следовтельно, отображение $f$ интегрируемо по множеству $X$ тогда и только тогда, когда отображение $f$ интегрируемо по каждому $X_i$. Равенство (4.3.3) следует из равенств (4.3.10), (4.3.14), (4.3.15). $\square$

Теорема 4.3.2 ($\sigma$-аддитивность интеграла Лебега). *Пусть*
$$f: X \to A$$
*$\mu$-измеримое отображение. Пусть*

(4.3.16) $$X = \bigcup_i X_i \quad i \neq j => X_i \cap X_j = \emptyset$$

*конечное или счётное объединение множеств $X_i$. Отображение $f$ интегрируемо по множеству $X$ тогда и только тогда, когда отображение $f$ интегрируемо по каждому $X_i$*

(4.3.17) $$\int_X d\mu(x) f(x) = \sum_i \int_{X_i} d\mu(x) f(x)$$

*где ряд в правой части сходится по норме.*

Доказательство. Согласно определению 4.2.1, для любого $\epsilon \in R$, $\epsilon > 0$, существует простое отображение
$$g: X \to A$$
интегрируемое на $X$ и такое, что для любого $x \in X$

(4.3.18) $$\|f(x) - g(x)\| < \frac{\epsilon}{2\mu(X)}$$

Согласно лемме 4.3.1,

(4.3.19) $$\int_X d\mu(x) g(x) = \sum_i \int_{X_i} d\mu(x) g(x)$$

где $g$ интегрируема на каждом $X_i$ и ряд (4.3.19) сходится. Согласно теореме 4.2.2, отображение $f$ интегрируема на каждом $X_i$. Согласно теореме 4.2.10,

(4.3.20) $$\left\| \int_X d\mu(x) f(x) - \int_X d\mu(x) g(x) \right\| = \left\| \int_X d\mu(x)(f(x) - g(x)) \right\|$$
$$< \frac{\epsilon}{2\mu(X)} \mu(X) = \frac{\epsilon}{2}$$

(4.3.21) $$\left\| \int_{X_i} d\mu(x) f(x) - \int_{X_i} d\mu(x) g(x) \right\| = \left\| \int_{X_i} d\mu(x)(f(x) - g(x)) \right\|$$
$$< \frac{\epsilon}{2\mu(X)} \mu(X_i)$$

следуют из (4.3.18). Согласно утверждениям [4]-2.1.9.3, 2.2.10.2,

$$(4.3.22)\quad \begin{aligned}&\left\|\sum_i \int_{X_i} d\mu(x)f(x) - \int_X d\mu(x)g(x)\right\|\\&=\left\|\sum_i \int_{X_i} d\mu(x)f(x) - \sum_i \int_{X_i} d\mu(x)g(x)\right\|\\&=\left\|\sum_i \left(\int_{X_i} d\mu(x)f(x) - \int_{X_i} d\mu(x)g(x)\right)\right\|\\&<\sum_i \left\|\int_{X_i} d\mu(x)(f(x)-g(x))\right\|\\&<\frac{\epsilon}{2\mu(X)}\sum_i \mu(X_i) = \frac{\epsilon}{2}\end{aligned}$$

следует из (4.2.9), (4.3.21). Неравенство

$$(4.3.23)\quad \begin{aligned}0\leq &\left\|\int_X d\mu(x)f(x) - \sum_i \int_{X_i} d\mu(x)f(x)\right\|\\=&\left\|\int_X d\mu(x)f(x) - \int_X d\mu(x)g(x)\right.\\&\left.+\int_X d\mu(x)g(x) - \sum_i \int_{X_i} d\mu(x)f(x)\right\|\\\leq &\left\|\int_X d\mu(x)f(x) - \int_X d\mu(x)g(x)\right\|\\&+\left\|\int_X d\mu(x)g(x) - \sum_i \int_{X_i} d\mu(x)f(x)\right\|\\<&\epsilon\end{aligned}$$

следует из (4.3.20), (4.3.22) и утверждения [4]-2.1.9.2. Так как $\epsilon$ произвольно, то равенство (4.3.17) следует из неравенства (4.3.23). $\square$

Следствие 4.3.3. *Если $\mu$-измеримое отображение*

$$f: X \to A$$

*интегрируемо на множестве $X$, то отображение $f$ интегрируемо на $\mu$-измеримом множестве $X' \subset X$.* $\square$

Теорема 4.3.4. *Пусть $\mu$-измеримые отображения*

$$f: X \to A$$
$$g: X \to A$$

*удовлетворяют условию*

$$(4.3.24)\qquad \|f(x)\| \leq M$$

$$(4.3.25)\qquad \|g(x)\| \leq M$$

*Если $f(x) = g(x)$ почти всюду, то*

$$\int_X d\mu(x) f(x) = \int_X d\mu(x) g(x) \tag{4.3.26}$$

Доказательство. Из равенств (4.3.24), (4.3.25) и теоремы 4.2.10 следует, что интегралы

$$\int_X d\mu(x) f(x)$$
$$\int_X d\mu(x) g(x)$$

существуют. Пусть

$$X_1 = \{x \in X : f(x) = g(x)\} \tag{4.3.27}$$
$$X_2 = \{x \in X : f(x) \neq g(x)\}$$

Так как $X = X_1 \cup X_2$, $X_1 \cap X_2 = \emptyset$, то, согласно теореме 4.3.2,

$$\int_X d\mu(x) f(x) = \int_{X_1} d\mu(x) f(x) + \int_{X_2} d\mu(x) f(x) \tag{4.3.28}$$

$$\int_X d\mu(x) g(x) = \int_{X_1} d\mu(x) g(x) + \int_{X_2} d\mu(x) g(x) \tag{4.3.29}$$

Согласно условию теоремы

$$\mu(X_2) = 0 \tag{4.3.30}$$

Из теоремы 4.2.10, утверждения [4]-2.1.9.1 и равенства (4.3.30) следует, что

$$0 \leq \left\| \int_{X_2} d\mu(x) f(x) \right\| \leq M\mu(X_2) = 0 \tag{4.3.31}$$

$$0 \leq \left\| \int_{X_2} d\mu(x) g(x) \right\| \leq M\mu(X_2) = 0 \tag{4.3.32}$$

Равенство

$$\int_{X_2} d\mu(x) f(x) = \int_{X_2} d\mu(x) g(x) = 0 \tag{4.3.33}$$

следует из (4.3.31), (4.3.32) и утверждения [4]-2.1.9.2. Равенство

$$\int_{X_1} d\mu(x) f(x) = \int_{X_1} d\mu(x) g(x) \tag{4.3.34}$$

следует из (4.3.27). Равенство (4.3.26) следует из (4.3.28), (4.3.29), (4.3.33), (4.3.34). . □

Теорема 4.3.5 (Неравенство Чебышева). *Если $c > 0$, то*[4.10]

$$\mu\{x \in X : \|f(x)\| \geq c\} \leq \frac{1}{c} \int_X d\mu(x) \|f(x)\| \tag{4.3.35}$$

---

[4.10]Смотри также [1], теорема на странице 300.

Доказательство. Пусть

(4.3.36) $$X' = \{x \in X : \|f(x)\| \geq c\}$$

Тогда

(4.3.37) $$\int_X d\mu(x)\|f(x)\| = \int_{X'} d\mu(x)\|f(x)\| + \int_{X-X'} d\mu(x)\|f(x)\|$$
$$\geq \int_{X'} d\mu(x)\|f(x)\| \geq c\mu(X')$$

(4.3.35) следует из (4.3.36), (4.3.37). $\square$

Теорема 4.3.6. *Если*

(4.3.38) $$\int_X d\mu(x)\|f(x)\| = 0$$

*то $f(x) = 0$ почти всюду.*

Доказательство. Из (4.3.35), (4.3.38) следует, что

(4.3.39) $$\mu\left\{x \in X : \|f(x)\| \geq \frac{1}{n}\right\} \leq n \int_X d\mu(x)\|f(x)\| = 0$$

для $n = 1, 2, \ldots$ . Равенство

$$\mu\{x \in X : f(x) \geq 0\} = \lim_{n \to \infty} \left\{x \in X : \|f(x)\| \geq \frac{1}{n}\right\} = 0$$

следует из (4.3.39). $\square$

Теорема 4.3.7 (Нормальная непрерывность интеграла Лебега). *Если*[4.11]

$$f : X \to A$$

*отображение, интегрируемое на множестве $X$, то для любого $\epsilon \in R$, $\epsilon > 0$, существует $\delta > 0$ такое, что*

(4.3.40) $$\left\|\int_E d\mu(x) f(x)\right\| < \epsilon$$

*для всякого $E \in \mathcal{C}_X$ такого, что $\mu(E) < \delta$.*

Доказательство. Согласно теореме 4.2.4, интеграл

$$\int_X d\mu(x) f(x)$$

существует тогда и только тогда, когда интеграл

$$\int_X d\mu(x)\|f(x)\|$$

существует. Пусть

(4.3.41) $$X_i = \{x \in X : i \leq \|f(x)\| < i+1\}$$

(4.3.42) $$Y_N = \bigcup_{i=1}^{N} X_i$$
$$Z_N = X \setminus Y_N$$

---

[4.11]Смотри также теорему [1]-5, страницы 301, 302.

Согласно теореме 4.3.2,

$$(4.3.43) \qquad \int_X d\mu(x)\|f(x)\| = \sum_{i=0}^{\infty} \int_{X_i} d\mu(x)\|f(x)\|$$

Так как ряд в (4.3.43) сходится по норме, то существует $N$ такое, что

$$(4.3.44) \qquad \sum_{i=N+1}^{\infty} \int_{X_i} d\mu(x)\|f(x)\| = \int_{Z_N} d\mu(x)\|f(x)\| < \frac{\epsilon}{2}$$

Пусть

$$(4.3.45) \qquad 0 < \delta < \frac{\epsilon}{2(N+1)}$$

Пусть $\mu(E) < \delta$. Тогда

$$(4.3.46) \qquad \begin{aligned} \left\|\int_E d\mu(x)f(x)\right\| &\leq \int_E d\mu(x)\|f(x)\| \\ &= \int_{E \cap Y_N} d\mu(x)\|f(x)\| + \int_{E \cap Z_N} d\mu(x)\|f(x)\| \end{aligned}$$

Согласно теореме 4.2.9, неравенство

$$(4.3.47) \qquad \int_{E \cap Y_N} d\mu(x)\|f(x)\| \leq (N+1)\mu(E) < (N+1)\delta = \frac{\epsilon}{2}$$

следует из (4.3.41), (4.3.42). Неравенство

$$(4.3.48) \qquad \int_{E \cap Z_N} d\mu(x)\|f(x)\| \leq \int_{Z_N} d\mu(x)\|f(x)\| < \frac{\epsilon}{2}$$

следует из (4.3.46). Неравенство (4.3.40) следует из неравенств (4.3.46), (4.3.47), (4.3.48). $\square$

## 4.4. Переход к пределу под знаком интеграла Лебега

Теорема 4.4.1. *Пусть отображение*

$$g : X \to R$$

*интегрируемо и почти всюду*

$$\|f(x)\| \leq g(x)$$

*Тогда отображение $f$ интегрируемо и*

$$(4.4.1) \qquad \left\|\int_X d\mu(x)f(x)\right\| \leq \int_X d\mu(x)g(x)$$

Доказательство. Согласно теореме 4.3.4,

$$\mu(X) = 0 \implies \int_X d\mu(x)f(x) = 0$$

Следовательно, мы можем предположить, что условие теоремы верно для любого $x \in X$. Теорема следует из теорем 4.2.4, [1]-$VII$ на странице 297. $\square$

Теорема 4.4.2 (теорема Лебега о мажорируемой сходимости). *Пусть*[4.12]

$$f_n : X \to A$$

*последовательность $\mu$-измеримых отображений в $\Omega$-группу $A$ и*

$$f(x) = \lim_{n\to\infty} f_n(x)$$

*почти всюду. Пусть отображение*

$$g : X \to R$$

*интегрируемо и для всех $n \in N$ почти всюду*

(4.4.2) $$\|f_n(x)\| \le g(x)$$

*Тогда отображение $f$ интегрируемо и*

(4.4.3) $$\int_X d\mu(x) f(x) = \lim_{n\to\infty} \int_X d\mu(x) f_n(x)$$

Доказательство. Согласно теореме 4.3.4,

$$\mu(X) = 0 \implies \int_X d\mu(x) f(x) = 0$$

Следовательно, мы можем предположить, что условие теоремы верно для любого $x \in X$. Согласно теореме [5]-2 на странице 56, неравенство

(4.4.4) $$\|f(x)\| \le g(x)$$

следует из неравенства (4.4.2). Согласно теореме 4.4.1, отображение $f$ интегрируемо. Согласно теореме 4.3.7, для любого $\epsilon \in R, \epsilon > 0$, существует $\delta > 0$ такое, что

(4.4.5) $$\int_Y d\mu(x) g(x) < \frac{\epsilon}{4}$$

для всякого $Y \in \mathcal{C}_X$ такого, что

(4.4.6) $$\mu(Y) < \delta$$

Согласно теореме 3.4.5, множество $Y$, удовлетворяющее условию (4.4.6), может быть выбрано таким образом, что последовательность $f_n$ сходится равномерно на множестве $Z = X \setminus Y$. Согласно определению [4]-2.6.3, найдётся $N$ такое, что при $n > N$, $x \in Z$ выполнено

(4.4.7) $$\|f(x) - f_n(x)\| < \frac{\epsilon}{2\mu(C)}$$

Неравенство

(4.4.8) $$\begin{aligned}&\left\|\int_X d\mu(x) f(x) - \int_X d\mu(x) f_n(x)\right\| \\ &\le \left\|\int_Z d\mu(x)(f(x) - f_n(x))\right\| + \left\|\int_Y d\mu(x) f(x)\right\| - \left\|\int_Y d\mu(x) f_n(x)\right\| \\ &< \frac{\epsilon}{2} + \frac{\epsilon}{4} + \frac{\epsilon}{4} = \epsilon\end{aligned}$$

---

[4.12]Смотри также теорему [1]-6 на страницах 302, 303.

следует из равенства

$$\int_X d\mu(x)f(x) - \int_X d\mu(x)f_n(x)$$
$$= \int_Y d\mu(x)f(x) + \int_Z d\mu(x)f(x) - \int_Y d\mu(x)f_n(x) - \int_Z d\mu(x)f_n(x)$$
$$= \int_Z d\mu(x)(f(x) - f_n(x)) + \int_Y d\mu(x)f(x) - \int_Y d\mu(x)f_n(x)$$

неравенств (4.4.2), (4.4.4), (4.4.7) и теоремы 4.2.10. Равенство (4.4.3) является следствием неравенства (4.4.8) и определения [4]-2.1.17. □

Следствие 4.4.3. *Если* $\|f_n(x)\| \leq M$ *и* $f_n \to f$, *то*

$$\int_X d\mu(x)f(x) = \lim_{n \to \infty} \int_X d\mu(x)f_n(x)$$

□

Теорема 4.4.4 (Беппо Леви). *Пусть отображения*[4.13]

(4.4.9) $$f_n : A \to R \quad i = 1, ...$$

*интегрируемы и*

(4.4.10) $$\int_A d\mu(x)f_n(x) \leq M \quad n = 1, ...$$

*для некоторой константы* $M$. *Пусть для любого* $n$

(4.4.11) $$f_n(x) \leq f_{n+1}(x)$$

*Тогда почти всюду на* $A$ *существует конечный предел*

(4.4.12) $$f(x) = \lim_{n \to \infty} f_n(x)$$

*Отображение* $f$ *интегрируемо на* $A$ *и*

(4.4.13) $$\int_A d\mu(x)f(x) = \lim_{i \to \infty} \int_A d\mu(x)f_i(x)$$

Доказательство. Будем предполагать $f_1(x) \geq 0$, так как в общем случае мы можем положить

$$f'_i(x) = f_i(x) - f_1(x)$$

Согласно этому предположению и условию (4.4.11)

(4.4.14) $$f_n(x) \geq 0$$

Из утверждения (4.4.14) и определений 4.1.2, 4.2.1 следует, что для любого измеримого множества $B \subseteq A$

(4.4.15) $$\int_B d\mu(x)f(x) \geq 0$$

Из теоремы 4.3.2 следует, что для любого измеримого множества $B \subseteq A$

(4.4.16) $$\int_A d\mu(x)f(x) = \int_B d\mu(x)f(x) + \int_{A \setminus B} d\mu(x)f(x)$$

---

[4.13]Смотри также теорему [1]-7 на странице 303.

Утверждение

$$(4.4.17) \qquad \int_B d\mu(x) f_n(x) \leq M \quad n = 1, ...$$

следует из (4.4.10), (4.4.15), (4.4.16).

Рассмотрим множество

$$\Omega = \{x \in A : f_n(x) \to \infty\}$$

Тогда

$$\Omega = \bigcap_r \bigcup_n \Omega_n(r)$$

где

$$(4.4.18) \qquad \Omega_n(r) = \{x \in A : f_n(x) > r\}$$

Согласно неравенству Чебышева (4.3.35),

$$(4.4.19) \qquad \mu(\Omega_n(r)) \leq \frac{M}{r}$$

следует из (4.4.10), (4.4.18). Так как

$$\Omega_1(r) \subseteq \Omega_2(r) \subseteq ... \subseteq \Omega_n(r) \subseteq ...$$

то

$$(4.4.20) \qquad \mu\left(\bigcup_n \Omega_n(r)\right) \leq \frac{M}{r}$$

следует из (4.4.19). Так как при любом $r$

$$\Omega \subseteq \bigcup_n \Omega_n(r)$$

то

$$(4.4.21) \qquad \mu(\Omega) \leq \frac{M}{r}$$

Так как $r$ произвольно, то $\mu(\Omega) = 0$ является следствием (4.4.21). Следовательно, монотонная последовательность $f_n(x)$ почти всюду на $A$ имеет конечный предел $f(x)$.

Пусть

$$(4.4.22) \qquad A(r) = \{x \in A : r - 1 \leq f(x) < r\}$$

$$(4.4.23) \qquad \phi : A \to R \quad x \in A(r) => \phi(x) = r$$

Положим

$$(4.4.24) \qquad B(s) = \bigcup_{r=1}^{s} A(r)$$

Из (4.4.22), (4.4.24) следует, что отображения $f_n$ и $f$ ограничены на $B(s)$. Из (4.4.22), (4.4.23) следует, что

$$(4.4.25) \qquad \phi(x) \leq f(x) + 1$$

Согласно теоремам 2.2.6, 4.4.1,

$$(4.4.26) \qquad \int_{B_s} d\mu(x) \phi(x) \leq \int_{B_s} d\mu(x) f(x) + \mu(A)$$

следует из (4.4.25). Согласно следствию 4.4.3,

$$(4.4.27) \qquad \int_{B_s} d\mu(x)\phi(x) \leq \lim_{n\to\infty} \int_{B_s} d\mu(x)f_n(x) + \mu(A) \leq M + \mu(A)$$

следует из (4.4.17), (4.4.22), (4.4.24), (4.4.26). Согласно равенствам (4.4.22), (4.4.23) и определению 3.2.1, отображение $\phi$ является простым. Согласно определению 4.1.2,

$$(4.4.28) \qquad \int_{B_s} d\mu(x)\phi(x) = \sum_{r=1}^{s} r\mu(A(r))$$

Сходимость ряда

$$\sum_{r=1}^{\infty} r\mu(A(r)) \leq M + \mu(A) < \infty$$

следует из (4.4.27), (4.4.28). Согласно определению 3.2.1, отображение $\phi$ является интегрируемым

$$(4.4.29) \qquad \int_A d\mu(x)\phi(x) = \sum_{r=1}^{\infty} r\mu(A(r))$$

Следовательно, утверждение теоремы является следствием теоремы 4.4.2, так как

$$f_n(x) \leq f(x) \leq \phi(x)$$

следует из (4.4.11), (4.4.13), (4.4.22), (4.4.23). $\qquad \square$

Глава 5

# Теорема Фубини

## 5.1. Произведение полуколец множеств

Определение 5.1.1. Пусть $\mathcal{L}_i$, $i = 1, ..., n$, - система подмножеств множества $A_i$. **Декартово произведение систем подмножеств**

$$\mathcal{M} = \mathcal{L}_1 \times ... \times \mathcal{L}_n$$

является системой подмножеств множества $A = A_1 \times ... \times A_n$, представимых в виде[5.1]

$$C = C_1 \times ... \times C_n \quad C_i \in \mathcal{L}_i$$

Если $\mathcal{L}_1 = ... = \mathcal{L}_n = \mathcal{L}$, то система множеств

$$\mathcal{M} = \mathcal{L}^n$$

называется **декартова степень $n$ систем подмножеств** □

Теорема 5.1.2. *Декартово произведение*[5.2]

$$\mathcal{S} = \mathcal{S}_1 \times ...\mathcal{S}_n$$

*полуколец множеств $\mathcal{S}_1, ..., \mathcal{S}_n$ является полукольцом множеств.*

Доказательство. Мы докажем теорему для случая $n = 2$. В общем случае доказательство аналогично.

5.1.2.1: Пусть $A, B \in \mathcal{S}_1 \times \mathcal{S}_2$. Следовательно

(5.1.1) $$A = A_1 \times A_2 \quad A_1 \in \mathcal{S}_1 \quad A_2 \in \mathcal{S}_2$$

(5.1.2) $$B = B_1 \times B_2 \quad B_1 \in \mathcal{S}_1 \quad B_2 \in \mathcal{S}_2$$

Очевидно, что

$$(x_1, x_2) \in (A_1 \times A_2) \cap (B_1 \times B_2)$$

тогда и только тогда, когда

$$x_1 \in A_1 \cap B_1 \quad x_2 \in A_2 \cap B_2$$

Следовательно

(5.1.3) $$(A_1 \times A_2) \cap (B_1 \times B_2) = (A_1 \cap B_1) \times (A_2 \cap B_2)$$

Из (5.1.1), (5.1.2) и утверждения 2.1.1.2 следует, что

(5.1.4) $$A_1 \cap B_1 \in \mathcal{S}_1 \quad A_2 \cap B_2 \in \mathcal{S}_2$$

---

[5.1]Смотри также определение в [1] на странице 310.
[5.2]Смотри также теорему [1]-1 на странице 311.

Из (5.1.1), (5.1.2), (5.1.4) следует, что

$$(A_1 \times A_2) \cap (B_1 \times B_2) \in \mathcal{S}_1 \times \mathcal{S}_2$$

Следовательно, утверждение 2.1.1.2 верно для $\mathcal{S}_1 \times \mathcal{S}_2$.

5.1.2.2: Пусть $B_1 \subset A_1$, $B_2 \subset A_2$. Из утверждения 2.1.1.3 следует, что

(5.1.5) $\qquad A_1 = B_1 \cup B_{1 \cdot 1} \cup ... \cup B_{1 \cdot k_1} \quad B_{1 \cdot i} \in \mathcal{S}_1$

(5.1.6) $\qquad A_2 = B_2 \cup B_{2 \cdot 1} \cup ... \cup B_{2 \cdot k_2} \quad B_{2 \cdot i} \in \mathcal{S}_2$

Из (5.1.5), (5.1.6) следует, что

$$(x_1, x_2) \in A = A_1 \times A_2$$

тогда и только тогда, когда $(x_1, x_2)$ принадлежит одному из следующих множеств: $B = B_1 \times B_2$, $B_1 \times B_{2 \cdot j} \in \mathcal{S}_1 \times \mathcal{S}_2$, $B_{1 \cdot i} \times B_2 \in \mathcal{S}_1 \times \mathcal{S}_2$, $B_{1 \cdot i} \times B_{2 \cdot j} \in \mathcal{S}_1 \times \mathcal{S}_2$. Следовательно,

$$A = B \cup (B_1 \times B_{2 \cdot 1}) \cup ... \cup (B_1 \times B_{2 \cdot k_2})$$
$$\cup (B_{1 \cdot 1} \times B_2) \cup (B_{1 \cdot 1} \times B_{2 \cdot 1}) \cup ... \cup (B_{1 \cdot 1} \times B_{2 \cdot k_2})$$
$$...$$
$$\cup (B_{k_1 \cdot 1} \times B_2) \cup (B_{k_1 \cdot 1} \times B_{2 \cdot 1}) \cup ... \cup (B_{k_1 \cdot 1} \times B_{2 \cdot k_2})$$

и утверждение 2.1.1.3 верно для $\mathcal{S}_1 \times \mathcal{S}_2$.

Из рассуждений 5.1.2.1, 5.1.2.2 следует, что $\mathcal{S}_1 \times \mathcal{S}_2$ является полукольцом множеств. $\square$

Однако декартово произведение колец множеств, вообще говоря, не является кольцом множеств.

Определение 5.1.3. Кольцо множеств

$$\mathcal{R} = \mathcal{R}_1 \otimes ... \otimes \mathcal{R}_n$$

порождённое полукольцом множеств

$$\mathcal{R}_1 \times ... \times \mathcal{R}_n$$

называется **произведением колец множеств** $\mathcal{R}_1, ..., \mathcal{R}_n$. $\square$

Теорема 5.1.4. *Произведение алгебр множеств является алгеброй множеств.*

Доказательство. Для $i = 1, ..., n$, пусть кольцо множеств $\mathcal{R}_i$ является алгеброй множеств. Согласно определению 2.1.2, алгебра множеств $\mathcal{R}_i$ имеет единицу $E_i$ такую, что

(5.1.7) $\qquad A_i \cap E_i = A_i$

для любого $A_i \in \mathcal{R}_i$. Из рассуждения 5.1.2.1 и из (5.1.7) следует, что

(5.1.8) $\qquad (A_1 \times ... \times A_n) \cap (E_1 \times ... \times E_n) = A_1 \times ... \times A_n$

Из теоремы 2.1.11 и равенства (5.1.8) следует, что для любого $A \in \mathcal{R}_1 \otimes ... \otimes R_n$ верно равенство

$$A \cap (E_1 \times ... \times E_n) = A$$

Следовательно, множество $E_1 \times ... \times E_n$ является единицей кольца множеств $\mathcal{R}_1 \times ... \times \mathcal{R}_n$. Согласно определению 2.1.2, кольцо множеств $\mathcal{R}_1 \times ... \times \mathcal{R}_n$ является алгеброй множеств. $\square$

## 5.2. Произведение мер

Теорема 5.2.1. *Пусть*

$$\mu_i : \mathcal{R}_i \to R$$

*мера на полукольце* $\mathcal{R}_i$, $i = 1, ..., n$. **Декартово произведение**

$$\mu = \mu_1 \times ... \times \mu_n$$

**мер** *определено формулой*[5.3]

$$(5.2.1) \qquad \mu(A_1 \times ... \times A_n) = \mu_1(A_1)...\mu_n(A_n)$$

Доказательство. Мы докажем теорему для случая $n = 2$. В общем случае доказательство аналогично. Пусть

$$(5.2.2) \qquad A = A_1 \times A_2 = \bigcup_{i=1}^{t} B^i$$

где $i \neq j => B^i \cap B^j = \emptyset$ и $B^i = B_1^i \times B_2^i$. Согласно лемме 2.1.10, существуют разложения

$$A_1 = \bigcup_{m=1}^{r} C_1^m \quad A_2 = \bigcup_{n=1}^{s} C_2^n$$

такие, что

$$(5.2.3) \qquad B_i^k = \bigcup_{m \in M_{ik}} C_i^m \quad i = 1, 2$$

$$M_{1k} \subseteq \{1, ..., r\} \quad M_{2k} \subseteq \{1, ..., s\}$$

Согласно утверждению 2.2.1.2 по отношению к мерам $\mu_1, \mu_2$, из (5.2.3) следует, что

$$
\begin{aligned}
\mu(A) &= \mu_1(A_1)\mu_2(A_2) = \sum_{m=1}^{r} \mu(C_1^m) \sum_{n=1}^{s} \mu(C_2^n) \\
&= \sum_{k=1}^{t} \sum_{m \in M_{1k}} C_1^m \sum_{n \in M_{2k}} C_2^n \\
&= \sum_{k=1}^{t} \mu_1(B_1^k)\mu_2(B_2^k) = \sum_{k=1}^{t} \mu_1(B^k)
\end{aligned}
\qquad (5.2.4)
$$

Из (5.2.2), (5.2.4) следует, что утверждение 2.2.1.2 верно для $\mu$. $\square$

Теорема 5.2.2. *Пусть* $\mu_i$, $i = 1, ..., n$, - *$\sigma$-аддитивная мера,*[5.4] *определённая на $\sigma$-алгебре* $\mathcal{C}_i$. *Декартово произведение мер* $\mu_1 \times ... \times \mu_n$ *является $\sigma$-аддитивной мерой на полукольце* $\mathcal{C}_1 \times ... \times \mathcal{C}_n$.

Доказательство. Мы докажем теорему для случая $n = 2$. В общем случае доказательство аналогично. Пусть

$$(5.2.5) \qquad C = \bigcup_{n=1}^{\infty} C_n \quad C, C_n \in \mathcal{C}_1 \times \mathcal{C}_2$$

---

[5.3]Смотри [1], страница 312, доказательство, что отображение, определённое формулой (5.2.1), является мерой.

[5.4]Смотри также теорему [1]-2 на странице 313.

где

(5.2.6) $$n \neq m => C_n \cap C_m = \emptyset$$

Согласно (5.2.5)

$$C = A \times B \quad\quad A \in \mathcal{C}_1 \quad\quad B \in \mathcal{C}_2$$

$$C_n = A_n \times B_n \quad A_n \in \mathcal{C}_1 \quad B_n \in \mathcal{C}_2$$

Рассмотрим множество отображений

$$f_n : X \to R$$

определённых правилом

(5.2.7) $$f_n(x) = \begin{cases} \mu_2(B_n) & x \in A_n \\ 0 & x \notin A_n \end{cases}$$

Пусть $x \in A$. Пусть

$$N_x = \{n : \exists y \in B, (x,y) \in A_n \times B_n\}$$

Согласно (5.2.6), для любого $y \in B$ существует единственное $n \in N_x$ такое, что

$$(x,y) \in A_n \times B_n \quad y \in B_n$$

При этом

$$n \in N_x, m \in N_x, n \neq m => B_n \cap B_m = \emptyset$$
$$\bigcup_{n \in N_x} B_n = B$$

Согласно утверждению 2.2.10.2 относительно меры $\mu_2$

(5.2.8) $$\sum_n f_n(x) = \mu_2(B)$$

Согласно теоремам 4.2.3, 4.2.7 и следствию 4.4.3, равенство

(5.2.9) $$\sum_n \int_A d\mu_1(x) f_n(x) = \int_A d\mu_1(x) \mu_2(B) = \mu_1(A)\mu_2(B) = \mu(C)$$

является следствием равенства (5.2.8). Согласно теореме 4.2.7, равенство

(5.2.10) $$\int_A d\mu_1(x) f_n(x) = \mu_1(A_n)\mu_2(B_n) = \mu(C_n)$$

является следствием равенства (5.2.7). Равенство

$$\mu(C) = \sum_n \mu(C_n)$$

является следствием равенств (5.2.9), (5.2.10). Следовательно, утверждение 2.2.10.2 верно для меры $\mu$.    $\square$

Лебегово продолжение меры $\mu_1 \times ... \times \mu_n$ определённое на $\sigma$-алгебре $\mathcal{C} \supset \mathcal{C}_1 \times ... \times \mathcal{C}_n$, называется **произведением мер**

$$\mu_1 \otimes ... \otimes \mu_n = \bigotimes \mu_i$$

Если $\mu_1 = ... = \mu_n = \mu$, то произведение мер называется **степенью меры** $\mu$

$$\mu^n = \mu_1 \otimes ... \otimes \mu_n$$

Теорема 5.2.3. *Пусть* $\mu_i$, $i = 1, 2$, *- $\sigma$-аддитивная мера,[5.5] определённая на $\sigma$-алгебре $\mathcal{C}_i$ подмножеств множества $X_i$. Пусть* $\mu = \mu_1 \otimes \mu_2$. *Пусть* $A \in \mathcal{C}_1 \otimes \mathcal{C}_2$.

*Пусть для каждого* $x_1 \in X_1$

$$(5.2.11) \qquad A_{2x_1} = \{x_2 \in X_2 : (x_1, x_2) \in A\} \in \mathcal{C}_2$$

*Если отображение*

$$x_1 \to \mu_2(A_{2x_1})$$

*$\mu_1$-интегрируемо, то*

$$(5.2.12) \qquad \mu(A) = \int_{X_1} d\mu_1(x_1) \mu_2(A_{2x_1})$$

*Пусть для каждого* $x_2 \in X_2$

$$A_{1x_2} = \{x_1 \in X_1 : (x_1, x_2) \in A\} \in \mathcal{C}_1$$

*Если отображение*

$$x_2 \to \mu_1(A_{1x_2})$$

*$\mu_2$-интегрируемо, то*

$$(5.2.13) \qquad \mu(A) = \int_{X_2} d\mu_2(x_2) \mu_1(A_{1x_2})$$

Доказательство. Мы докажем равенство (5.2.12). Доказательство равенства (5.2.13) аналогично.

Лемма 5.2.4. *Равенство (5.2.12) верно для множества вида*

$$(5.2.14) \qquad A = A_1 \times A_2 \quad A_1 \in \mathcal{C}_1 \quad A_2 \in \mathcal{C}_2$$

Доказательство. Согласно определению 5.2.1,

$$(5.2.15) \qquad \mu(A) = \mu_1(A_1) \mu_2(A_2)$$

Согласно (5.2.11)

$$(5.2.16) \qquad A_{2x_1} = \begin{cases} A_2 & x_1 \in A_1 \\ \emptyset & x_1 \notin A_1 \end{cases}$$

Из (5.2.16) и теоремы 2.2.2 следует, что

$$(5.2.17) \qquad \mu_2(A_{2x_1}) = \begin{cases} \mu_2(A_2) & x_1 \in A_1 \\ 0 & x_1 \notin A_1 \end{cases}$$

Из (5.2.14), (5.2.17) и определения 3.2.1 следует, что отображение

$$x_1 \to \mu_2(A_{2x_1})$$

является простым. Согласно определению 4.1.2,

$$(5.2.18) \qquad \int_{X_1} d\mu_1(x_1) \mu_2(A_{2x_1}) = \mu_1(A_1) \mu_2(A_2)$$

Лемма следует из равенств (5.2.15), (5.2.18). $\odot$

---

[5.5]Смотри также теорему [1]-3 на странице 314. Я не предполагаю, что $X_i$ - множество действительных чисел.

Лемма 5.2.5. Равенство (5.2.12) верно для множества вида

$$(5.2.19) \qquad A = \bigcup_{i=1}^{n} A_{1i} \times A_{2i} \quad A_{1i} \in \mathcal{C}_1 \quad A_{2i} \in \mathcal{C}_2$$

Доказательство. Согласно лемме 2.1.10, существует конечная система множеств $B_{11}, ..., B_{1t_1} \in \mathcal{C}_1$ такая, что

$$(5.2.20) \qquad i \neq j \implies B_{1i} \cap B_{1j} = \emptyset$$

$$(5.2.21) \qquad A_{1i} = \bigcup_{j \in M_{1i}} B_{1j}$$

где $M_{1i} \subset \{1, ..., t_1\}$. Согласно лемме 2.1.10, существует конечная система множеств $B_{21}, ..., B_{2t_2} \in \mathcal{C}_2$ такая, что

$$(5.2.22) \qquad i \neq j \implies B_{2i} \cap B_{2j} = \emptyset$$

$$(5.2.23) \qquad A_{2i} = \bigcup_{j \in M_{2i}} B_{2j}$$

где $M_{2i} \subset \{1, ..., t_2\}$. Из (5.2.20), (5.2.22) следует, что

$$(5.2.24) \qquad i \neq j \land k \neq l \implies (B_{1i} \times B_{2k}) \cap (B_{1j} \times B_{2l}) = \emptyset$$

Из равенств (5.2.21), (5.2.23) следует, что для любого $i$

$$(5.2.25) \qquad A_{1i} \times A_{2i} = \bigcup_{k \in M_{1i}} \bigcup_{l \in M_{2i}} B_{1k} \times B_{2l}$$

Пусть

$$(5.2.26) \qquad M_i = M_{1i} \times M_{2i}$$

Равенство

$$(5.2.27) \qquad A_{1i} \times A_{2i} = \bigcup_{(k,l) \in M_i} B_{1k} \times B_{2l}$$

следует из (5.2.25), (5.2.26). Пусть

$$(5.2.28) \qquad M = \bigcup_i M_i$$

Равенство

$$(5.2.29) \qquad A = \bigcup_{(k,l) \in M} B_{1k} \times B_{2l}$$

следует из (5.2.19), (5.2.27), (5.2.28). Согласно утверждению 2.2.1.2 и определению 5.2.1,

$$(5.2.30) \qquad \mu(A) = \sum_{(k,l) \in M} \mu(B_{1k} \times B_{2l}) = \sum_{(k,l) \in M} \mu_1(B_{1k})\mu_2(B_{2l})$$

следует из (5.2.24), (5.2.29).

Согласно (5.2.11), (5.2.29)

$$(5.2.31) \qquad A_{2x_1} = \begin{cases} \displaystyle\bigcup_{(k,l) \in M} B_{2l} & x_1 \in B_{1k} \\ \emptyset & x_1 \notin \displaystyle\bigcup_k B_{1k} \end{cases}$$

Из равенства (5.2.31), утверждения 2.2.1.2, и теоремы 2.2.2 следует, что

$$(5.2.32) \quad \mu_2(A_{2x_1}) = \begin{cases} \sum_{(k,l) \in M} \mu_2(B_{2l}) & x_1 \in B_{1k} \\ 0 & x_1 \notin \bigcup_k B_{1k} \end{cases}$$

Из (5.2.29), (5.2.32) и определения 3.2.1 следует, что отображение

$$x_1 \to \mu_2(A_{2x_1})$$

является простым. Согласно определению 4.1.2,

$$(5.2.33) \quad \int_{X_1} d\mu_1(x_1)\mu_2(A_{2x_1}) = \sum_{(k,l) \in M} \mu_1(B_{1k})\mu_2(B_{2l})$$

Лемма следует из равенств (5.2.30), (5.2.33). $\odot$

Согласно теореме 2.3.13, для любого множества $A \in \mathcal{C}_{\mu_1 \otimes \mu_2}$ существует множество $B$ такое, что

$$(5.2.34) \quad A \subseteq B$$

$$(5.2.35) \quad \mu(A) = \mu(B)$$

$$(5.2.36) \quad B = \bigcap_n B_n$$

$$(5.2.37) \quad B_1 \supseteq B_2 \supseteq ... \supseteq B_n \supseteq ...$$

$$(5.2.38) \quad B_n = \bigcup_k B_{nk}$$

$$(5.2.39) \quad B_{nk} \in \mathcal{R}(\mathcal{C}_{\mu_1} \times \mathcal{C}_{\mu_2})$$

$$(5.2.40) \quad \mu(B_{nk}) < \mu(A) + \frac{1}{n}$$

$$(5.2.41) \quad B_{n1} \subseteq B_{n2} \subseteq ... \subseteq B_{nk} \subseteq ...$$

Согласно леммам 5.2.4, 5.2.5, из утверждения (5.2.39) следует, что равенство (5.2.12) верно для множеств $B_{nk}$

$$(5.2.42) \quad \mu(B_{nk}) = \int_{X_1} d\mu_1(x_1)\mu_2(B_{nk \cdot 2x_1})$$

где

$$(5.2.43) \quad B_{nk \cdot 2x_1} = \{x_2 \in X_2 : (x_1, x_2) \in B_{nk}\}$$

Утверждение

$$(5.2.44) \quad B_{n1 \cdot 2x_1} \subseteq B_{n2 \cdot 2x_1} \subseteq ... \subseteq B_{nk \cdot 2x_1} \subseteq ...$$

является следствием (5.2.41), (5.2.43). Рассмотрим множество отображений

$$f_{nk} : X_1 \to R \quad f_{nk}(x_1) = \mu_2(B_{nk \cdot 2x_1})$$

Согласно теореме 2.2.6, утверждение

$$(5.2.45) \quad f_{nk}(x_1) \leq f_{nk+1}(x_1)$$

является следствием утверждения (5.2.44). Утверждение

$$(5.2.46) \quad \int_{X_1} d\mu_1(x_1) f_{nk}(x_1) = \int_{X_1} d\mu_1(x_1) \mu_2(B_{nk \cdot 2 x_1}) < \mu(A) + \frac{1}{n}$$

является следствием утверждений (5.2.40), (5.2.42). Согласно теореме 4.4.4,

- почти всюду на $X_1$ существует конечный предел

$$(5.2.47) \quad f_n(x_1) = \lim_{k \to \infty} f_{nk}(x_1) = \lim_{k \to \infty} \mu_2(B_{nk \cdot 2 x_1})$$

- отображение $f_n$ интегрируемо на $X_1$ и

$$(5.2.48) \quad \begin{aligned} \int_{X_1} d\mu_1(x_1) f_n(x_1) &= \lim_{k \to \infty} \int_{X_1} d\mu_1(x_1) f_{nk}(x_1) \\ &= \lim_{k \to \infty} \int_{X_1} d\mu_1(x_1) \mu_2(B_{nk \cdot 2 x_1}) \end{aligned}$$

Утверждение

$$(5.2.49) \quad \lim_{k \to \infty} \int_{X_1} d\mu_1(x_1) \mu_2(B_{nk \cdot 2 x_1}) = \int_{X_1} d\mu_1(x_1) \lim_{k \to \infty} \mu_2(B_{nk \cdot 2 x_1})$$

является следствием утверждений (5.2.47), (5.2.48). Утверждение

$$(5.2.50) \quad B_{n \cdot 2 x_1} = \bigcup_k B_{nk \cdot 2 x_1}$$

где

$$B_{n \cdot 2 x_1} = \{x_2 \in X_2 : (x_1, x_2) \in B_n\}$$

следует из цепочки утверждений

$$x_2 \in B_{n \cdot 2 x_1} <=> (x_1, x_2) \in B_n <=> \exists k, (x_1, x_2) \in B_{nk}$$
$$<=> \exists k, x_2 \in B_{nk \cdot 2 x_1} <=> x_2 \in \bigcup_k B_{nk \cdot 2 x_1}$$

Согласно теореме 2.2.13,

- равенства

$$(5.2.51) \quad \lim_{k \to \infty} \mu_2(B_{nk \cdot 2 x_1}) = \mu_2(B_{n \cdot 2 x_1})$$

$$(5.2.52) \quad \lim_{k \to \infty} \int_{X_1} d\mu_1(x_1) \mu_2(B_{nk \cdot 2 x_1}) = \int_{X_1} d\mu_1(x_1) \mu_2(B_{n \cdot 2 x_1})$$

  являются следствием утверждений (5.2.44), (5.2.49), (5.2.50).
- равенство

$$(5.2.53) \quad \lim_{k \to \infty} \mu(B_{nk}) = \mu(B_n)$$

  является следствием утверждений (5.2.38), (5.2.41).

Из равенств (5.2.42), (5.2.52), (5.2.53) следует, что равенство (5.2.12) верно для множеств $B_n$

$$(5.2.54) \quad \mu(B_n) = \int_{X_1} d\mu_1(x_1) \mu_2(B_{n \cdot 2 x_1})$$

Утверждение

$$(5.2.55) \quad B_{1 \cdot 2 x_1} \supseteq B_{2 \cdot 2 x_1} \supseteq ... \supseteq B_{k \cdot 2 x_1} \supseteq ...$$

является следствием (5.2.37), (5.2.43). Равенство
$$f_n : X_1 \to R \quad f_n(x_1) = \mu_2(B_{n \cdot 2 x_1})$$
следует из равенств (5.2.47), (5.2.51). Согласно теореме 2.2.6, утверждение
$$(5.2.56) \quad f_n(x_1) \geq f_{n+1}(x_1)$$
является следствием утверждения (5.2.55). Согласно теореме 4.2.10, утверждение
$$(5.2.57) \quad \int_{X_1} d\mu_1(x_1) f_n(x_1) = \int_{X_1} d\mu_1(x_1) \mu_2(B_{n \cdot 2 x_1}) > 0$$
является следствием утверждения 2.2.1.1. Согласно теореме 4.4.4,
- почти всюду на $X_1$ существует конечный предел
$$(5.2.58) \quad f(x_1) = \lim_{n \to \infty} f_n(x_1) = \lim_{n \to \infty} \mu_2(B_{n \cdot 2 x_1})$$
- отображение $f$ интегрируемо на $X_1$ и
$$(5.2.59) \quad \begin{aligned} \int_{X_1} d\mu_1(x_1) f(x_1) &= \lim_{n \to \infty} \int_{X_1} d\mu_1(x_1) f_n(x_1) \\ &= \lim_{n \to \infty} \int_{X_1} d\mu_1(x_1) \mu_2(B_{n \cdot 2 x_1}) \end{aligned}$$

Утверждение
$$(5.2.60) \quad \lim_{n \to \infty} \int_{X_1} d\mu_1(x_1) \mu_2(B_{n \cdot 2 x_1}) = \int_{X_1} d\mu_1(x_1) \lim_{n \to \infty} \mu_2(B_{n \cdot 2 x_1})$$
является следствием утверждений (5.2.58), (5.2.59). Утверждение
$$(5.2.61) \quad B_{2 x_1} = \bigcap_k B_{n \cdot 2 x_1}$$
где
$$B_{2 x_1} = \{x_2 \in X_2 : (x_1, x_2) \in B\}$$
следует из цепочки утверждений
$$x_2 \in B_{2 x_1} <=> (x_1, x_2) \in B <=> \forall n, (x_1, x_2) \in B_n$$
$$<=> \forall n, x_2 \in B_{n \cdot 2 x_1} <=> x_2 \in \bigcap_n B_{n \cdot 2 x_1}$$

Согласно теореме 2.2.12,
- равенство
$$(5.2.62) \quad \lim_{n \to \infty} \int_{X_1} d\mu_1(x_1) \mu_2(B_{n \cdot 2 x_1}) = \int_{X_1} d\mu_1(x_1) \mu_2(B_{2 x_1})$$
является следствием утверждений (5.2.55), (5.2.60), (5.2.61).
- равенство
$$(5.2.63) \quad \lim_{n \to \infty} \mu(B_n) = \mu(B)$$
является следствием утверждений (5.2.36), (5.2.37).

Из равенств (5.2.54), (5.2.62), (5.2.63) следует, что равенство (5.2.12) верно для множества $B$
$$(5.2.64) \quad \mu(B) = \int_{X_1} d\mu_1(x_1) \mu_2(B_{2 x_1})$$

Лемма 5.2.6. Если

$$\mu(A) = 0 \tag{5.2.65}$$

то

$$\int_{X_1} d\mu_1(x_1)\mu_2(A_{2x_1}) = 0 \tag{5.2.66}$$

Доказательство. Равенство

$$\int_{X_1} d\mu_1(x_1)\mu_2(B_{2x_1}) = 0 \tag{5.2.67}$$

следует из равенств (5.2.35), (5.2.65). Согласно теореме 4.3.6 и утверждению 2.2.1.1, равенство

$$\mu_2(B_{2x_1}) = 0 \tag{5.2.68}$$

почти всюду является следствием равенства (5.2.67). Утверждение

$$A_{2x_1} \subseteq B_{2x_1} \tag{5.2.69}$$

является следствием утверждения (5.2.34). Согласно теореме 2.2.8, равенство

$$\mu_2(A_{2x_1}) = 0 \tag{5.2.70}$$

почти всюду является следствием (5.2.68), (5.2.69). Согласно определению 4.1.2, равенство (5.2.66) является следствием равенства (5.2.70). ⊙

Пусть $C = B \setminus A$. Тогда

$$B = A \cup C \quad A \cap C = \emptyset \tag{5.2.71}$$

$$B_{2x_1} = A_{2x_1} \cup C_{2x_1} \quad A_{2x_1} \cap C_{2x_1} = \emptyset \tag{5.2.72}$$

Равенства

$$\mu(B) = \mu(A) + \mu(C) \tag{5.2.73}$$

$$\mu_2(B_{2x_1}) = \mu_2(A_{2x_1}) + \mu_2(C_{2x_1}) \tag{5.2.74}$$

являются следствием равенств (5.2.71), (5.2.72). Равенство

$$\mu(A) + \mu(C) = \int_{X_1} d\mu_1(x_1)\mu_2(A_{2x_1}) + \int_{X_1} d\mu_1(x_1)\mu_2(C_{2x_1}) \tag{5.2.75}$$

является следствием равенств (5.2.64), (5.2.73), (5.2.74) и теоремы 4.2.3. Равенство

$$\mu(C) = 0 \tag{5.2.76}$$

является следствием равенств (5.2.35), (5.2.73). Согласно лемме 5.2.6, равенство

$$\int_{X_1} d\mu_1(x_1)\mu_2(C_{2x_1}) = 0 \tag{5.2.77}$$

является следствием равенства (5.2.76). Равенство (5.2.12) является следствием равенств (5.2.75), (5.2.76), (5.2.77). □

## 5.3. Теорема Фубини

Теорема 5.3.1 (Фубини). *Пусть $\mu_1$ и $\mu_2$ - $\sigma$-аддитивные полные меры, определеные на $\sigma$-алгебрах.*[5.6] *Пусть $\mu = \mu_1 \otimes \mu_2$. Пусть $A$ - полная абелевая $\Omega$-группа. Пусть отображение*

$$f : X_1 \times X_2 \to A$$

*с компактным множеством значений $\mu$-интегрируемо на множестве $B \subseteq X_1 \times X_2$. Тогда*

$$(5.3.1) \qquad \int_B d\mu(x_1, x_2) f(x_1, x_2) = \int_{X_1} d\mu(x_1) \int_{B_{2x_1}} d\mu(x_2) f(x_1, x_2)$$

$$(5.3.2) \qquad \int_B d\mu(x_1, x_2) f(x_1, x_2) = \int_{X_2} d\mu(x_2) \int_{B_{1x_2}} d\mu(x_1) f(x_1, x_2)$$

Доказательство. Мы докажем равенство (5.3.1). Доказательство равенства (5.3.2) аналогично.

Лемма 5.3.2. *Равенство (5.3.1) верно для простого отображения $f$.*

Доказательство. Пусть $y_1, y_2, ...$ - область значений отображения $f$. Пусть

$$(5.3.3) \qquad B_i = \{(x_1, x_2) \in B : f(x_1, x_2) = y_i\}$$

Так как

$$(5.3.4) \qquad B = \bigcup_i B_i \quad i \neq j => B_i \cap B_j = \emptyset$$

то равенство

$$(5.3.5) \qquad \int_B d\mu(x_1, x_2) f(x_1, x_2) = \sum_n \mu(B_n) y_n$$

является следствием определения 4.1.2. Равенство

$$(5.3.6) \qquad \mu(B_n) = \int_{X_1} d\mu_1(x_1) \mu_2(B_{n \cdot 2 x_1})$$

где

$$(5.3.7) \qquad B_{n \cdot 2 x_1} = \{x_2 \in X_2 : (x_1, x_2) \in B_n\}$$

является следствием теоремы 5.2.3. Равенство

$$(5.3.8) \qquad \int_B d\mu(x_1, x_2) f(x_1, x_2) = \sum_n \left( \int_{X_1} d\mu_1(x_1) \mu_2(B_{n \cdot 2 x_1}) \right) y_n$$

является следствием равенств (5.3.5), (5.3.6). Согласно определению 4.1.2, равенство

$$(5.3.9) \qquad \int_{B_{2x_1}} d\mu_2(x_2) f(x_1, x_2) = \sum_n \mu_2(B_{n \cdot 2 x_1}) y_n$$

---

[5.6] Смотри также теорему [1]-5 на странице 317.

является следствием равенств (5.3.3), (5.3.7). Согласно теоремам 4.2.3, 4.2.8, равенство

$$\int_{X_1} d\mu_1(x_1) \int_{B_{2x_1}} d\mu(x_2) f(x_1, x_2) = \int_{X_1} d\mu_1(x_1) \left( \sum_n \mu_2(B_{n \cdot 2x_1}) y_n \right)$$
(5.3.10)
$$= \sum_n \int_{X_1} d\mu_1(x_1)(\mu_2(B_{n \cdot 2x_1}) y_n)$$
$$= \sum_n \left( \int_{X_1} d\mu_1(x_1) \mu_2(B_{n \cdot 2x_1}) \right) y_n$$

является следствием равенства (5.3.9). Равенство (5.3.1) является следствием равенств (5.3.8), (5.3.10). $\odot$

Согласно теореме 3.3.2, существует последовательность простых отображений

$$f_n : X_1 \times X_2 \to A$$

сходящаяся равномерно к $f$

(5.3.11) $$f(x) = \lim_{n \to \infty} f_n(x)$$

Согласно лемме 5.3.2, равенство

(5.3.12) $$\int_B d\mu(x_1, x_2) f_n(x_1, x_2) = \int_{X_1} d\mu(x_1) \int_{B_{2x_1}} d\mu(x_2) f_n(x_1, x_2)$$

верно для любого отображения $f_n$. Равенства

(5.3.13) $$\int_B d\mu(x_1, x_2) f(x_1, x_2) = \lim_{n \to \infty} \int_B d\mu(x_1, x_2) f_n(x_1, x_2)$$

(5.3.14) $$\int_{B_{2x_1}} d\mu_2(x_2) f(x_1, x_2) = \lim_{n \to \infty} \int_{B_{2x_1}} d\mu_2(x_2) f_n(x_1, x_2)$$

являются следствием определения 4.2.1. Так как отображение $f$ имеет компактное множество значений, то существует $M > 0$ такое, что

(5.3.15) $$\|f(x_1, x_2)\| \leq M$$

Согласно построению в доказательстве теоремы 3.3.2, утверждение

(5.3.16) $$\|f_n(x_1, x_2)\| < M + \frac{1}{n} \leq M + 1$$

являются следствием утверждения (5.3.15). Согласно теореме 4.2.10, утверждение

(5.3.17) $$\left| \int_{B_{2x_1}} d\mu_2(x_2) f_n(x_1, x_2) \right| \leq \mu_2(B_{2x_1})(M + 1)$$

являются следствием утверждения (5.3.16). Согласно теоремам 4.2.8, 5.2.3, отображение

$$x_1 \to \mu_2(B_{2x_1})(M + 1)$$

$\mu_1$-интегрируемо. Согласно теореме 4.4.2, из (5.3.11), (5.3.14) следует, что отображение

$$x_1 \to \int_{B_{2x_1}} d\mu_2(x_2) f(x_1, x_2)$$

$\mu_1$-интегрируемо и

$$
\begin{aligned}
&\int_{X_1} d\mu_1(x_1) \int_{B_{2x_1}} d\mu_2(x_2) f(x_1, x_2) \\
&= \int_{X_1} d\mu_1(x_1) \lim_{n \to \infty} \int_{B_{2x_1}} d\mu_2(x_2) f_n(x_1, x_2) \\
&= \lim_{n \to \infty} \int_{X_1} d\mu_1(x_1) \int_{B_{2x_1}} d\mu_2(x_2) f_n(x_1, x_2)
\end{aligned}
\tag{5.3.18}
$$

является следствием (5.3.14). Равенство (5.3.1) является следствием равенств (5.3.13), (5.3.18). □

Глава 6

# Список литературы

# Глава 7

# Предметный указатель

# Глава 8

# Специальные символы и обозначения